\documentclass{article}

\usepackage{arxiv}
\usepackage{graphicx, color}
\usepackage{amsmath,amssymb}
\usepackage[utf8]{inputenc} 
\usepackage[T1]{fontenc}    
\usepackage{hyperref}       
\usepackage{url}            
\usepackage{booktabs}       
\usepackage{amsfonts}       
\usepackage{microtype}      
\usepackage{lipsum}

\title{Damped and Divergence Exact Solutions for the Duffing Equation using Leaf Functions and Hyperbolic Leaf Functions}

\author{
  Kazunori Shinohara\thanks{10-3 Takiharu-cho, Minami-ku, Nagoya 457-8530, Japan} \\
  Department of Mechanical Systems Engineering\\
  Daido University\\
  10-3 Takiharu-cho, Minami-ku, Nagoya 457-8530, Japan \\
  \texttt{shinohara@06.alumni.u-tokyo.ac.jp} \\
}

\begin{document}
\maketitle

\begin{abstract}
According to the wave power rule, the second derivative of a function $x(t)$ with respect to the variable $t$ is equal to negative $n$ times the function $x(t)$ raised to the power of $2n-1$. Solving the ordinary differential equations numerically results in waves appearing in the figures. The ordinary differential equation is very simple; however, waves, including the regular amplitude and period, are drawn in the figure. In this study, the function for obtaining the wave is called the leaf function. Based on the leaf function, the exact solutions for the undamped and unforced Duffing equations are presented. In the ordinary differential equation, in the positive region of the variable $x(t)$, the second derivative $\mathrm{d}^2x(t)/\mathrm{d}t^2$ becomes negative. Therefore, in the case that the curves vary with the time $t$ under the condition $x(t)>0$, the gradient $\mathrm{d}x(t)/\mathrm{d}t$ constantly decreases as time $t$ increases. That is, the tangential vector on the curve of the graph (with the abscissa $t$ and the ordinate $x(t)$) changes from the upper right direction to the lower right direction as time $t$ increases. On the other hand, in the negative region of the variable $x(t)$, the second derivative $\mathrm{d}^2x(t)/\mathrm{d}t^2$ becomes positive. The gradient $\mathrm{d}x(t)/\mathrm{d}t$ constantly increases as time $t$ decreases. That is, the tangent vector on the curve changes from the lower right direction to the upper right direction as time $t$ increases. Since the behavior occurring in the positive region of the variable $x(t)$ and the behavior occurring in the negative region of the variable $x(t)$ alternately occur in regular intervals, waves appear by these interactions. 
In this paper, I present seven types of damped and divergence exact solutions by combining trigonometric functions, hyperbolic functions, hyperbolic leaf functions, leaf functions, and exponential functions. In each type, I show the derivation method and numerical examples, as well as describe the features of the waveform.
\end{abstract}


\section{Introduction}
\label{Introduction}
\subsection{Hyperbolic leaf functions}
\label{Hyperbolic Leaf Functions}
In our previous report, the exact solutions for the Duffing equation under free vibration were presented using leaf functions \cite{Kaz_duf}. Leaf functions comprise two types of functions with initial conditions. One of them satisfies the following ordinary differential equation and initial conditions.
\begin{equation}
\frac{\mathrm{d}^2x(t) }{\mathrm{d}t^2}=-nx(t)^{2n-1} \label{1.1.1}
\end{equation}
\begin{equation}
x(0)=0 \label{1.1.2}
\end{equation}
\begin{equation}
\frac{\mathrm{d}x(0)}{\mathrm{d}t}=1\label{1.1.3}
\end{equation}
The variable $n$ represents an integer. In the paper, the variable was named as the basis. A function that satisfies the above equations is defined as follows.
\begin{equation}
x(t)=\mathrm{sleaf}_{n}(t) \label{1.1.4}
\end{equation}
Another function satisfies the following ordinary differential equations and initial conditions.
\begin{equation}
\frac{\mathrm{d}^2x(t)}{\mathrm{d}t^2}=-nx(t)^{2n-1} \label{1.1.5}
\end{equation}
\begin{equation}
x(0)=1 \label{1.1.6}
\end{equation}
\begin{equation}
\frac{\mathrm{d}x(0)}{\mathrm{d}t}=0 \label{1.1.7}
\end{equation}
The variable $n$ represents an integer. A function that satisfies the above equations is defined as follows.
\begin{equation}
x(t)=\mathrm{cleaf}_{n}(t) \label{1.1.8}
\end{equation}
For the basis $n=1$, the function $\mathrm{sleaf}_{1}(t)$ and $\mathrm{cleaf}_{1}(t)$ represent the trigonometric functions $\mathrm{sin}(t)$ and $\mathrm{cos}(t)$, respectively. For the basis $n=2$, the functions $\mathrm{sleaf}_{2}(t)$ and $\mathrm{cleaf}_{2}(t)$ represent the lemniscate functions $\mathrm{sl}(t)$ and $\mathrm{cl}(t)$, respectively. Based on results obtained from solving the ordinary differential equation numerically, Eqs. (\ref{1.1.1})-(\ref{1.1.3}) or Eqs. (\ref{1.1.5})-(\ref{1.1.7}) yield waves with respect to an arbitrary $n$. The ordinary differential equation is very simple. However, waves including the regular amplitude and period are drawn in the figure. 
On the right side of Eq. (\ref{1.1.1}) and Eq. (\ref{1.1.5}), if the minus changes to plus, the waves disappear with respect to any arbitrary $n$. The variable $x(t)$ increases monotonically as the variable $t$ increases. These ordinary differential equations are as follows:
\begin{equation}
\frac{\mathrm{d}^2x(t)}{\mathrm{d}t^2}=nx(t)^{2n-1} \label{1.1.9}
\end{equation}
\begin{equation}
x(0)=0 \label{1.1.10}
\end{equation}
\begin{equation}
\frac{\mathrm{d}x(0)}{\mathrm{d}t}=1\label{1.1.11}
\end{equation}
In this study, a function that satisfies the above equations is defined as $\mathrm{sleafh}_{n}(x)$ \cite{Kaz_sh}. The basis $n$ represents an integer.
\begin{equation}
x(t)=\mathrm{sleafh}_{n}(t)
\end{equation}
For the basis $n=1$, the function $\mathrm{sleafh}_{1}(t)$ represents the hyperbolic function $sinh(t)$. In case of the basis $n=2$, the function $\mathrm{sleafh}_{2}(t)$ represents the Hyperbolic Lemniscate Function $\mathrm{slh}(t)$ or $\mathrm{sinhlemn}(t)$ \cite{Ber, Wei}. We further discuss combinations of the following differential equations and initial conditions.
\begin{equation}
\frac{\mathrm{d}^2x(t)}{\mathrm{d}t^2}=nx(t)^{2n-1} \label{1.1.13}
\end{equation}
\begin{equation}
x(0)=1 \label{1.1.14}
\end{equation}
\begin{equation}
\frac{\mathrm{d}x(0)}{\mathrm{d}t}=0 \label{1.1.15}
\end{equation}
The basis $n$ represents an integer. A function that satisfies the above equations is defined as follows \cite{Kaz_ch}.
\begin{equation}
x(t)=\mathrm{cleafh}_{n}(t)
\end{equation}
For the basis $n=1$, the function $\mathrm{cleafh}_{1}(t)$ represents the hyperbolic function $\mathrm{cosh}(t)$. For the basis $n=2$, the leaf function that satisfies Eq. (\ref{1.1.13}) - (\ref{1.1.15}) is $\mathrm{cleafh}_{2}(t)$. In the literature, the corresponding functions with respect to the function $\mathrm{cleafh}_{2}(t)$ cannot be found \cite{Wei}. For $n=3$, Ramanujan suggested an inverse function using a series \cite{Ber}. Apart from this type of solution, exact solutions have not been presented. By applying an imaginary number to the phase of a trigonometric function, the relation between the trigonometric function and the hyperbolic function can be obtained. Similar analogies exist between leaf functions and hyperbolic leaf functions. Using the imaginary number for the phase of this leaf function, the relation between the leaf function and the hyperbolic leaf function can be derived \cite{Kaz_sh}. First, we describe the relation between the leaf function $\mathrm{sleaf}_{n}(t)$ and the hyperbolic leaf function $\mathrm{sleafh}_{n}(t)$. When the basis $n(=2m-1)$ is an odd number, the following relation is obtained.
\begin{equation}
\mathrm{sleaf}_{2m-1}(i \cdot t)=i \cdot \mathrm{sleafh}_{2m-1}(t) \\ (m=1,2,3, \cdots)
\end{equation}
If the basis $n(=2m)$ is an even number, the following relation is obtained.
\begin{equation}
\mathrm{sleaf}_{2m}(i \cdot t)=i \cdot \mathrm{sleaf}_{2m}(t) \\ (m=1,2,3, \cdots)
\end{equation}
\begin{equation}
\mathrm{sleafh}_{2m}(i \cdot t)=i \cdot \mathrm{sleafh}_{2m}(t) \\ (m=1,2,3, \cdots)
\end{equation}
With respect to an arbitrary basis $n$, the relation between the leaf function $\mathrm{cleaf}_{n}(t)$ and the hyperbolic function $\mathrm{cleafh}_{n}(t)$ is given as follows:
\begin{equation}
\mathrm{cleafh}_{n}(i \cdot t)=\mathrm{cleaf}_{n}(t) \\ (n=1,2,3, \cdots)
\end{equation}
Using the above relations, we can derive the relation equation between the hyperbolic leaf function $\mathrm{sleafh}_n(t)$ and $\mathrm{cleafh}_n(t)$ using the relation equation between the leaf functions $\mathrm{sleaf}_n(t)$ and $\mathrm{cleaf}_n(t)$ at the basis $n=1, 2$ and $3$.  We can also derive the addition theorem for the hyperbolic leaf functions $\mathrm{sleafh}_n(t)$ and $\mathrm{cleafh}_n(t)$ using the addition theorem between $\mathrm{sleaf}_n(t)$ and $\mathrm{cleaf}_n(t)$ \cite{Kaz_sh, Kaz_add}. These relation equations and addition theorems are both theoretically and numerically consistent.

\subsection{Comparison of legacy functions with both leaf functions and hyperbolic leaf functions}
\label{Comparison}
The leaf functions and hyperbolic leaf functions based on the basis $n=1$ are as follows:
\begin{equation}
\mathrm{sleaf}_{1}(t)=\mathrm{sin}(t)
\end{equation}
\begin{equation}
\mathrm{cleaf}_{1}(t)=\mathrm{cos}(t)
\end{equation}
\begin{equation}
\mathrm{sleafh}_{1}(t)=\mathrm{sinh}(t)
\end{equation}
\begin{equation}
\mathrm{cleafh}_{1}(t)=\mathrm{cosh}(t)
\end{equation}
The lemniscate functions are presented by Gauss \cite{Gauss, roy}. The relation equations between this function and the leaf function are as follows:
\begin{equation}
\mathrm{sleaf}_{2}(t)=\mathrm{sl}(t)
\end{equation}
\begin{equation}
\mathrm{cleaf}_{2}(t)=\mathrm{cl}(t)
\end{equation}
\begin{equation}
\mathrm{sleafh}_{2}(t)=\mathrm{slh}(t) \label{slh}
\end{equation}
The definition of the function $\mathrm{slh}(t)$ of Eq.(\ref{slh}) can be confirmed by references \cite{Car, Neu}. For the function corresponding to the hyperbolic leaf function $\mathrm{cleafh}_2(t)$, no clear description can be found in the literature. In the case that the basis of a leaf function or a hyperbolic leaf function becomes $n=3$ or higher, it cannot be represented by the lemniscate function. 
We can obtain the second derivatives of the Jacobi elliptic function, the lemniscate function, the leaf function, or the hyperbolic leaf function. These functions are defined as $x(t)$. These second derivatives of the function $x(t)$ can be described by the sum of terms obtained by raising the original function, $x(t)$.
\begin{equation}
\frac{\mathrm{d}^2x(t)}{\mathrm{d}t^2}=c_0+c_1x(t)+c_2x(t)^2+c_3x(t)^3 (c_0 \sim c_3:constants)
\end{equation}
By deciding the values of the coefficients $c_0 \sim c_3$, one of either functions in the lemniscate function, the leaf function and the hyperbolic function is determined. The maximum value of the exponent in the above equation is 3. Historically, discussions have not been had in the literature with respect to exponents above 4. Using definitions based on the leaf function or the hyperbolic leaf function, we can discuss differential equations composed of higher order exponents such as the following
\begin{equation}
\begin{split}
\frac{\mathrm{d}^2x(t)}{\mathrm{d}t^2}=c_0+c_1x(t)+c_2x(t)^2+c_3x(t)^3+c_4x(t)^4 +c_5x(t)^5  \\
+ \cdots +c_{n-1}x(t)^{n-1}+c_nx(t)^n  (c_0 \sim c_n:constants)
\end{split}
\end{equation}

\subsection{The Duffing equation}
\label{Duffing equation}
For the basis $n=2$, both leaf functions and hyperbolic functions are applied to the damped or divergence Duffing equation. The exact solutions consist of seven types of damped and divergence exact solutions. In each case, the derivation method and numerical examples are shown and the features of the waveform described. The Duffing equation consists of the second derivative of the unknown function, the first derivative of the unknown function, the unknown function, the cube of the unknown function and a trigonometric function. The equation is presented as follows:
\begin{equation}
\frac{\mathrm{d}^2x(t)}{\mathrm{d}t^2} + \delta \frac{\mathrm{d}x(t)}{\mathrm{d}t}+ \alpha x(t)+ \beta x(t)^3=F\mathrm{cos}(\omega t) \label{Duffing}
\end{equation}
The variable $x(t)$ in the above equation represents an unknown function. The unknown function depends on time, $t$. The sign $\mathrm{d}x(t)/dt$ and the sign $\mathrm{d}^2x(t)/\mathrm{d}t^2$ represent the first and second derivatives of $x(t)$,respectively. The coefficients $\delta, \alpha, \beta$, $F$ and $\omega$ do not depend on time, $t$. If the above equation is regarded as a mathematical model of mechanical vibration, the first, second, third and fourth terms on the left side of Eq. (\ref{Duffing}) represent the inertial, damping, rigidity, and nonlinear stiffness terms, respectively. The term on the right side of Eq. (\ref{Duffing}) represents the external force term.

\subsection{Solving the Duffing equation}
\label{Solving the Duffing Equation}
The method for finding the solution of the damped Duffing equation is roughly divided into two types: numerical solutions and exact solutions. In numerical solutions, a simple technique based on Taylor expansion is applied to determine an approximate solution for a nonlinear Duffing oscillator with damping effect under different initial conditions \cite{Nad}. The nonlinear problem of the Duffing equation has been considered using the homotopy analysis technique. Compared to the perturbation methods \cite{Nay}, homotopy treatment does not require any small parameters \cite{Tur}.The dynamic behaviors and effects of random parameters are investigated. The Chebyshev orthogonal polynomial approximation method is applied to reduce the random parameter \cite{Zha}. The harmonic oscillations of freedom of a Duffing oscillator with large damping are investigated using a simple point collocation method \cite{Dai}.The improved constrained optimization harmonic balance method is employed to solve the Duffing oscillator. Analytical gradients of the object function are formulated. The sensitivity information of the Fourier coefficients can also be obtained \cite{Liao}. Conversely, exact solutions of the Duffing equation using the Jacobi elliptic function are discussed. The exact solution for the cubic-quintic Duffing oscillator based on the use of Jacobi elliptic functions is presented \cite{Alex}. Furthermore, the exact solutions for the damped Duffing equation are presented by extending the parameters of the exact solution \cite{Alex2}. The authors in Ref. \cite{Bel} do not assume any expression for the solution but exactly solved the nonlinear differential equation, unlike the Elias-Zuniga procedure. An analytical solution for the damped Duffing equation was derived by setting the time-dependent modulus in the Jacobi elliptic function \cite{Kim}.

\subsection{Originality and purpose}
\label{Originality and Purpose}
In our previous paper, exact solutions for undamped and unforced Duffing equations were presented by Ref. \cite{Kaz_duf} using the leaf function. Conversely, the divergence phenomena (or damped vibration) can be obtained in the Duffing equation. The purpose of this paper is to present more exact solutions for the Duffing equation by combining hyperbolic leaf functions, leaf functions, and exponential functions. These combining functions make it possible to produce divergence phenomena or damped vibration for the Duffing equation. To represent the exact solution that decay with time, the exponential function is placed in the phase of the leaf function. In this type of solution, the description of the exact solutions could not be found in literature. In contrast, in the exact solutions that diverge with time, the original functions in the exact solutions also need to diverge with time. However, there are no examples defined as solutions for the Duffing equation using the lemniscate function $\mathrm{slh}(t)$. In this paper, to represent the exact solution for the Duffing equation, the integral function of the leaf function or the hyperbolic leaf function is placed in the phase of the trigonometric functions or hyperbolic functions. The description of the exact solutions could also not be found in the literature. 
In this paper, seven types of the exact solution are presented, out of which five are divergence solutions and two are damped solutions without external forces. In an exact solution using the Jacobi elliptic function, multiple parameters of the Jacobi elliptic function cumulatively influence the periods and amplitudes. Therefore, the parameter based on a period is hard to separate from the parameter based on an amplitude. To determine these parameters, we need to solve a sixth order equation \cite{Alex2}. The exact solutions presented in this paper have the feature that the coefficients of the Duffing equation can be represented in both amplitude and phase by two parameters using a leaf function, a hyperbolic leaf function, and an exponential function.

\section{Numerical data for the hyperbolic leaf function}
\label{Numerical data for the Hyperbolic leaf function}
In the hyperbolic leaf function described in the previous section, curves and numerical data were described. The exact solutions are applied to the hyperbolic leaf functions or leaf functions based on the basis, $n = 2$. The graph of the hyperbolic leaf function based on the basis $n=2$ is shown in Fig. 1. The vertical and horizontal axes represent the variables $x(t)$ and $t$, respectively. For the basis $n=2$, the numerical data are summarized in Table 1. The following equations are obtained from the initial conditions in Eq. (\ref{1.1.10}) and Eq. (\ref{1.1.14}).
\begin{equation}
x(0)=\mathrm{sleafh}_2(0)=0 \label{2.1}
\end{equation}
\begin{equation}
x(0)=\mathrm{cleafh}_2(0)=1 \label{2.2}
\end{equation}
The hyperbolic leaf function $\mathrm{sleafh}_2(t)$ is an odd function whereas the hyperbolic leaf function $\mathrm{cleafh}_2(t)$ is an even function. The following relations are obtained.
\begin{equation}
\mathrm{sleafh}_2(-t)=-\mathrm{sleafh}_2(t) \label{2.3}
\end{equation}
\begin{equation}
\mathrm{cleafh}_2(-t)=\mathrm{cleafh}_2(t) \label{2.4}
\end{equation}
The hyperbolic leaf functions $\mathrm{sleafh}_2(t)$ and $\mathrm{cleafh}_2(t)$ have limits. Given that each limit is $\zeta_2$ and $\eta_2$, they are obtained by the following equations \cite{Kaz_sh, Kaz_ch}.
\begin{equation}
\zeta_2=\int_{0}^{\infty} \frac{1}{ \sqrt{1+u^4} } du=1.85407 \dotsb
\label{2.5}
\end{equation}
\begin{equation}
\eta_2=\int_{1}^{\infty} \frac{1}{ \sqrt{u^4-1} } du=1.31102 \dotsb
\label{2.6}
\end{equation}
From these extreme values, the following relations are obtained by the limits $\zeta_2$ and $\eta_2$.
\begin{equation}
\lim_{t \to \zeta_2} \mathrm{sleafh}_2(t)=+\infty \label{2.7}
\end{equation}
\begin{equation}
\lim_{t \to -\zeta_2} \mathrm{sleafh}_2(t)=-\infty \label{2.8}
\end{equation}\begin{equation}
\lim_{t \to \pm\eta_2} \mathrm{cleafh}_2(t)=+\infty \label{2.9}
\end{equation}
\begin{table}
\caption{ Numerical data for $\mathrm{sleafh}_2(t)$, $\mathrm{cleafh}_2(t)$, $\int_{0}^{t} \mathrm{sleafh}_2(u)\mathrm{d}u$ and $\int_{0}^{t} \mathrm{cleafh}_2(u)\mathrm{d}u$ }
\label{tab:1}  
\centering
\begin{tabular}{ccccc}
\hline\noalign{\smallskip}
$t$ & $\mathrm{sleafh}_2(t)$ & $\mathrm{cleafh}_2(t)$  & $\int_{0}^{t} \mathrm{sleafh}_2(u)\mathrm{d}u$ & $\int_{0}^{t} \mathrm{cleafh}_2(u)\mathrm{d}u$ \\
\noalign{\smallskip}\hline\noalign{\smallskip}
0.0	&	0.0	&	1.0	&	0.0	&	0.0	\\
0.1	&	0.100001 $\cdots$	&	1.010050302$\cdots$	&	0.005000017$\cdots$	&	0.100334338$\cdots$	\\
0.2	&	0.200032004 $\cdots$	&	1.040819659$\cdots$	&	0.020001067$\cdots$	&	0.202699225$\cdots$	\\
0.3	&	0.300243164 $\cdots$	&	1.094280815$\cdots$	&	0.045012155$\cdots$	&	0.309252773$\cdots$	\\
0.4	&	0.401026189$\cdots$	&	1.174155243$\cdots$	&	0.080068354$\cdots$	&	0.422433069$\cdots$	\\
0.5	&	0.503141363 $\cdots$	&	1.286737281$\cdots$	&	0.125261234$\cdots$	&	0.545169614$\cdots$	\\
0.6	&	0.607860912 $\cdots$	&	1.442513783$\cdots$	&	0.180782679$\cdots$	&	0.681212481$\cdots$	\\
0.7	&	0.717150254$\cdots$	&	1.659450438$\cdots$	&	0.246984703$\cdots$	&	0.835695885$\cdots$	\\
0.8	&	0.833926638$\cdots$	&	1.970197693$\cdots$	&	0.324460822$\cdots$	&	1.016197007$\cdots$	\\
0.9	&	0.962467275$\cdots$	&	2.439867095$\cdots$	&	0.414159959$\cdots$	&	1.234951913$\cdots$	\\
1	&	1.10910365$\cdots$	&	3.218145911$\cdots$	&	0.517553929$\cdots$	&	1.514209452$\cdots$	\\
1.1	&	1.283479121$\cdots$	&	4.739630017$\cdots$	&	0.636899254$\cdots$	&	1.90228478$\cdots$	\\
1.2	&	1.500980201$\cdots$	&	9.006810867$\cdots$	&	0.775675993$\cdots$	&	2.544535649$\cdots$	\\
1.3	&	1.787827508$\cdots$	&	90.67188404$\cdots$	&	0.939383753$\cdots$	&	4.853820909$\cdots$	\\
1.4	&	2.192925266$\cdots$	&		&	1.137129615$\cdots$	&		\\
1.5	&	2.819825139$\cdots$	&		&	1.385213908$\cdots$	&		\\
1.6	&	3.934210676$\cdots$	&		&	1.716804817$\cdots$	&		\\
1.7	&	6.489993451$\cdots$	&		&	2.216905555$\cdots$	&		\\
1.8	&	18.49292858$\cdots$	&		&	3.263963079$\cdots$	&		\\

\noalign{\smallskip}\hline
\end{tabular}
\\ (Note) Numerical data under the inequality $t<0$ are obtained using Eqs. (\ref{2.3}) and (\ref{2.4}).
\end{table}

%
\begin{figure*}[h]
\begin{center}
\includegraphics[width=0.75\textwidth]{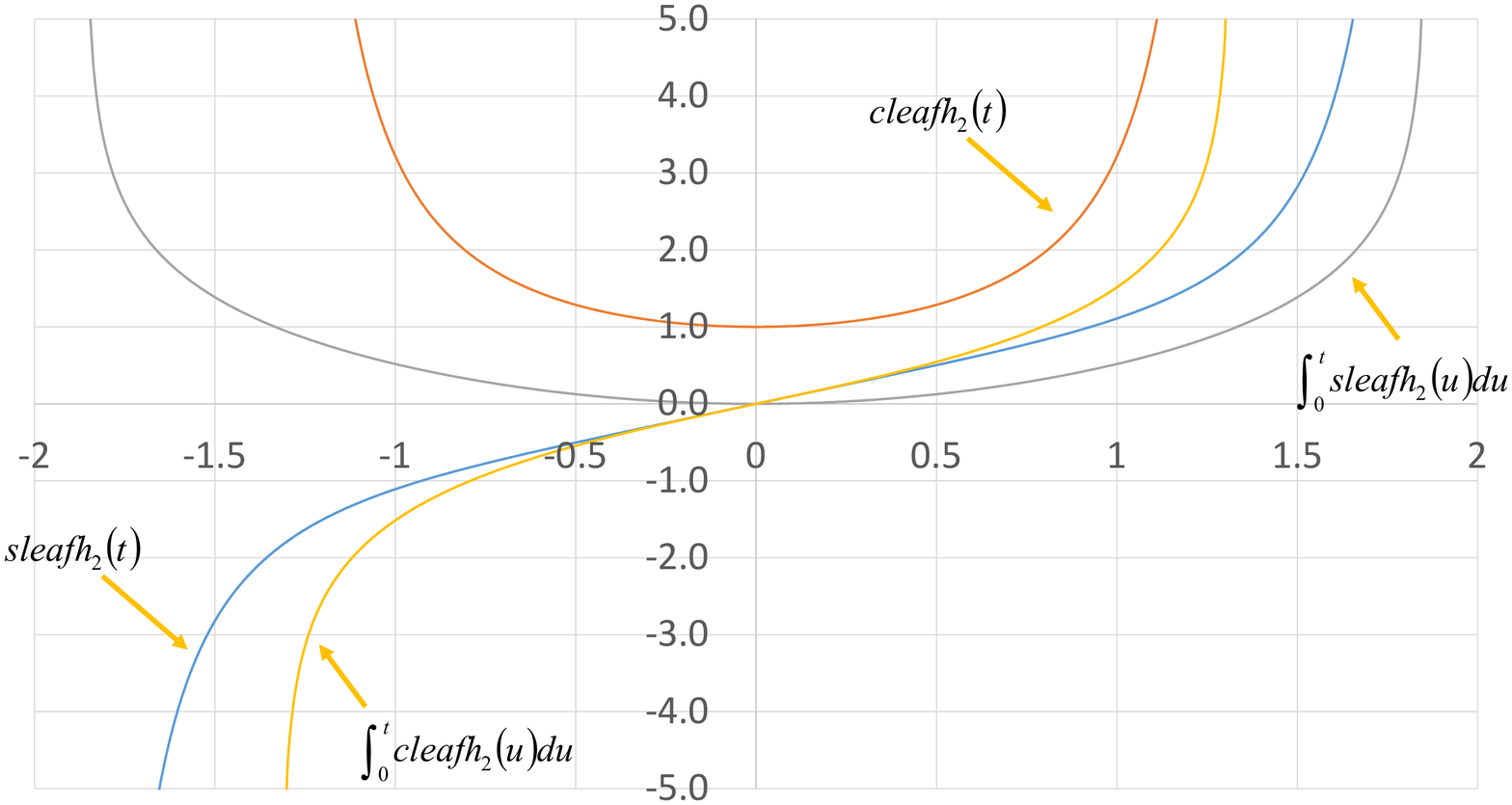}
\end{center}
\caption{Curves of hyperbolic leaf functions and their integral functions}
\label{fig:1}       
\end{figure*}

\section{Exact solution for the Duffing equation}
\label{Exact solution for the Duffing equation}
We consider the case where $\delta=0$ and $F=0$ in the Duffing equation of (\ref{Duffing}). Three types of divergence exact solutions were presented in the above equation. We define those ordinary differential equations and initial conditions as type (V\hspace{-.1em}I\hspace{-.1em}I\hspace{-.1em}I) to type (XIV). (See Ref. \cite{Kaz_duf}  for types (I) to (V\hspace{-.1em}I\hspace{-.1em}I)). 
The integral function of the leaf functions and the hyperbolic leaf functions are defined as follows:
\begin{equation}
\mathrm{SL}_2(t)=\int_{0}^{t} \mathrm{sleaf}_2(u) \mathrm{d}u \label{3.1}
\end{equation}
\begin{equation}
\mathrm{SLH}_2(t)=\int_{0}^{t} \mathrm{sleafh}_2(u)\mathrm{d}u \label{3.2}
\end{equation}
\begin{equation}
\mathrm{CLH}_2(t)=\int_{0}^{t} \mathrm{cleafh}_2(u)\mathrm{d}u \label{3.3}
\end{equation}
The variables $x(t)$, $A$, $t$, $\omega$, $u$, $\phi$ and $B$ represent the displacement, amplitude, time, angular frequency, parameter, initial phase and constant, respectively.

\subsection{Type (V\hspace{-.1em}I\hspace{-.1em}I\hspace{-.1em}I)  (See Appendix (V\hspace{-.1em}I\hspace{-.1em}I\hspace{-.1em}I )in\ detail)}
Exact\ Solution: 
\begin{equation}
 x(t)=A \cdot \mathrm{cosh}(\mathrm{CLH}_2(\omega t + \phi)) \label{3.1.1}
\end{equation}
Ordinary\ Differential\ Equation: 
\begin{equation}
 \frac{\mathrm{d}^2x(t)}{\mathrm{d}t^2}+3 \omega^2 x(t)- 4 (\frac{\omega}{A})^2 x(t)^3=0 \label{3.1.2}
\end{equation}
Initial\ Condition:  
\begin{equation}
x(0)=A \cdot \mathrm{cosh}(\mathrm{CLH}_2(\phi))  \label{3.1.3}
\end{equation}
Initial\ Velocity:  
\begin{equation}
 \frac{\mathrm{d}x(0)}{\mathrm{d}t}=A \cdot \mathrm{sinh}(\mathrm{CLH}_2(\phi)) \cdot \omega \cdot \mathrm{cleafh}_2(\phi)  \label{3.1.4}
\end{equation}

\subsection{Type (I\hspace{-.1em}X)  (See Appendix (I\hspace{-.1em}X)\ in detail)}
Exact\ Solution: 
\begin{equation}
 x(t)=A \cdot \mathrm{sinh}(\mathrm{CLH}_2(\omega t + \phi)) \label{3.2.1}
\end{equation}
Ordinary\ Differential\ Equation:  
\begin{equation}
\frac{\mathrm{d}^2x(t)}{\mathrm{d}t^2}-3 \omega^2 x(t)- 4 (\frac{\omega}{A})^2 x(t)^3=0
\label{3.2.2}
\end{equation}
Initial\ Condition:  
\begin{equation}
 x(0)=A \cdot \mathrm{sinh}(\mathrm{CLH}_2(\phi)) \label{3.2.3}
\end{equation}
Initial\ Velocity:  
\begin{equation}
 \frac{\mathrm{d}x(0)}{\mathrm{d}t}=A \cdot \mathrm{cosh}(\mathrm{CLH}_2(\phi)) \cdot \omega \cdot \mathrm{cleafh}_2(\phi) \label{3.2.4}
\end{equation}

\subsection{ Type (X)  (See Appendix (X)\ in detail)}
Exact\ Solution: 
\begin{equation}
\begin{split}
 x(t)=A \cdot \mathrm{cos}(\mathrm{SL}_2(\omega t + \phi))-A \cdot \mathrm{sin}(\mathrm{SL}_2(\omega t + \phi))+\sqrt{2} A \cdot \mathrm{cosh}(\mathrm{CLH}_2(\omega t + \phi)) 
\label{3.3.1}
\end{split}
\end{equation}
Ordinary\ Differential\ Equation: 
\begin{equation}
 \frac{\mathrm{d}^2x(t)}{\mathrm{d}t^2}+3 \omega^2 (1+2\sqrt{2})x(t)- 2 (\frac{\omega}{A})^2 x(t)^3=0
\label{3.3.2}
\end{equation}
Initial\ Condition:  
\begin{equation}
 x(0)=\sqrt{2} A \cdot \Bigl\{  \mathrm{cos}(\mathrm{SL}_2(\phi)+\frac{\pi}{4})
+\mathrm{cosh}(\mathrm{CLH}_2(\phi)) \Bigr\}
\label{3.3.3}
\end{equation}
Initial\ Velocity: 
\begin{equation}
\begin{split}
&\frac{\mathrm{d}x(0)}{\mathrm{d}t}=-A \cdot \mathrm{sin}(\mathrm{SL}_2(\phi)) \cdot \omega \cdot \mathrm{sleaf}_2(\phi) -A \cdot \mathrm{cos}(\mathrm{SL}_2(\phi)) \cdot \omega \cdot \mathrm{sleaf}_2(\phi) \\
&+\sqrt{2}A \cdot \mathrm{sinh}(\mathrm{CLH}_2(\phi)) \cdot \omega \cdot \mathrm{cleafh}_2(\phi) \label{3.3.4}
\end{split}
\end{equation}

\subsection{ Type (X\hspace{-.1em}I)  (See Appendix (X\hspace{-.1em}I)\ in detail) }
Next, consider the case where $F = 0$ in the Duffing equation(Eq. (\ref{Duffing})). In the exact solution that satisfies the Duffing equation, negative damping (divergence) based on type (X\hspace{-.1em}I) and type (X\hspace{-.1em}I\hspace{-.1em}I) are presented. The ordinary differential equations and initial conditions are as follows.

Exact\ Solution:
\begin{equation}
x(t)=A \cdot \mathrm{e}^{\omega t} \cdot \mathrm{sleafh}_2(B \cdot \mathrm{e}^{\omega t} + \phi)
\label{3.4.1}
\end{equation}
Ordinary\ Differential\ Equation: 
\begin{equation}
\frac{\mathrm{d}^2x(t)}{\mathrm{d}t^2}-3 \omega \frac{\mathrm{d}x(t)}{\mathrm{d}t}
+2 \omega^2 x(t)- 2 (\frac{B\omega}{A})^2 x(t)^3=0
\label{3.4.2}
\end{equation}
Initial\ Condition: 
\begin{equation}
x(0)=A \cdot \mathrm{sleafh}_2(B + \phi)
\label{3.4.3}
\end{equation}
Initial\ Velocity:
\begin{equation}
 \frac{\mathrm{d}x(0)}{\mathrm{d}t}=A \omega \cdot \mathrm{sleafh}_2(B + \phi) + A B \omega \cdot \sqrt{1+(\mathrm{sleafh}_2(B+\phi))^4} 
\label{3.4.4}
\end{equation}

\subsection{ Type (X\hspace{-.1em}I\hspace{-.1em}I)  (See Appendix (X\hspace{-.1em}I\hspace{-.1em}I)\ in detail)}
Exact\ Solution:
\begin{equation}
x(t)=A \cdot \mathrm{e}^{\omega t} \cdot \mathrm{cleafh}_2(B \cdot \mathrm{e}^{\omega t} + \phi)
\label{3.5.1}
\end{equation}
Ordinary\ Differential\ Equation: 
\begin{equation}
\frac{\mathrm{d}^2x(t)}{\mathrm{d}t^2}-3 \omega \frac{\mathrm{d}x(t)}{\mathrm{d}t}
+2 \omega^2 x(t)- 2 (\frac{B\omega}{A})^2 x(t)^3=0
\label{3.5.2}
\end{equation}
Initial\ Condition: 
\begin{equation}
x(0)=A \cdot \mathrm{cleafh}_2(B + \phi)
\label{3.5.3}
\end{equation}
Initial\ Velocity:
\begin{equation}
 \frac{\mathrm{d}x(0)}{\mathrm{d}t}=A \omega \cdot \mathrm{cleafh}_2(B + \phi) + A B \omega \cdot \sqrt{(\mathrm{cleafh}_2(B+\phi))^4-1} 
\label{3.5.4}
\end{equation}

\subsection{ Type (X\hspace{-.1em}I\hspace{-.1em}I\hspace{-.1em}I)  (See Appendix (X\hspace{-.1em}I\hspace{-.1em}I\hspace{-.1em}I)\ in detail) }
In the exact solution that satisfies Eq. (\ref{Duffing}) under the condition $F=0$, damping vibrations based on type (X\hspace{-.1em}I\hspace{-.1em}I\hspace{-.1em}I) and type (X\hspace{-.1em}I\hspace{-.1em}V) are presented. The ordinary differential equations and initial conditions are as follows.

Exact\ Solution:
\begin{equation}
x(t)=A \mathrm{e}^{\omega t} \cdot \mathrm{sleaf}_2(B \cdot \mathrm{e}^{\omega t} + \phi)
\label{3.6.1}
\end{equation}
Ordinary\ Differential\ Equation: 
\begin{equation}
\frac{\mathrm{d}^2x(t)}{\mathrm{d}t^2}-3 \omega \frac{\mathrm{d}x(t)}{\mathrm{d}t}
+2 \omega^2 x(t) + 2 (\frac{B\omega}{A})^2 x(t)^3=0
\label{3.6.2}
\end{equation}
Initial\ Condition: 
\begin{equation}
x(0)=A \cdot \mathrm{sleaf}_2(B + \phi)
\label{3.6.3}
\end{equation}
Initial\ Velocity:
\begin{equation}
 \frac{\mathrm{d}x(0)}{\mathrm{d}t}=A \omega \cdot \mathrm{sleaf}_2(B + \phi) + A B \omega \cdot \sqrt{1-(\mathrm{sleaf}_2(B+\phi))^4} 
\label{3.6.4}
\end{equation}

\subsection{Type (X\hspace{-.1em}I\hspace{-.1em}V)  (See Appendix (X\hspace{-.1em}I\hspace{-.1em}V)\ in detail) }
Exact\ Solution:
\begin{equation}
x(t)=A \mathrm{e}^{\omega t} \cdot \mathrm{cleaf}_2(B \cdot \mathrm{e}^{\omega t} + \phi)
\label{3.7.1}
\end{equation}
Ordinary\ Differential\ Equation: 
\begin{equation}
\frac{\mathrm{d}^2x(t)}{\mathrm{d}t^2}-3 \omega \frac{\mathrm{d}x(t)}{\mathrm{d}t}
+2 \omega^2 x(t) + 2 (\frac{B\omega}{A})^2 x(t)^3=0
\label{3.7.2}
\end{equation}
Initial\ Condition: 
\begin{equation}
x(0)=A \cdot \mathrm{cleaf}_2(B + \phi)
\label{3.7.3}
\end{equation}
Initial\ Velocity:
\begin{equation}
 \frac{\mathrm{d}x(0)}{\mathrm{d}t}= A \omega \cdot \mathrm{cleaf}_2(B + \phi) - A B \omega \cdot \sqrt{1-(\mathrm{cleaf}_2(B+\phi))^4} 
\label{3.7.4}
\end{equation}


\section{Analysis of results}
\label{Analysis of results}
\subsection{Divergence solution of type(V\hspace{-.1em}I\hspace{-.1em}I\hspace{-.1em}I)}
\label{Divergence solution1}
Let us consider a case where $A =1$, $\omega=1$ and $\phi=0$ in Eq. (\ref{3.1.1}), the integral function, the hyperbolic function $\mathrm{cosh}(t)$ and the hyperbolic leaf function $\mathrm{cleafh}_2(t)$ are shown in Fig. 2. The vertical and horizontal axes represent the variables $x(t)$ and $t$, respectively. As shown in Eq. (\ref{2.6}) and (\ref{2.9}), the hyperbolic leaf function $\mathrm{cleafh}_2(t)$ has limits \cite{Kaz_ch}. Therefore, the hyperbolic leaf function $\mathrm{cleafh}_2(t)$ increases sharply near the limit $t=\pm\eta_2$. The integral of the hyperbolic leaf function $\mathrm{cleafh}_2(t)$ is obtained as follows \cite{Kaz_ch}:

\begin{equation}
\begin{split}
\mathrm{CLH}_2(t)&=\int_{0}^{t} \mathrm{cleafh}_2(u) \mathrm{d}u \\
&=\mathrm{ln}(\sqrt{(\mathrm{cleafh}_2(t))^2+1}+\sqrt{(\mathrm{cleafh}_2(t))^2-1}) -\mathrm{ln}(\sqrt{2})
\end{split}
\end{equation}

The function $\mathrm{cleafh}_2(t)$ is included in the integral function, so that the integral function also increases sharply near the limit $t=\pm\eta_2$. In the exact solution of type (V\hspace{-.1em}I\hspace{-.1em}I\hspace{-.1em}I) including the integral function $\int_{0}^{l} \mathrm{cleafh}_2(t)\mathrm{d}t$, the (V\hspace{-.1em}I\hspace{-.1em}I\hspace{-.1em}I) exact solution increases monotonically near the limit $t=\pm\eta_2$. 
In the case that the curves vary with the amplitude $A$ under the conditions $\omega=1$ and $\phi=0$, the curves obtained by the (V\hspace{-.1em}I\hspace{-.1em}I\hspace{-.1em}I) solution are shown in Fig. 3. The limits $t=\pm\eta_2$ do not vary even if the parameter $A$ varies. For the inequality $A> 0$, the (V\hspace{-.1em}I\hspace{-.1em}I\hspace{-.1em}I) solution increases monotonically with time, $t$. For the inequality $A <0$, the curve of the (V\hspace{-.1em}I\hspace{-.1em}I\hspace{-.1em}I) solution decreases monotonically with $t$. The initial position $x(0)$ at $t=0$ (Eq. (\ref{3.1.3})) also varies according to the parameter $A$. In case that the curves vary with the phase $\omega$ under the conditions $A=1$ and $\phi=0$, the curves obtained by the (V\hspace{-.1em}I\hspace{-.1em}I\hspace{-.1em}I) solution are shown in Fig.4. The limit is $t=\pm\frac{\eta_2}{\omega}$. The limit varies according to the phase $\omega$. As the phase $\omega$ increases, the absolute value of the limit becomes smaller and the possible domain of the variable $t$ becomes narrower. As the phase $\omega$ decreases, the absolute value of the limit becomes larger and the possible domain of the variable $t$ become wider.

%
\begin{figure*}[tb]
\begin{center}
\includegraphics[width=0.75 \textwidth]{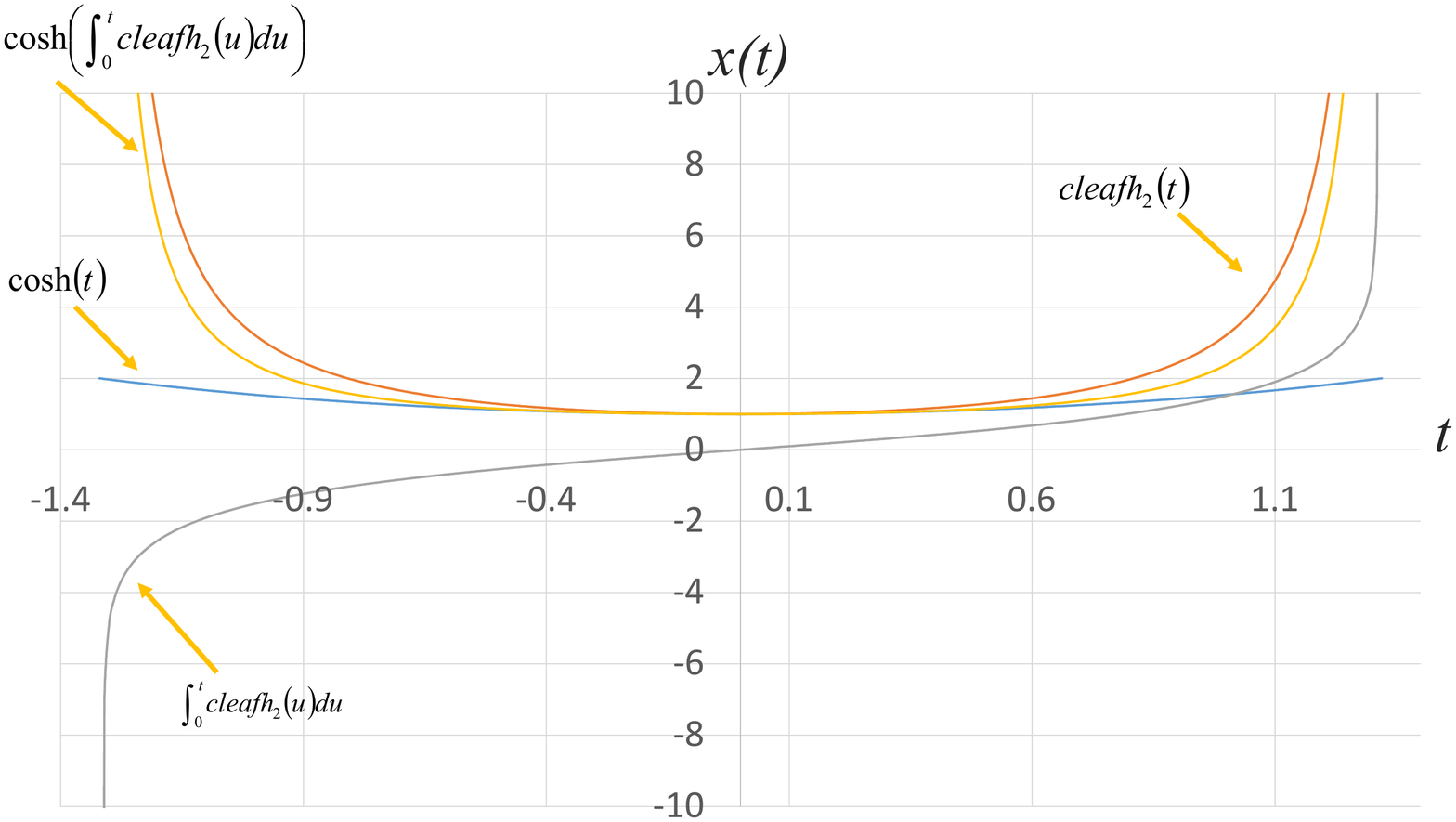}
\end{center}
\caption{ Curves obtained by the (V\hspace{-.1em}I\hspace{-.1em}I\hspace{-.1em}I) solution, the hyperbolic function $\mathrm{cosh}(t)$, the hyperbolic leaf function $\mathrm{cleafh}_2(t)$ and the function $\int_{0}^{t} \mathrm{cleafh}_2(t)\mathrm{d}t$ }
\label{fig:2}       
\end{figure*}

%
\begin{figure*}[tb]
\begin{center}
\includegraphics[width=0.75 \textwidth]{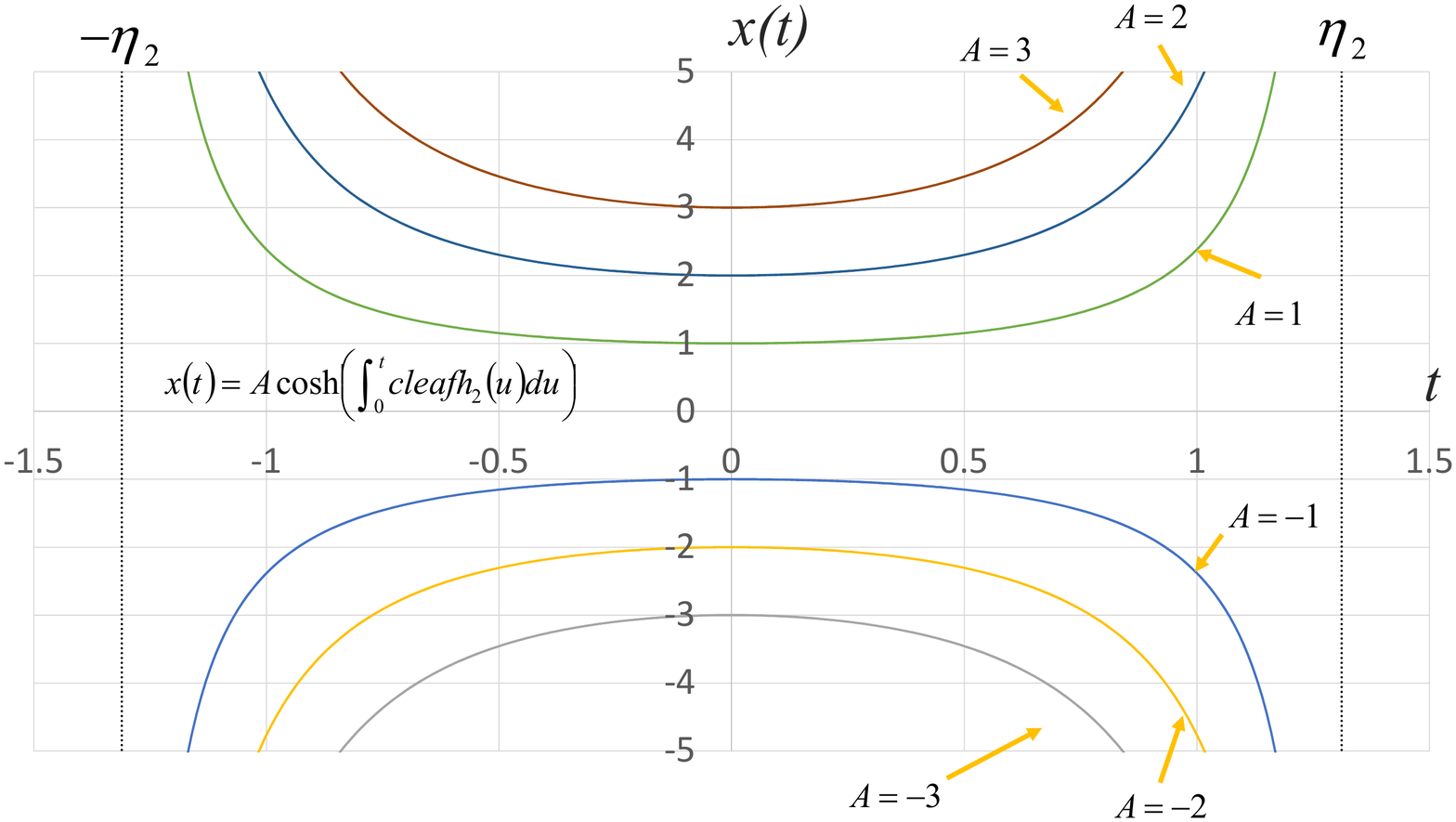}
\end{center}
\caption{Curves obtained by the (V\hspace{-.1em}I\hspace{-.1em}I\hspace{-.1em}I) exact solution as a function of variations in the amplitude $A$ ($A=\pm1$, $\pm2$, $\pm3$) (Set $\omega=1$ and $\phi=0$) }
\label{fig:3}       
\end{figure*}

%
\begin{figure*}[tb]
\begin{center}
\includegraphics[width=0.75 \textwidth]{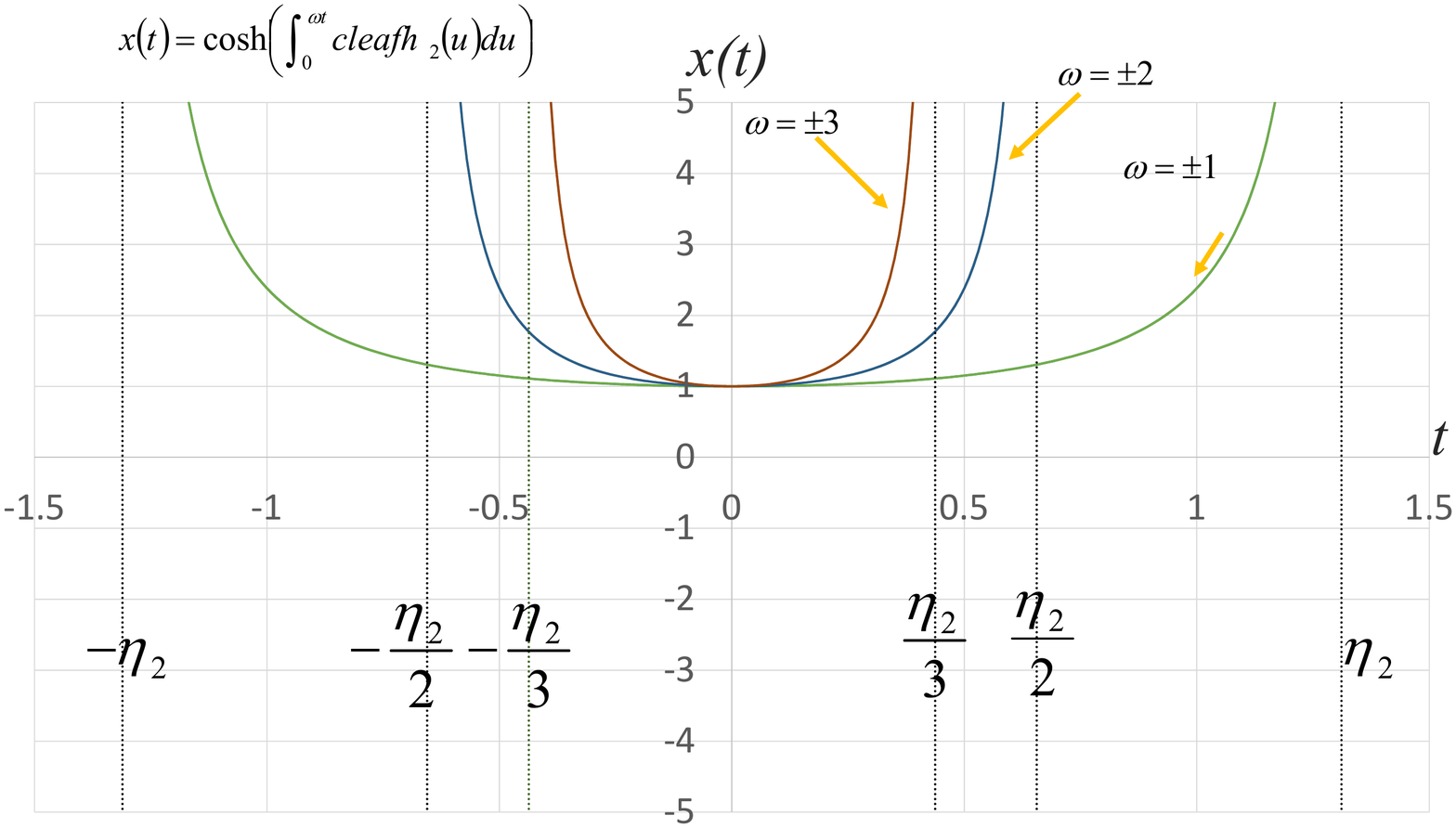}
\end{center}
\caption{Curves obtained by the (V\hspace{-.1em}I\hspace{-.1em}I\hspace{-.1em}I) exact solution as a function of variations in the angular frequency $\omega$  ($\omega=\pm1$, $\pm2$, $\pm3$) (Set $A=1$ and $\phi=0$)}
\label{fig:4}       
\end{figure*}

\begin{table}
\caption{ Numerical data for the (V\hspace{-.1em}I\hspace{-.1em}I\hspace{-.1em}I) solution under the conditions $A=1$, $\phi=0$ and $\omega=1$ }
\label{tab:2}      
\centering
\begin{tabular}{cccc}
\hline\noalign{\smallskip}
$t$ & $x(t)$ & $x(t)^3$  & $\frac{\mathrm{d}^2x(t)}{\mathrm{d}t^2}$ \\
\noalign{\smallskip}\hline\noalign{\smallskip}
-1.3	&	64.11860322 $\cdots$	&	263604.0988 $\cdots$	&	1054224.039 $\cdots$	\\
-1.2	&	6.407910814 $\cdots$	&	263.1172829 $\cdots$	&	1033.245399 $\cdots$	\\
-1.1	&	3.42520749 $\cdots$	&	40.18469304 $\cdots$	&	150.4631497 $\cdots$	\\
-1	&	2.382904017 $\cdots$	&	13.53068078 $\cdots$	&	46.97401106 $\cdots$	\\
-0.9	&	1.864530965 $\cdots$	&	6.48199663 $\cdots$	&	20.33439362 $\cdots$	\\
-0.8	&	1.562318622 $\cdots$	&	3.813368964 $\cdots$	&	10.56651999 $\cdots$	\\
-0.7	&	1.369995576 $\cdots$	&	2.571328091 $\cdots$	&	6.175325634 $\cdots$	\\
-0.6	&	1.241137787 $\cdots$	&	1.911877201 $\cdots$	&	3.924095444 $\cdots$	\\
-0.5	&	1.152322184 $\cdots$	&	1.530106881 $\cdots$	&	2.663460974 $\cdots$	\\
-0.4	&	1.090559612 $\cdots$	&	1.297024649 $\cdots$	&	1.916419762 $\cdots$	\\
-0.3	&	1.048200959 $\cdots$	&	1.151684863 $\cdots$	&	1.462136572 $\cdots$	\\
-0.2	&	1.020613924 $\cdots$	&	1.063125332 $\cdots$	&	1.190659556 $\cdots$	\\
-0.1	&	1.005037714 $\cdots$	&	1.015189405 $\cdots$	&	1.045644478 $\cdots$	\\
0	&	1	&	1	&	1	\\
0.1	&	1.005037714 $\cdots$	&	1.015189405 $\cdots$	&	1.045644478 $\cdots$	\\
0.2	&	1.020613924 $\cdots$	&	1.063125332 $\cdots$	&	1.190659556 $\cdots$	\\
0.3	&	1.048200959 $\cdots$	&	1.151684863 $\cdots$	&	1.462136572 $\cdots$	\\
0.4	&	1.090559612 $\cdots$	&	1.297024649 $\cdots$	&	1.916419762 $\cdots$	\\
0.5	&	1.152322184 $\cdots$	&	1.530106881 $\cdots$	&	2.663460974 $\cdots$	\\
0.6	&	1.241137787 $\cdots$	&	1.911877201 $\cdots$	&	3.924095444 $\cdots$	\\
0.7	&	1.369995576 $\cdots$	&	2.571328091 $\cdots$	&	6.175325634 $\cdots$	\\
0.8	&	1.562318622 $\cdots$	&	3.813368964 $\cdots$	&	10.56651999 $\cdots$	\\
0.9	&	1.864530965 $\cdots$	&	6.48199663 $\cdots$	&	20.33439362 $\cdots$	\\
1	&	2.382904017 $\cdots$	&	13.53068078 $\cdots$	&	46.97401106 $\cdots$	\\
1.1	&	3.42520749 $\cdots$	&	40.18469304 $\cdots$	&	150.4631497 $\cdots$	\\
1.2	&	6.407910814 $\cdots$	&	263.1172829 $\cdots$	&	1033.245399 $\cdots$	\\
1.3	&	64.11860322 $\cdots$	&	263604.0988 $\cdots$	&	1054224.039 $\cdots$	\\
\noalign{\smallskip}\hline
\end{tabular}
\end{table}

\subsection{Divergence solutions of the type(I\hspace{-.1em}X)}
\label{Divergence solution9}
Let us consider a case where $A =1$, $\omega=1$ and $\phi=0$ in Eq. (\ref{3.2.1}), the integral function, the hyperbolic function $\mathrm{sinh}(t)$, and the hyperbolic leaf function $\mathrm{cleafh}_2(t)$ are shown in Fig. 5. The vertical and horizontal axes represent the variables $x(t)$ and $t$, respectively. The hyperbolic leaf function $\mathrm{cleafh}_2(t)$ has limits \cite {Kaz_ch}. Therefore, the exact solution of type (I\hspace{-.1em}X) increases monotonically near the limit $t=\eta_2$, and decreases monotonically near the limit $t=-\eta_2$. In case that curves vary with the amplitude $A$ under the conditions $\omega=1$ and $\phi=0$, the curves obtained by the (I\hspace{-.1em}X) solution are shown as in Fig.6. The limits of $t=\pm\eta_2$ do not vary even if the parameter $A$ varies. For the inequality $A>0$, the (I\hspace{-.1em}X) solution increases monotonically with time, $t$. For the inequality $A<0$, the (I\hspace{-.1em}X) solution decreases monotonically with $t$. The curve passes through $x(0)=0$. In case that the curves vary with the amplitude $\omega$ under the conditions $A=1$ and $\phi=0$, the curves obtained by the (I\hspace{-.1em}X) solution are shown as in Fig.7. The limit is $t=\pm\frac{\eta_2}{\omega}$. The limit varies according to the phase $\omega$. As the phase $\omega$ increases, the absolute value of the limit decreases and the possible domain of the variable $t$ becomes narrower. As the phase $\omega$ decreases, the absolute value of the limit increases and the possible domain of the variable $t$ becomes wider.

%
\begin{figure*}[tb]
\begin{center}
\includegraphics[width=0.75 \textwidth]{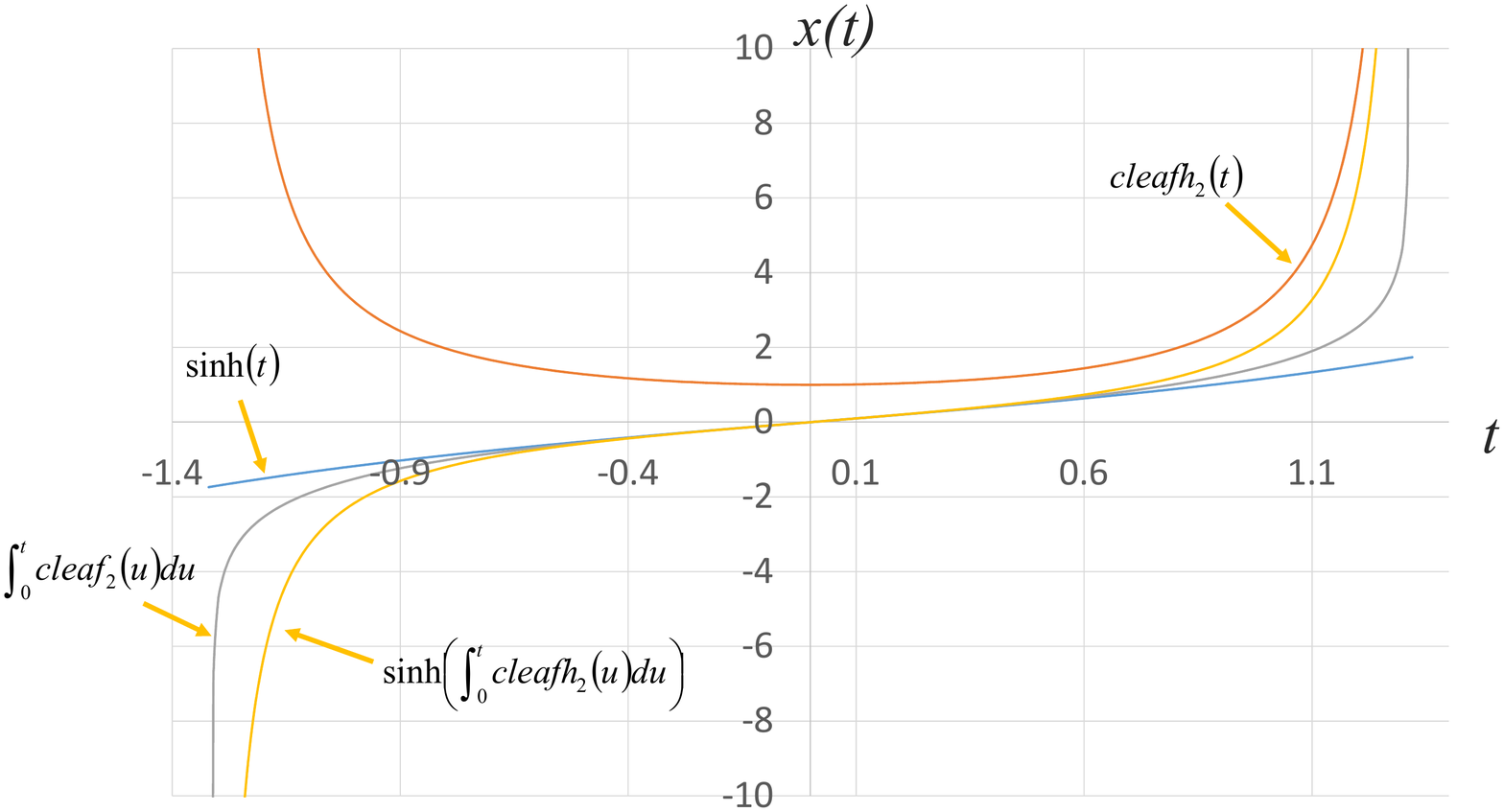}
\end{center}
\caption{Curves obtained by the (I\hspace{-.1em}X) exact solution, the hyperbolic function $\mathrm{sinh}(t)$, the leaf function $\mathrm{cleafh}_2(t)$ and the function $\int_{0}^{t} \mathrm{cleafh}_2(u)\mathrm{d}u$ }
\label{fig:5}       
\end{figure*}

%
\begin{figure*}[tb]
\begin{center}
\includegraphics[width=0.75 \textwidth]{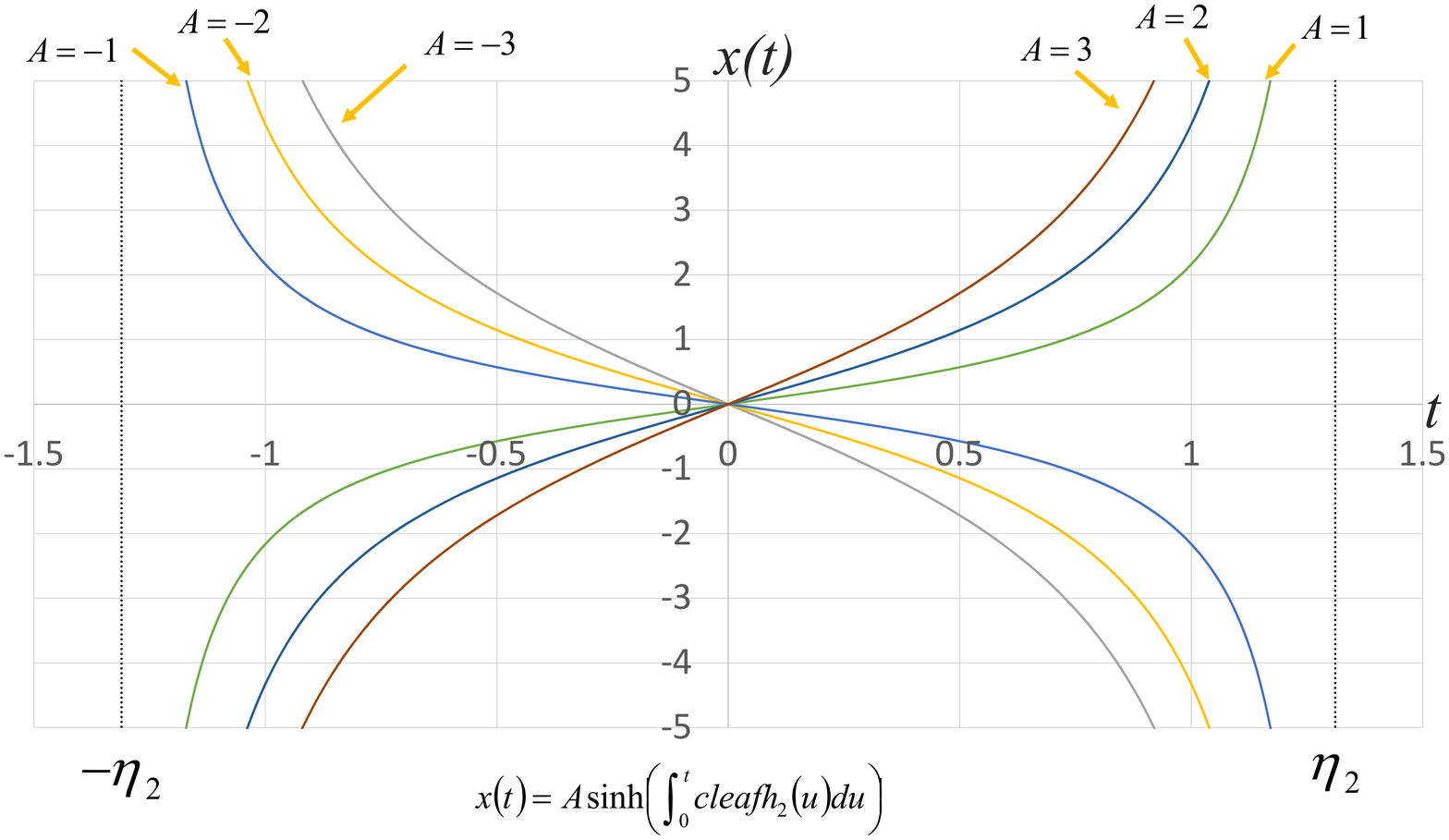}
\end{center}
\caption{Curves obtained by the (I\hspace{-.1em}X) exact solution as a function of variations in the amplitude $A$  ($A=\pm1$, $\pm2$, $\pm3$) (Set $\omega=1$ and $\phi=0$)}
\label{fig:6}       
\end{figure*}

%
\begin{figure*}[tb]
\begin{center}
\includegraphics[width=0.75 \textwidth]{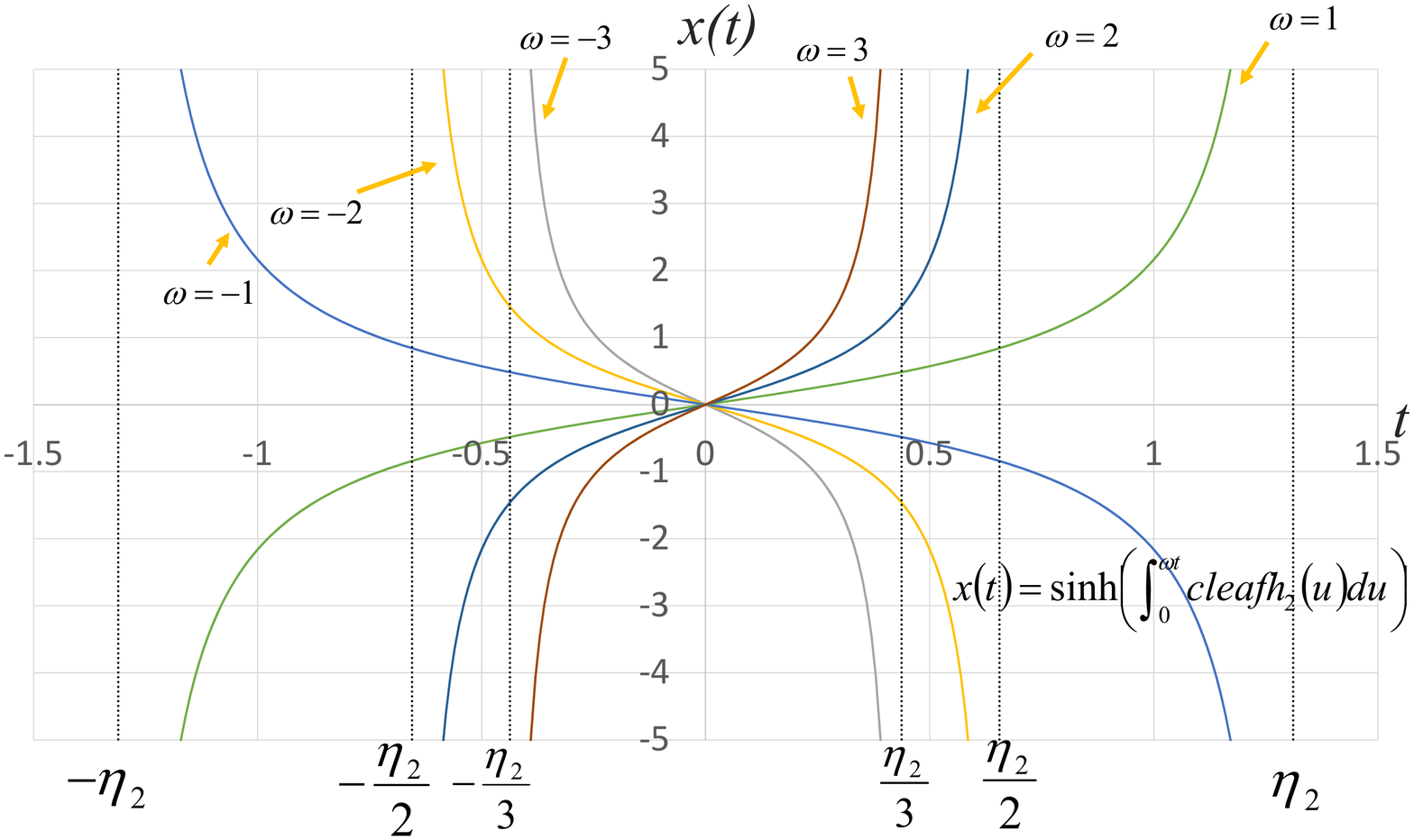}
\end{center}
\caption{Curves obtained by the (I\hspace{-.1em}X) exact solution as a function of variations in the amplitude $\omega$  ($\omega=\pm1$, $\pm2$, $\pm3$) (Set $A=1$ and $\phi=0$)}
\label{fig:7}       
\end{figure*}

\begin{table}
\caption{ Numerical data for the (I\hspace{-.1em}X) solution under the conditions $A=1$, $\phi=0$ and $\omega=1$  }
\label{tab:3}     
\centering
\begin{tabular}{cccc}
\hline\noalign{\smallskip}
$t$ & $x(t)$ & $x(t)^3$  & $\frac{\mathrm{d}^2x(t)}{\mathrm{d}t^2}$ \\
\noalign{\smallskip}\hline\noalign{\smallskip}
-1.3	&	-64.11080469 $\cdots$	&	-263507.9268	 $\cdots$&	-1054272.011	 $\cdots$\\
-1.2	&	-6.329401315 $\cdots$	&	-253.5641778	 $\cdots$&	-1033.249417 $\cdots$	\\
-1.1	&	-3.275980212 $\cdots$	&	-35.15797148	 $\cdots$&	-150.4601587	 $\cdots$\\
-1	&	-2.162921994 $\cdots$	&	-10.11864992	 $\cdots$&	-46.96343323	 $\cdots$\\
-0.9	&	-1.573682217 $\cdots$	&	-3.897185801	 $\cdots$&	-20.30981105 $\cdots$	\\
-0.8	&	-1.20034973 $\cdots$	&	-1.729511276 $\cdots$	&	-10.51910279	 $\cdots$\\
-0.7	&	-0.936422917 $\cdots$	&	-0.821137905 $\cdots$	&	-6.093824366	 $\cdots$\\
-0.6	&	-0.735134686 $\cdots$	&	-0.397283697 $\cdots$	&	-3.794540983	 $\cdots$\\
-0.5	&	-0.572578742 $\cdots$	&	-0.187717888	 $\cdots$&	-2.468609061	 $\cdots$\\
-0.4	&	-0.435109489 $\cdots$	&	-0.082375045 $\cdots$	&	-1.634829526	 $\cdots$\\
-0.3	&	-0.314205747 $\cdots$	&	-0.031020041 $\cdots$	&	-1.066698123	 $\cdots$\\
-0.2	&	-0.204090129 $\cdots$	&	-0.008500921 $\cdots$	&	-0.646274802	 $\cdots$\\
-0.1	&	-0.100502766 $\cdots$	&	-0.001015159 $\cdots$	&	-0.305570034 $\cdots$	\\
0	&	0	&	0	&	0	\\
0.1	&	0.100502766 $\cdots$	&	0.001015159 $\cdots$	&	0.305570034	 $\cdots$\\
0.2	&	0.204090129 $\cdots$	&	0.008500921 $\cdots$	&	0.646274802	 $\cdots$\\
0.3	&	0.314205747 $\cdots$	&	0.031020041 $\cdots$	&	1.066698123	 $\cdots$\\
0.4	&	0.435109489 $\cdots$	&	0.082375045 $\cdots$	&	1.634829526	 $\cdots$\\
0.5	&	0.572578742 $\cdots$	&	0.187717888 $\cdots$	&	2.468609061	 $\cdots$\\
0.6	&	0.735134686 $\cdots$	&	0.397283697 $\cdots$	&	3.794540983	 $\cdots$\\
0.7	&	0.936422917 $\cdots$	&	0.821137905 $\cdots$	&	6.093824366	 $\cdots$\\
0.8	&	1.20034973 $\cdots$	&	1.729511276 $\cdots$	&	10.51910279	 $\cdots$\\
0.9	&	1.573682217 $\cdots$	&	3.897185801 $\cdots$	&	20.30981105	 $\cdots$\\
1	&	2.162921994 $\cdots$	&	10.11864992 $\cdots$	&	46.96343323	 $\cdots$\\
1.1	&	3.275980212 $\cdots$	&	35.15797148 $\cdots$	&	150.4601587	 $\cdots$\\
1.2	&	6.329401315 $\cdots$	&	253.5641778 $\cdots$	&	1033.249417	 $\cdots$\\
1.3	&	64.11080469 $\cdots$	&	263507.9268 $\cdots$	&	1054272.011	 $\cdots$\\

\noalign{\smallskip}\hline
\end{tabular}
\end{table}

\subsection{Divergence solution of the type (X)}
\label{Divergence solution10}
Let us consider a case where $A =1$, $\omega=1$ and $\phi=0$ in Eq. (\ref{3.3.1}). The (X) exact solution, the hyperbolic function $\mathrm{cosh}(t)$ and the hyperbolic leaf function $\mathrm{cleafh}_2(t)$ are shown in Fig. 8. The vertical and horizontal axes represent the variables $x(t)$ and $t$, respectively. As shown in Fig. 8, these functions are even functions that are symmetrical about the $x(t)$ axis. The limits exist in the (X) exact solution. The (X) exact solution is transformed as follows:
\begin{equation}
x_1(t)=\mathrm{cos}(\mathrm{SL}_2(t))-\mathrm{sin}(\mathrm{SL}_2(t)) \label{4.3.1}
\end{equation}
Using the Eq. (III.4) in Ref. \cite{Kaz_duf}, squaring both sides of the above equation yields the following equation.
\begin{equation}
\begin{split}
(x_1(t))^2&=(\mathrm{cos}(\mathrm{SL}_2(t))^2+(\mathrm{sin}(\mathrm{SL}_2(t))^2 -2 \mathrm{sin}(\mathrm{SL}_2(t)) \cdot \mathrm{cos}(\mathrm{SL}_2(t)) \\
&=1-\mathrm{sin}(2 \cdot \mathrm{SL}_2(t))=1-(\mathrm{sleaf}_2(t))^2 \label{4.3.2}
\end{split}
\end{equation}
The variable $x_1(t)$ is obtained as follows:
\begin{equation}
x_1(t)=\pm\sqrt{1-(\mathrm{sleaf}_2(t))^2} \label{4.3.3}
\end{equation}
Next, the following equation is transformed:
\begin{equation}
x_2(t)=\sqrt{2}\mathrm{cosh}(\mathrm{CLH}_2(t)) \label{4.3.4}
\end{equation}
Using the Eq. (H5) in Ref. \cite{Kaz_ch}, squaring both sides of the above equation gives the following equation.
\begin{equation}
\begin{split}
(x_2(t))^2 =2 (\mathrm{cosh}( \mathrm{CLH}_2(t)) )^2 =1+\mathrm{cosh}(2 \mathrm{CLH}_2(t) )=1+(\mathrm{cleafh}_2(t))^2
\label{4.3.5}
\end{split}
\end{equation}
As shown in Fig. 8, the inequality $x_2(t)>0$ is obvious. The following equation is obtained.
\begin{equation}
x_2(t)=\sqrt{1 + (\mathrm{cleafh}_2(t))^2} \label{4.3.6}
\end{equation}
Therefore, the (X) exact solution is obtained as follows:
\begin{equation}
x(t)=x_1(t)+x_2(t)=\pm\sqrt{1 - (\mathrm{sleaf}_2(t))^2}+\sqrt{1 + (\mathrm{cleafh}_2(t))^2} \label{4.3.7}
\end{equation}
Because the (X) exact solution contains the function $\mathrm{cleafh}_2(t)$, its exact solution also has limits according to the limits of $\mathrm{cleafh}_2(t)$. The domain of (X) is as follows.
\begin{equation}
-\eta_2<t<\eta_2 \label{4.3.8}
\end{equation}
The sign $\eta_2$ is a constant described by Eq. (\ref{2.6}). The sign $\eta_2$ has the following relation given by Eq.  (24) in Ref. \cite{Kaz_ch}.
\begin{equation}
\eta_2=\frac{\pi_2}{2} \ \ \ \ \ (\int_{1}^{\infty} \frac{1}{\sqrt{t^4-1}}\mathrm{d}t = \int_{0}^{1} \frac{1}{\sqrt{1-t^4}}\mathrm{d}t ) \label{4.3.9}
\end{equation}
In the case that the curves vary with the amplitude $A$ under the conditions $\omega=1$ and $\phi=0$, the curves obtained by the (X) solution are shown in Fig. 9. The limit $t=\pm\eta_2$ does not vary even if the parameter $A$ varies. Conversely, the initial condition (Eq. (\ref{3.3.3})) varies with the amplitude, $A$. In case that the curves vary with the amplitude $\omega$ under the conditions $A=1$ and $\phi=0$, the curves obtained by the (X) solution are shown in Fig. 10. The limit is $t=\pm\frac{\eta_2}{\omega}$. The limit varies with the phase, $\omega$. As the phase $\omega$ increases, the absolute value of the limit decreases and the possible domain of the variable $t$ becomes narrower. As the phase $\omega$ decreases, the absolute value of the limit increases and the possible domain of the variable $t$ becomes larger.

%
\begin{figure*}[tb]
\begin{center}
\includegraphics[width=0.75 \textwidth]{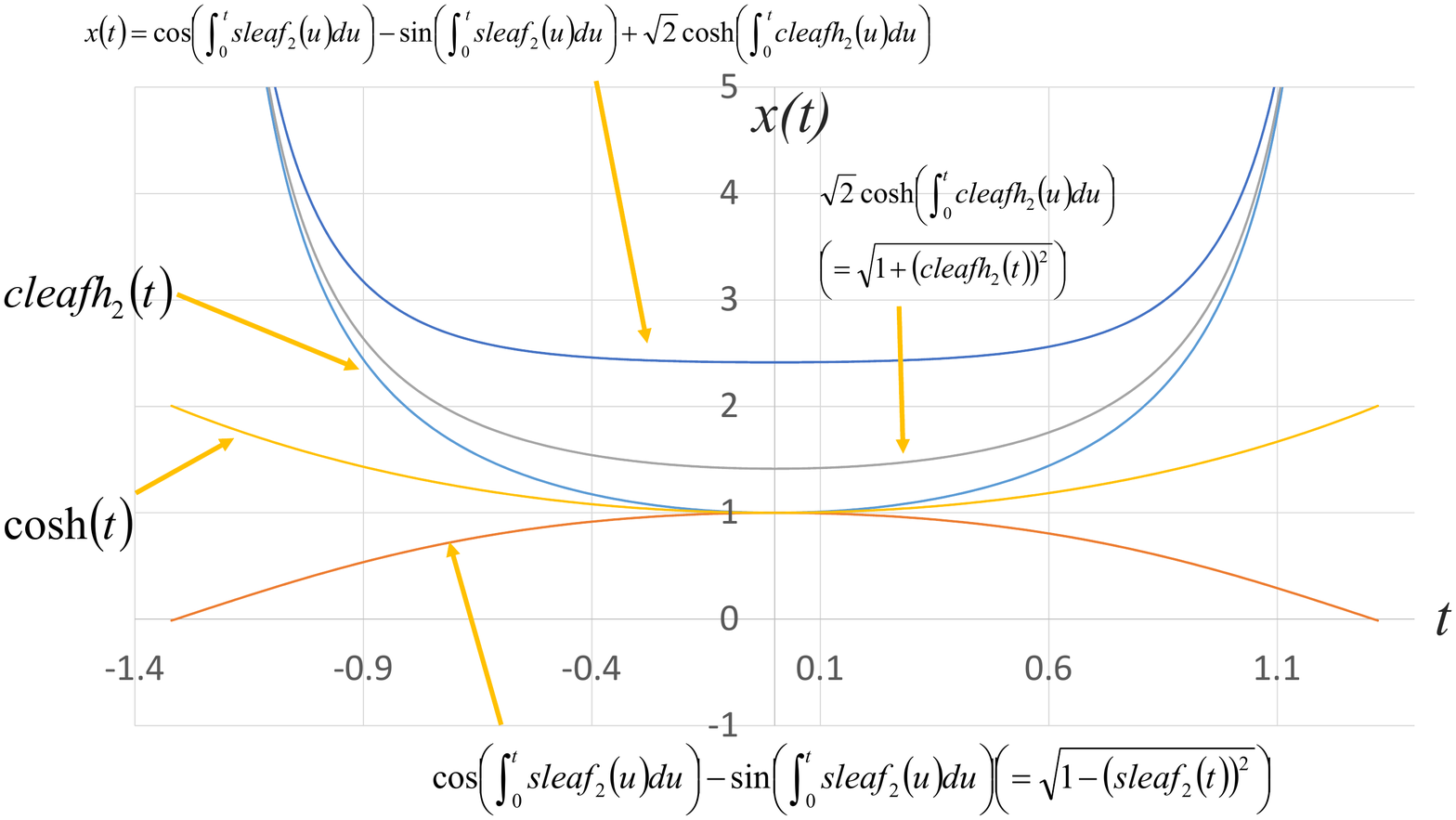}
\end{center}
\caption{ Curves obtained by the (X) exact solution, the hyperbolic function $\mathrm{cosh}(t)$, the leaf function $\mathrm{cleafh}_2(t)$, the function
$\sqrt{1+(\mathrm{cleafh}_2(t))^2}$  and the function $\sqrt{1-(\mathrm{sleafh}_2(t))^2}$  }
\label{fig:8}       
\end{figure*}

%
\begin{figure*}[tb]
\begin{center}
\includegraphics[width=0.75 \textwidth]{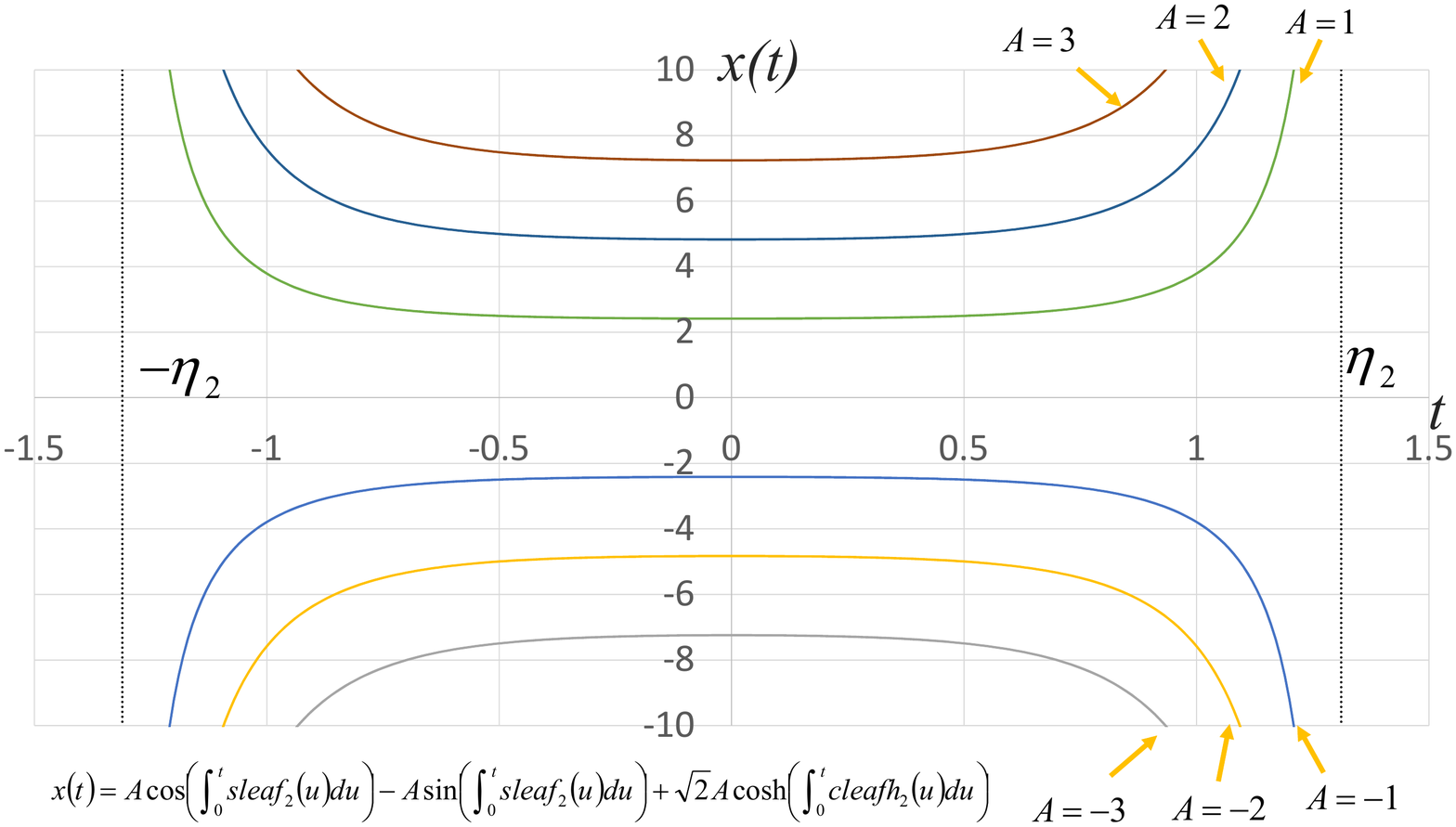}
\end{center}
\caption{ Curves obtained by the (X) exact solution as a function of variations in the amplitude $A$ ($A=\pm1$, $\pm2$, $\pm3$) (Set $\omega=1$ and $\phi=0$)}
\label{fig:9}       
\end{figure*}

%
\begin{figure*}[tb]
\begin{center}
\includegraphics[width=0.75 \textwidth]{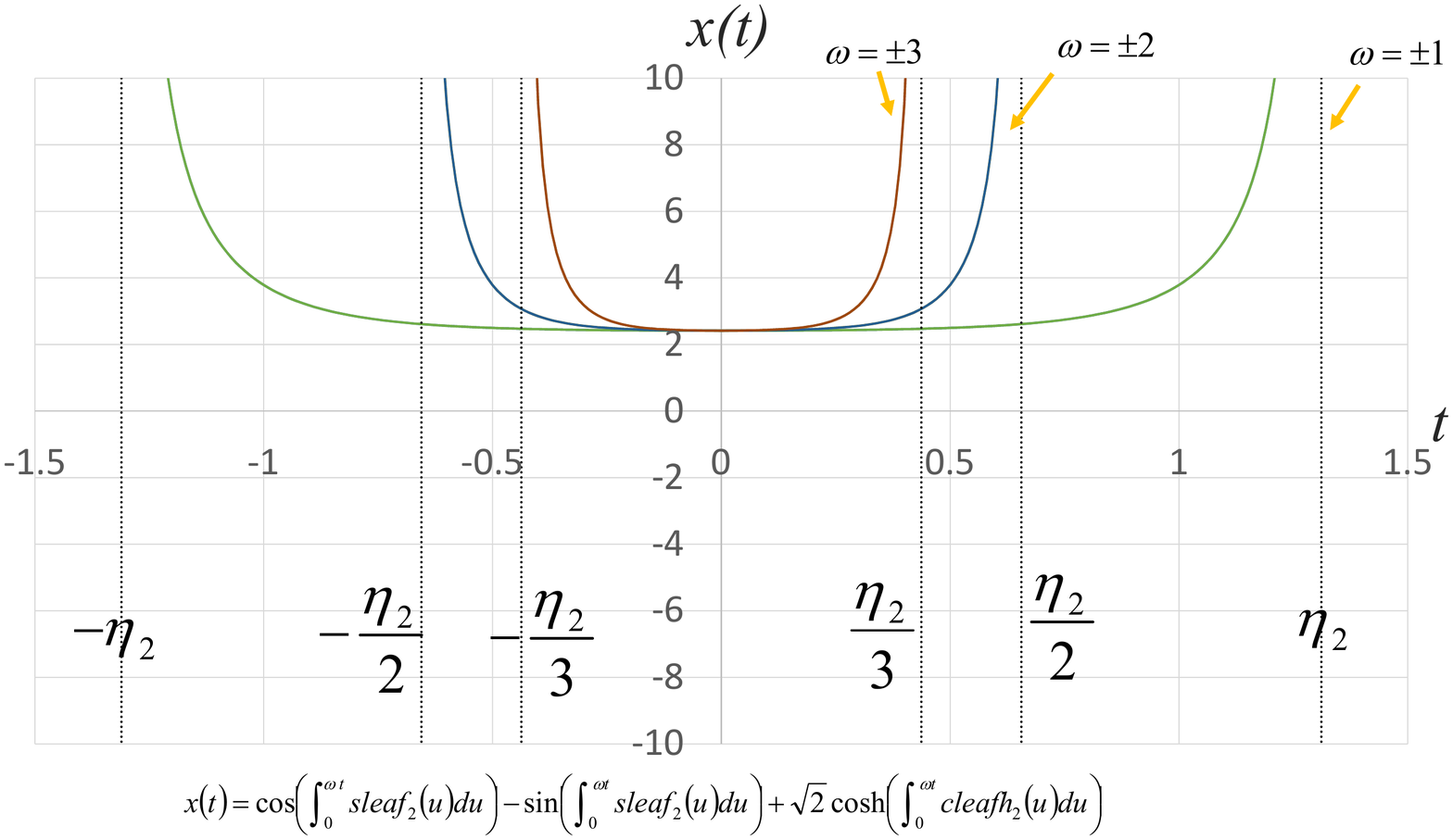}
\end{center}
\caption{ Curves obtained by the (X) exact solution as a function of variations in the parameter $\omega$ ($\omega=\pm1$, $\pm2$, $\pm3$) (Set $A=1$ and $\phi=0$)}
\label{fig:10}       
\end{figure*}

\begin{table}
\caption{ Numerical data for the (X) solution under the conditions $A=1$, $\phi=0$ and $\omega=1$ }
\label{tab:4}    
\centering
\begin{tabular}{cccc}
\hline\noalign{\smallskip}
$t$ & $x(t)$ & $x(t)^3$  & $\frac{\mathrm{d}^2x(t)}{\mathrm{d}t^2}$ \\
\noalign{\smallskip}\hline\noalign{\smallskip}
-1.3	&	90.69299413 $\cdots$	&	745969.7552 $\cdots$	&	1490897.888 $\cdots$	\\
-1.2	&	9.218211398 $\cdots$	&	783.3213979 $\cdots$	&	1460.769087 $\cdots$	\\
-1.1	&	5.135927878 $\cdots$	&	135.4742481 $\cdots$	&	211.9609381 $\cdots$	\\
-1	&	3.789591734  $\cdots$	&	54.42234777 $\cdots$	&	65.32012665 $\cdots$	\\
-0.9	&	3.173173112  $\cdots$	&	31.95076764 $\cdots$	&	27.45673813  $\cdots$	\\
-0.8	&	2.849526806 $\cdots$	&	23.13759636	 $\cdots$&	13.54756343 $\cdots$	\\
-0.7	&	2.667396157 $\cdots$	&	18.97852968	 $\cdots$&	7.321249434 $\cdots$	\\
-0.6	&	2.560946671 $\cdots$	&	16.79583519	 $\cdots$&	4.178464757 $\cdots$	\\
-0.5	&	2.497442952 $\cdots$	&	15.57710437	 $\cdots$&	2.470362978 $\cdots$	\\
-0.4	&	2.459244636 $\cdots$	&	14.87322673	 $\cdots$&	1.501336344 $\cdots$	\\
-0.3	&	2.436395606 $\cdots$	&	14.46250169	 $\cdots$&	0.94231307 $\cdots$	\\
-0.2	&	2.423168514 $\cdots$	&	14.22822918	 $\cdots$&	0.625685442  $\cdots$	\\
-0.1	&	2.416325506 $\cdots$	&	14.10802806	 $\cdots$&	0.463877721 $\cdots$	\\
0	&	2.414213562 $\cdots$	&	14.07106781  $\cdots$	&	0.414213562 $\cdots$	\\
0.1	&	2.416325506  $\cdots$	&	14.10802806  $\cdots$	&	0.463877721 $\cdots$	\\
0.2	&	2.423168514 $\cdots$	&	14.22822918	 $\cdots$&	0.625685442 $\cdots$	\\
0.3	&	2.436395606 $\cdots$	&	14.46250169	 $\cdots$&	0.94231307 $\cdots$	\\
0.4	&	2.459244636 $\cdots$	&	14.87322673	 $\cdots$&	1.501336344 $\cdots$	\\
0.5	&	2.497442952 $\cdots$	&	15.57710437	 $\cdots$&	2.470362978 $\cdots$	\\
0.6	&	2.560946671 $\cdots$	&	16.79583519	 $\cdots$&	4.178464757 $\cdots$	\\
0.7	&	2.667396157 $\cdots$	&	18.97852968	 $\cdots$&	7.321249434 $\cdots$	\\
0.8	&	2.849526806 $\cdots$	&	23.13759636	 $\cdots$&	13.54756343 $\cdots$	\\
0.9	&	3.173173112 $\cdots$	&	31.95076764	 $\cdots$&	27.45673813 $\cdots$	\\
1	&	3.789591734 $\cdots$	&	54.42234777 $\cdots$	&	65.32012665 $\cdots$	\\
1.1	&	5.135927878  $\cdots$	&	135.4742481  $\cdots$	&	211.9609381 $\cdots$	\\
1.2	&	9.218211398 $\cdots$	&	783.3213979  $\cdots$	&	1460.769087 $\cdots$	\\
1.3	&	90.69299413 $\cdots$	&	745969.7552 $\cdots$	&	1490897.888  $\cdots$	\\

\noalign{\smallskip}\hline
\end{tabular}
\end{table}

\subsection{Divergence solution of the type (X\hspace{-.1em}I)}
\label{Divergence solution11}
Let us consider a case where $A =1$, $\omega=1$, $\phi=-1$ and $B=1$ in Eq. (\ref{3.4.1}), the (X\hspace{-.1em}I) exact solution, the hyperbolic function $\mathrm{sinh}(t)$ and the hyperbolic leaf function $\mathrm{sleafh}_2(t)$ are shown in Fig. 11. Based on Eq. (\ref{3.4.3}), the curve passes through $x(0)=0$. For the inequality $t<0$, the variable $x(t)$ converges asymptotically to zero. For $t>0$, the variable $x(t)$ increases monotonically. The limit exists for the hyperbolic leaf function $\mathrm{sleafh}_2(\mathrm{e}^t-1)$ (See Ref. \cite{Kaz_sh}). The possible range of $\mathrm{e}^t-1$ is as follows:
\begin{equation}
-\zeta_2<\mathrm{e}^t-1<\zeta_2 \label{4.4.1}
\end{equation}
The following equation is obtained from the above equation.
\begin{equation}
-\zeta_2+1<\mathrm{e}^t<\zeta_2+1 \label{4.4.2}
\end{equation}
The constant $\zeta_2$ is obtained by the Eq. (\ref{2.5}) (See Ref. \cite{Kaz_sh}) .
In Eq. (\ref{4.4.2}), the inequality $ - \zeta_2+1=-0.85407\cdots<\mathrm{e}^t$ always holds for any arbitrary time, $t$. We can take the logarithm of both sides in the inequality $\mathrm{e}^t<\zeta_2+1$ in Eq. (\ref{4.4.2}). Solving for the variable $t$ yields the following equation.
\begin{equation}
t<\mathrm{ln}(\zeta_2+1) \ \ (=\mathrm{ln}(1.85407 \cdots+1)=1.0487 \cdots) \label{4.4.3}
\end{equation}
The above inequality represents the domain of the variable $t$ in Eq. (\ref{4.4.1}). In case that curves vary with the amplitude $A$ under the conditions $\omega=1$, $\phi=-1$ and $B=1$, the curves obtained by the (X\hspace{-.1em}I) solution are  shown in Fig.12. The limits do not vary even if the parameter $A$ varies.
\begin{equation}
t=\mathrm{ln}(\zeta_2+1)=\mathrm{ln}(1.85407 \cdots+1)=1.0487 \cdots
\end{equation}
For the inequality $A>0$, the exact solution $x(t)$ of (X\hspace{-.1em}I) diverges to plus infinity. 
On the other hand, for the inequality $A<0$, the exact solution $x(t)$ of (X\hspace{-.1em}I) diverges to minus infinity. In case that the curves vary with the phase $\omega$ under the conditions $A=1$, $\phi=-1$ and $B=1$, the curves obtained by the (X\hspace{-.1em}I) solution are shown in Fig. 13. 
Let us consider the limit of the (X\hspace{-.1em}I) exact solution when the parameter $\omega$ varies. The limit of the variable $t$ is obtained by the following equation.
\begin{equation}
\mathrm{e}^{\omega t}-1=\pm\zeta_2=\pm 1.85407 \cdots \label{4.4.4}
\end{equation}
When the right side of Eq. (\ref{4.4.4}) is $1.85407 \cdots$, solving the above equation for the variable $t$ yields the following equation.
\begin{equation}
t=\frac{1}{\omega} \mathrm{ln}(\zeta_2+1) \label{4.4.5}
\end{equation}
Let us consider a case where the right side of Eq. (\ref{4.4.4}) is -$1.85407 \cdots$.
\begin{equation}
\mathrm{e}^{\omega t}-1=-1.85407 \cdots \label{4.4.6}
\end{equation}
The variable $t$ that satisfies the above equation does not exist. Therefore, the condition for the existence of the limit is $\zeta_2=1.85407 \cdots$. When $\omega> 0$ in Eq. (\ref{4.4.5}), the limit of the variable $t$ is generated in the domain where $t> 0$. When $\omega< 0$ in Eq. (\ref{4.4.5}), the limit of the variable $t$ is produced in the domain $t< 0$.  In case that the curves vary with the phase $A$ under the conditions $\omega=1$, $\phi=-2$ and $B=2$, the curves obtained by the (X\hspace{-.1em}I) solution are shown in Fig.14. The possible range of  $2\mathrm{e}^t-2$ in the above equation is obtained as follows:
\begin{equation}
-\zeta_2<2\mathrm{e}^t-2<\zeta_2 \label{4.4.8}
\end{equation}
The above equation is transformed as an inequality of the variable $t$. The following equation is obtained.
\begin{equation}
\mathrm{ln}(\frac{2-\zeta_2}{2})<t<\mathrm{ln}(\frac{2+\zeta_2}{2}) \label{4.4.9}
\end{equation}
Because the constant $\zeta_2$ is defined by the Eq. (\ref{2.5}), the following equation is obtained by a numerical value.
\begin{equation}
-2.617 \cdots<t<0.655 \cdots \label{4.4.10}
\end{equation}
Therefore, one or two limits depend on the variables $B$ and $\phi$. 
For $A\geqq0$, the exact solution $x(t)$ of the type (X\hspace{-.1em}I) passes through zero and increases monotonously. For $A\leqq0$, the exact solution $x(t)$ of the type (X\hspace{-.1em}I) passes through zero and decreases monotonously. The limits $t =-2.617 \cdots$ and $t=0.6559 \cdots$ do not vary even if the parameter $A$ varies. In case that the curves vary with the amplitude $\omega$ under the conditions $A=1$, $\phi=-2$ and $B=2$, the curves obtained by the (X\hspace{-.1em}I) solution are shown in Fig.15. 
The possible range of the variable $t$ in the above equation is as follows:
\begin{equation}
\frac{1}{\omega} \mathrm{ln}(\frac{2-\zeta_2}{2})<t<\frac{1}{\omega} \mathrm{ln}(\frac{2+\zeta_2}{2}) \label{4.4.11}
\end{equation}
As the absolute value of $\omega$ increases in the above equation, the limit of the variable $t$ approaches zero, and the domain that the variable $t$ can take is narrowed. As the absolute value of $\omega$ decreases, the limit of the variable $t$ goes away from $0$. Therefore, the domain of the variable $t$ widens.

%
\begin{figure*}[tb]
\begin{center}
\includegraphics[width=0.75 \textwidth]{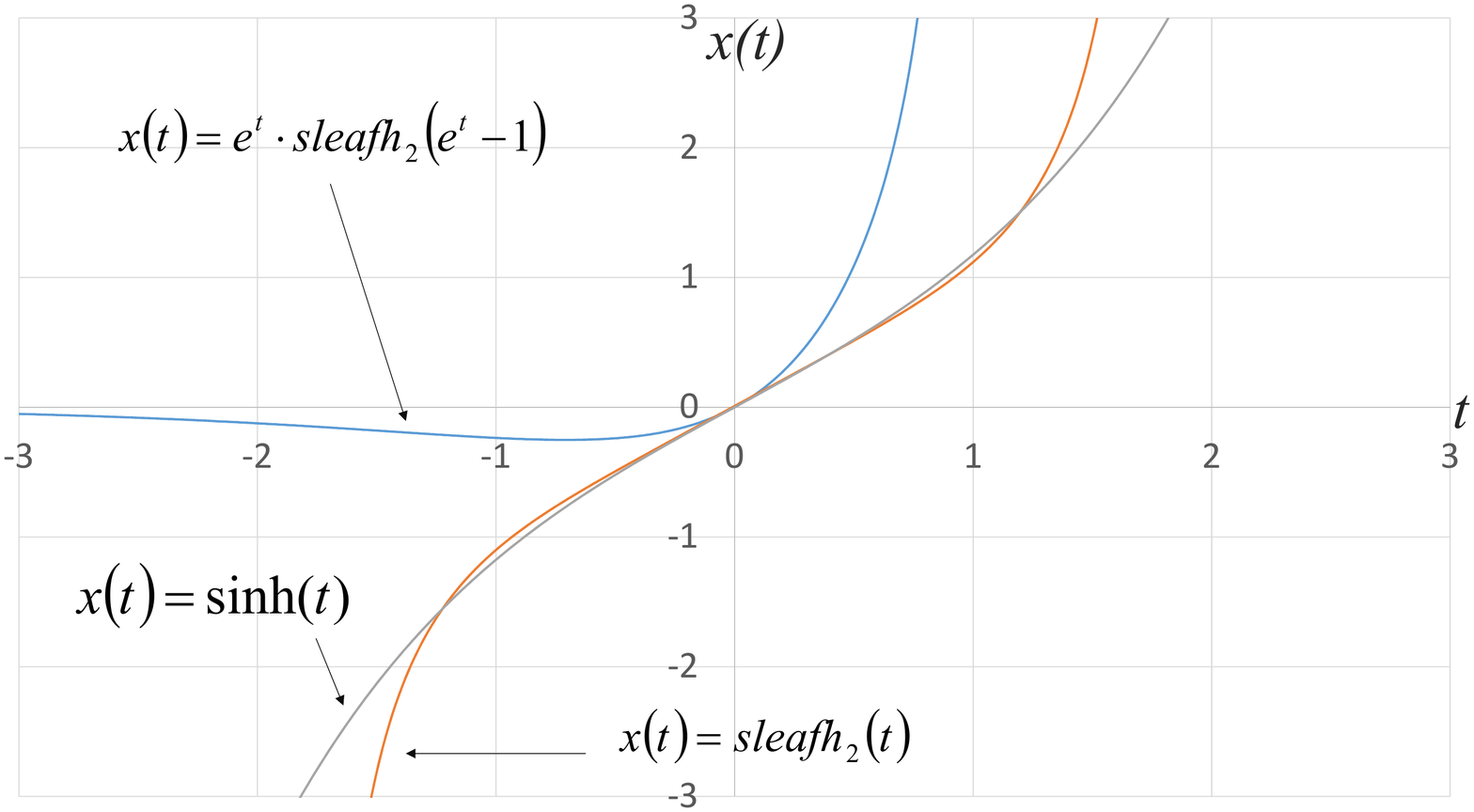}
\end{center}
\caption{ Curves obtained by the (X\hspace{-.1em}I) exact solution, the hyperbolic function $\mathrm{sinh}(t)$ and the leaf function $\mathrm{sleafh}_2(t)$ }
\label{fig:11}       
\end{figure*}
%
\begin{figure*}[tb]
\begin{center}
\includegraphics[width=0.75 \textwidth]{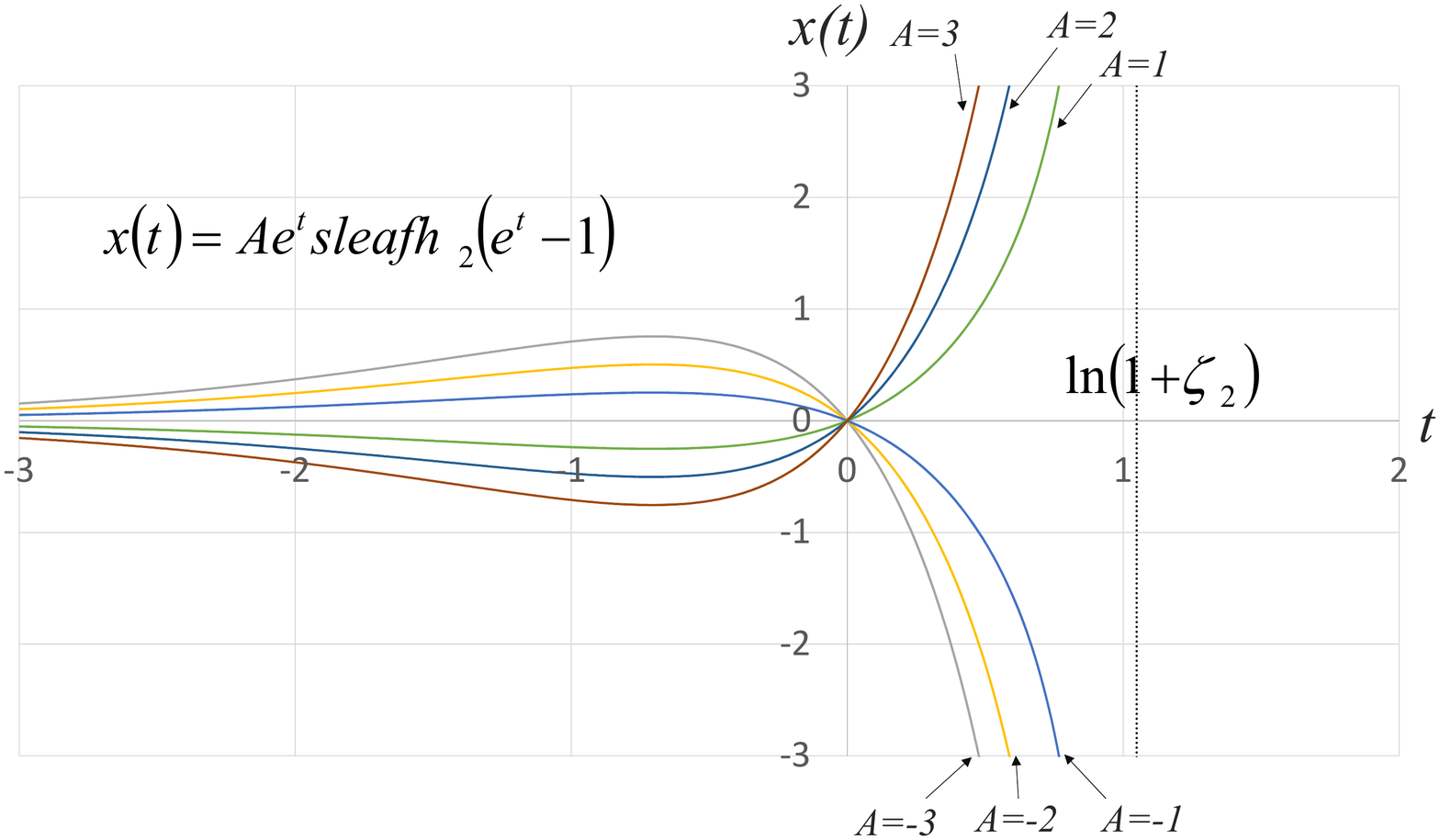}
\end{center}
\caption{ Curves obtained by the (X\hspace{-.1em}I) exact solution as a function of variations in the parameter $A$ 
($A=\pm1$, $\pm2$, $\pm3$) (Set $B=1$, $\omega=1$ and $\phi=-1$) }
\label{fig:12}       
\end{figure*}
%
\begin{figure*}[tb]
\begin{center}
\includegraphics[width=0.75 \textwidth]{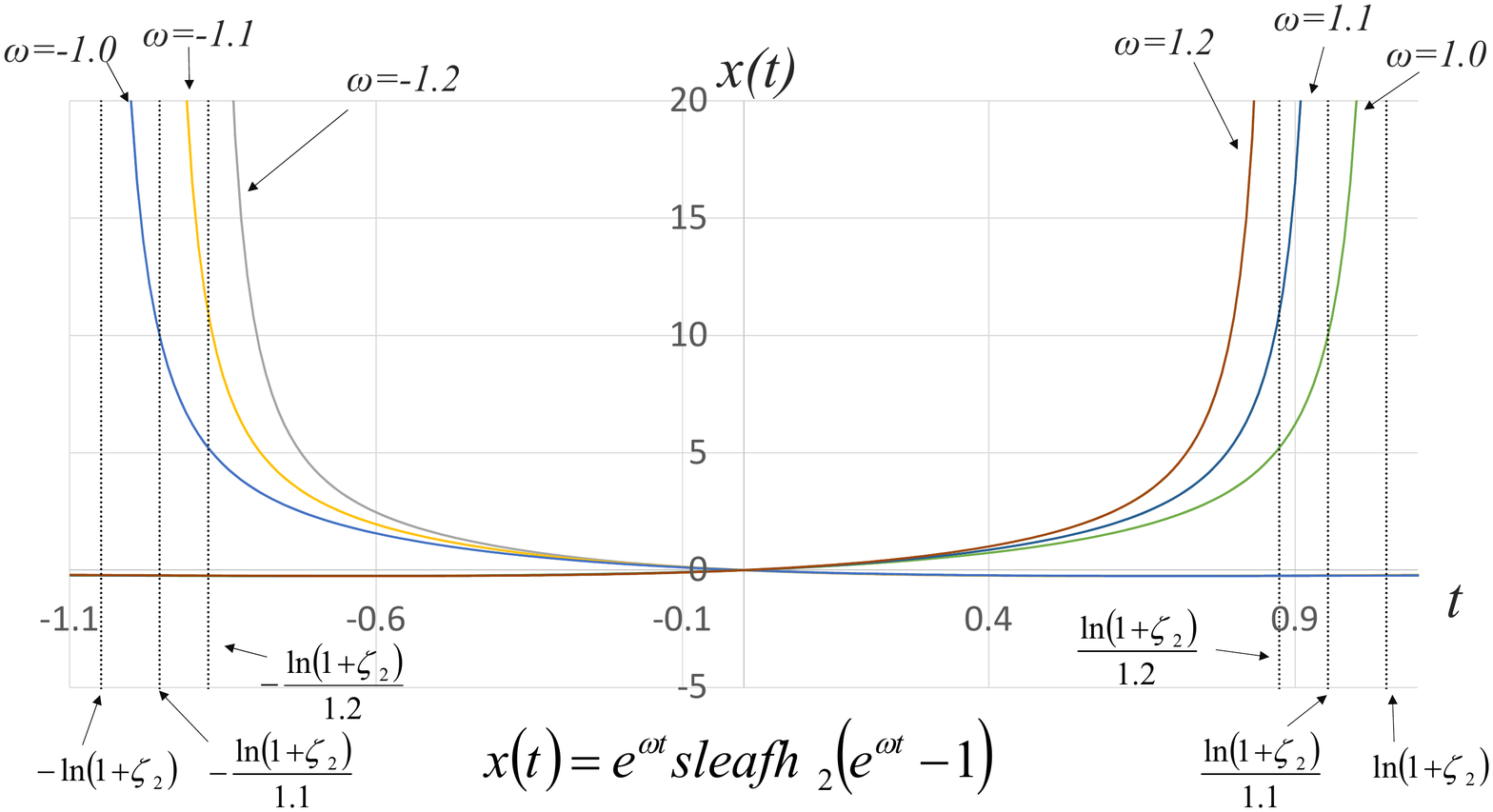}
\end{center}
\caption{ Curves obtained by the (X\hspace{-.1em}I) exact solution as a function of variations in the parameter $\omega$ ($\omega=\pm1.1$, $\pm1.2$, $\pm1.3$)(Set $B=1$, $A=1$ and $\phi=-1$) }
\label{fig:13}       
\end{figure*}
%
\begin{figure*}[tb]
\begin{center}
\includegraphics[width=0.75 \textwidth]{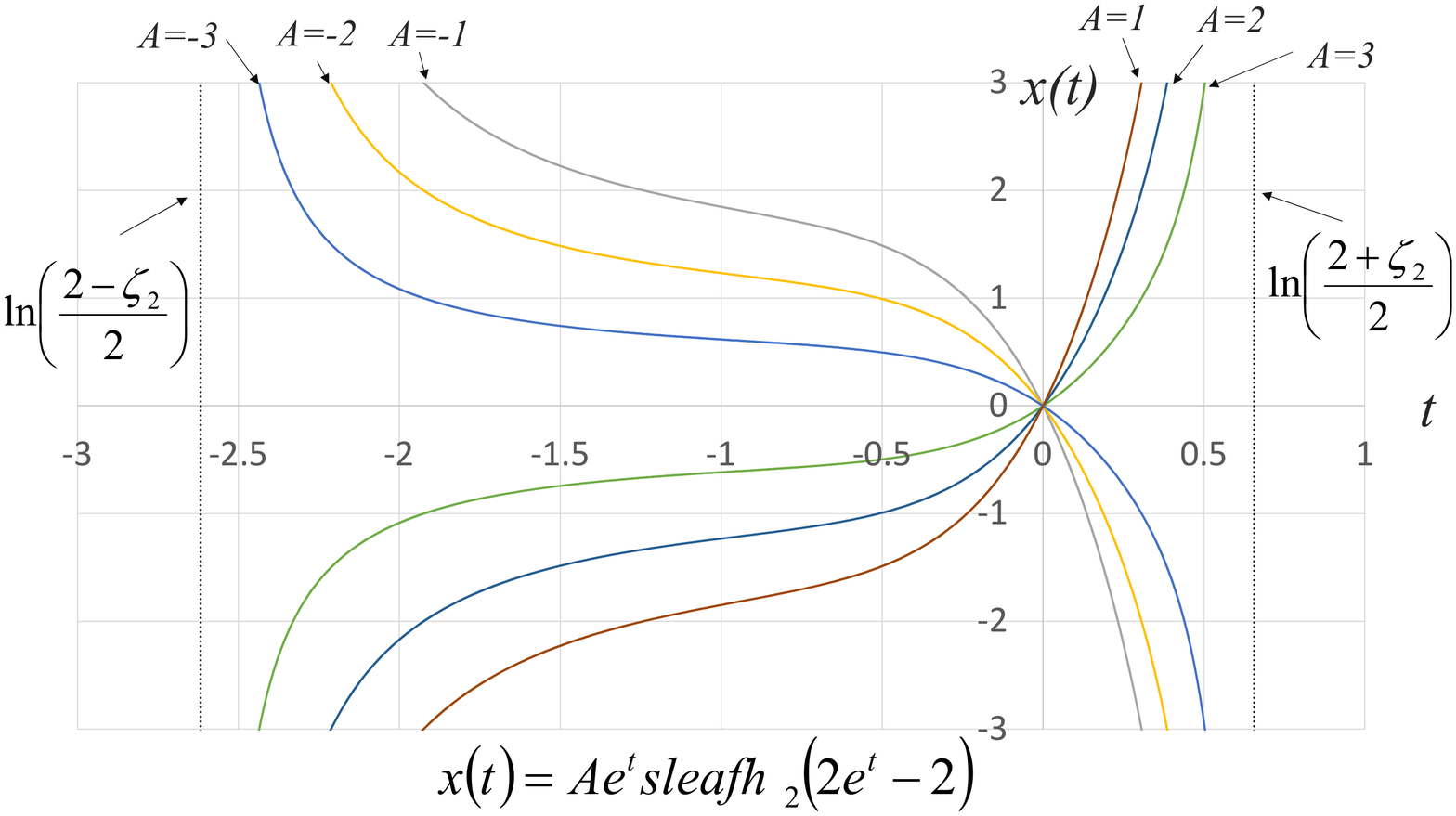}
\end{center}
\caption{ Curves obtained by the (X\hspace{-.1em}I) exact solution as a function of variations in the parameter $A$ ($A=\pm1$, $\pm2$, $\pm3$) (Set $B=2$, $\omega=1$ and $\phi=-2$) }
\label{fig:14}       
\end{figure*}
%
\begin{figure*}[tb]
\begin{center}
\includegraphics[width=0.75 \textwidth]{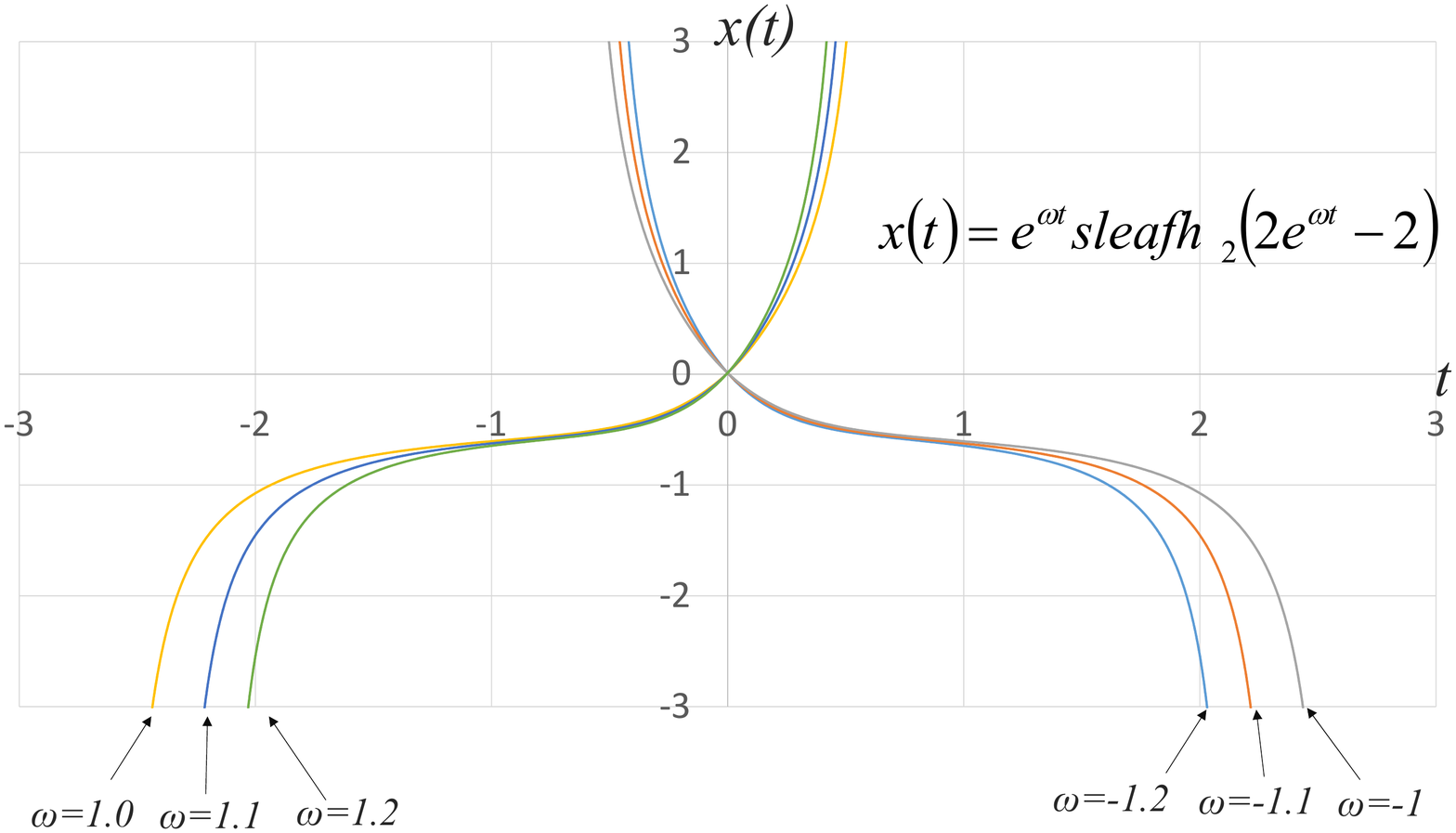}
\end{center}
\caption{ Curves obtained by the (X\hspace{-.1em}I) exact solution as a function of variations in the parameter $\omega$ ($\omega=\pm1.1$, $\pm1.2$, $\pm1.3$) (Set $B=2$, $A=1$ and $\phi=-2$) }
\label{fig:15}       
\end{figure*}

\begin{table}
\caption{ Numerical data for the (X\hspace{-.1em}I) solution under the conditions $A=1$, $B=1$, $\phi=-1$ and $\omega$=1}
\label{tab:5} 
\centering
\begin{tabular}{cccc}
\hline\noalign{\smallskip}
$t$ & $x(t)$ & $x(t)^3$  & $\frac{\mathrm{d}^2x(t)}{\mathrm{d}t^2}$ \\
\noalign{\smallskip}\hline\noalign{\smallskip}
-1.3	&	-0.203943407 $\cdots$	&	-0.0084826$\cdots$	&	0.034470582	$\cdots$\\
-1.2	&	-0.215597164	$\cdots$&	-0.010021417$\cdots$	&	0.070159104$\cdots$	\\
-1.1	&	-0.226540444	$\cdots$&	-0.011626185$\cdots$	&	0.116551151$\cdots$	\\
-1.0	&	-0.236307021	$\cdots$&	-0.013195622$\cdots$	&	0.176516436$\cdots$	\\
-0.9	&	-0.244294135	$\cdots$&	-0.014579382$\cdots$	&	0.253778736$\cdots$	\\
-0.8	&	-0.249725215	$\cdots$&	-0.015573535$\cdots$	&	0.353055955$\cdots$	\\
-0.7	&	-0.251602503	$\cdots$&	-0.0159274	$\cdots$&	0.480320133	$\cdots$\\
-0.6	&	-0.248647118	$\cdots$&	-0.015372705$\cdots$	&	0.643034898$\cdots$	\\
-0.5	&	-0.239224366	$\cdots$&	-0.013690403$\cdots$	&	0.850315776$\cdots$	\\
-0.4	&	-0.221252392	$\cdots$&	-0.010830884$\cdots$	&	1.113050698$\cdots$	\\
-0.3	&	-0.192093247	$\cdots$&	-0.007088205$\cdots$	&	1.443881365$\cdots$	\\
-0.2	&	-0.148426715	$\cdots$&	-0.003269913$\cdots$	&	1.857080067$\cdots$	\\
-0.1	&	-0.086107351	$\cdots$&	-0.000638441$\cdots$	&	2.368914385$\cdots$	\\
0	&	0	&	0	&	3	\\
0.1	&	0.116233275	$\cdots$&	0.001570332	$\cdots$&	3.783778377	$\cdots$\\
0.2	&	0.270486954	$\cdots$&	0.019789689	$\cdots$&	4.790941844	$\cdots$\\
0.3	&	0.472968457	$\cdots$&	0.105802647	$\cdots$&	6.191958894	$\cdots$\\
0.4	&	0.738030422	$\cdots$&	0.401996982	$\cdots$&	8.415748922	$\cdots$\\
0.5	&	1.088787044	$\cdots$&	1.290710471	$\cdots$&	12.56683153	$\cdots$\\
0.6	&	1.569138605	$\cdots$&	3.863526737	$\cdots$&	21.69657816	$\cdots$\\
0.7	&	2.277895182	$\cdots$&	11.81955724	$\cdots$&	45.67358469	$\cdots$\\
0.8	&	3.485728688	$\cdots$&	42.3526649	$\cdots$&	127.5548632$\cdots$	\\
0.9	&	6.220097075	$\cdots$&	240.6531152$\cdots$	&	605.0071132$\cdots$	\\
1.0	&	20.01721695	$\cdots$&	8020.678131$\cdots$	&	17263.60074$\cdots$	\\

\noalign{\smallskip}\hline
\end{tabular}
\end{table}

\subsection{Negative damping solution of the type (X\hspace{-.1em}I\hspace{-.1em}I)}
\label{Divergence solution12}
\begin{figure*}[tb]
\begin{center}
\includegraphics[width=0.75 \textwidth]{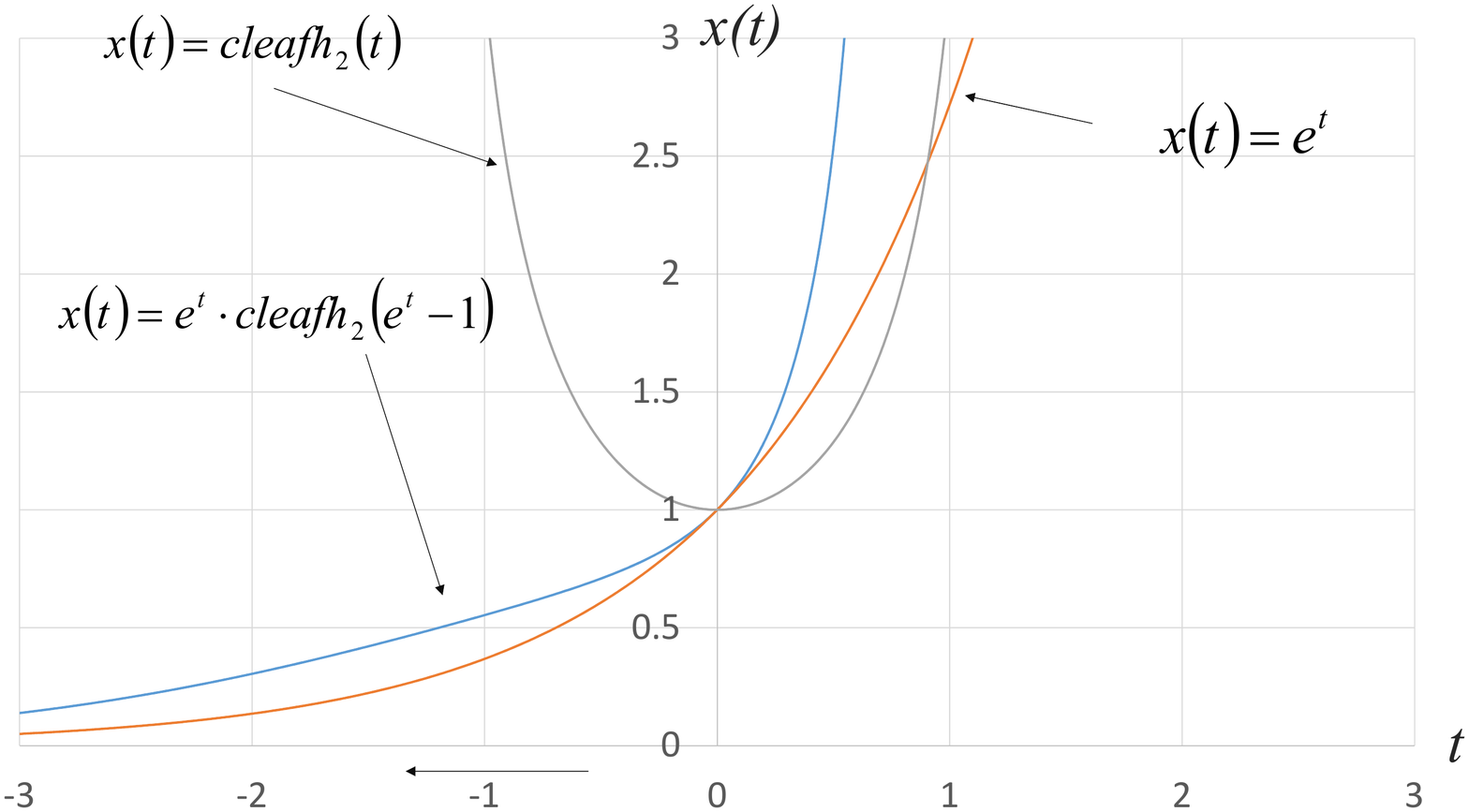}
\end{center}
\caption{ Curves obtained by the (X\hspace{-.1em}I\hspace{-.1em}I) exact solution, the exponential function $e^t$ and the leaf function $\mathrm{cleafh}_2(t)$ ($B=1$, $A=1$, $\omega=1$ and $\phi=-1$) }
\label{fig:16}       
\end{figure*}

\begin{figure*}[tb]
\begin{center}
\includegraphics[width=0.75 \textwidth]{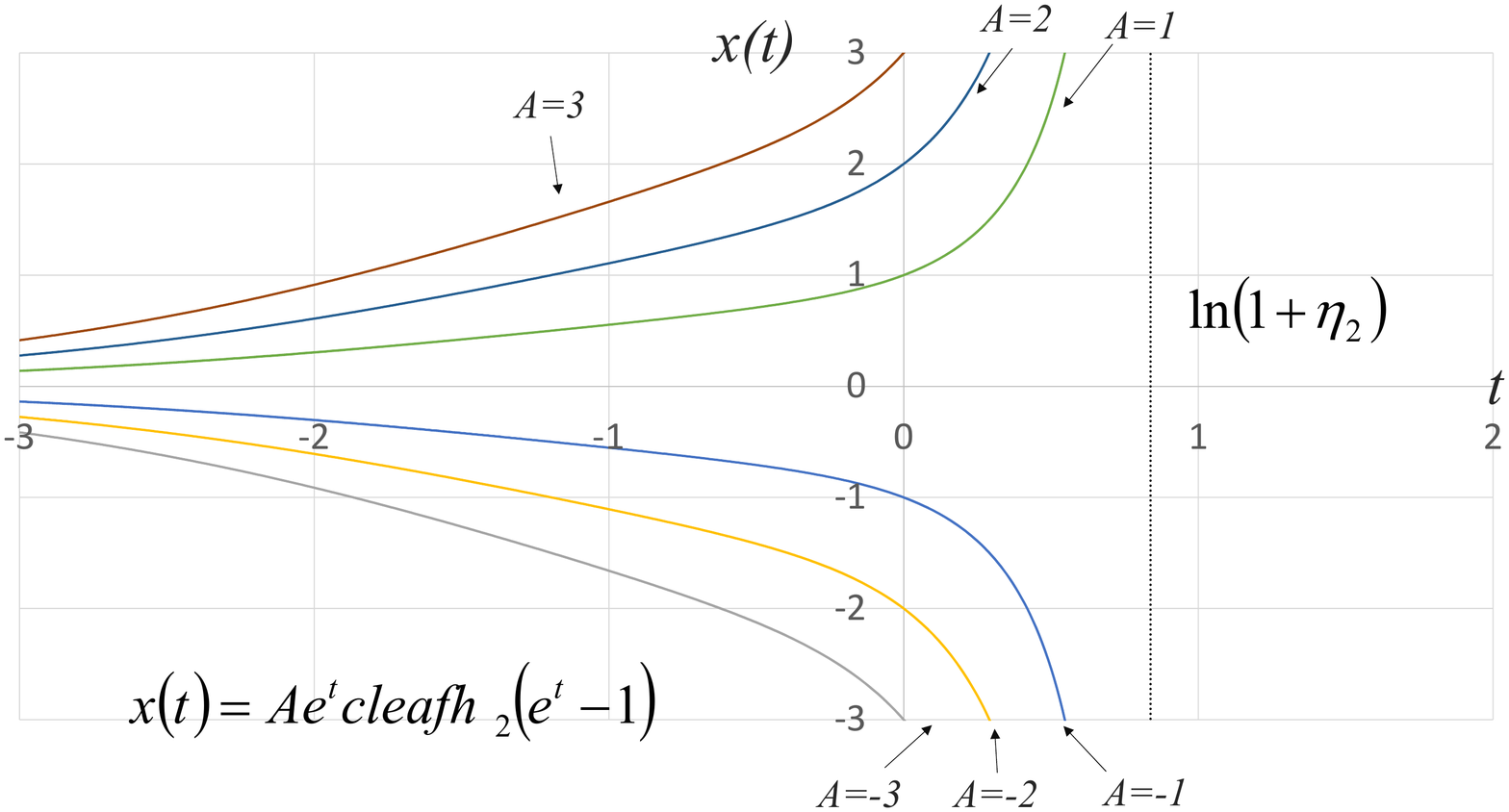}
\end{center}
\caption{ Curves obtained by the (X\hspace{-.1em}I\hspace{-.1em}I) exact solution as a function of variations in the parameter $A$ ($A=\pm1$, $\pm2$, $\pm3$) (Set $B=1$, $\omega=1$ and $\phi=-1$)}
\label{fig:17}       
\end{figure*}

\begin{figure*}[tb]
\begin{center}
\includegraphics[width=0.75 \textwidth]{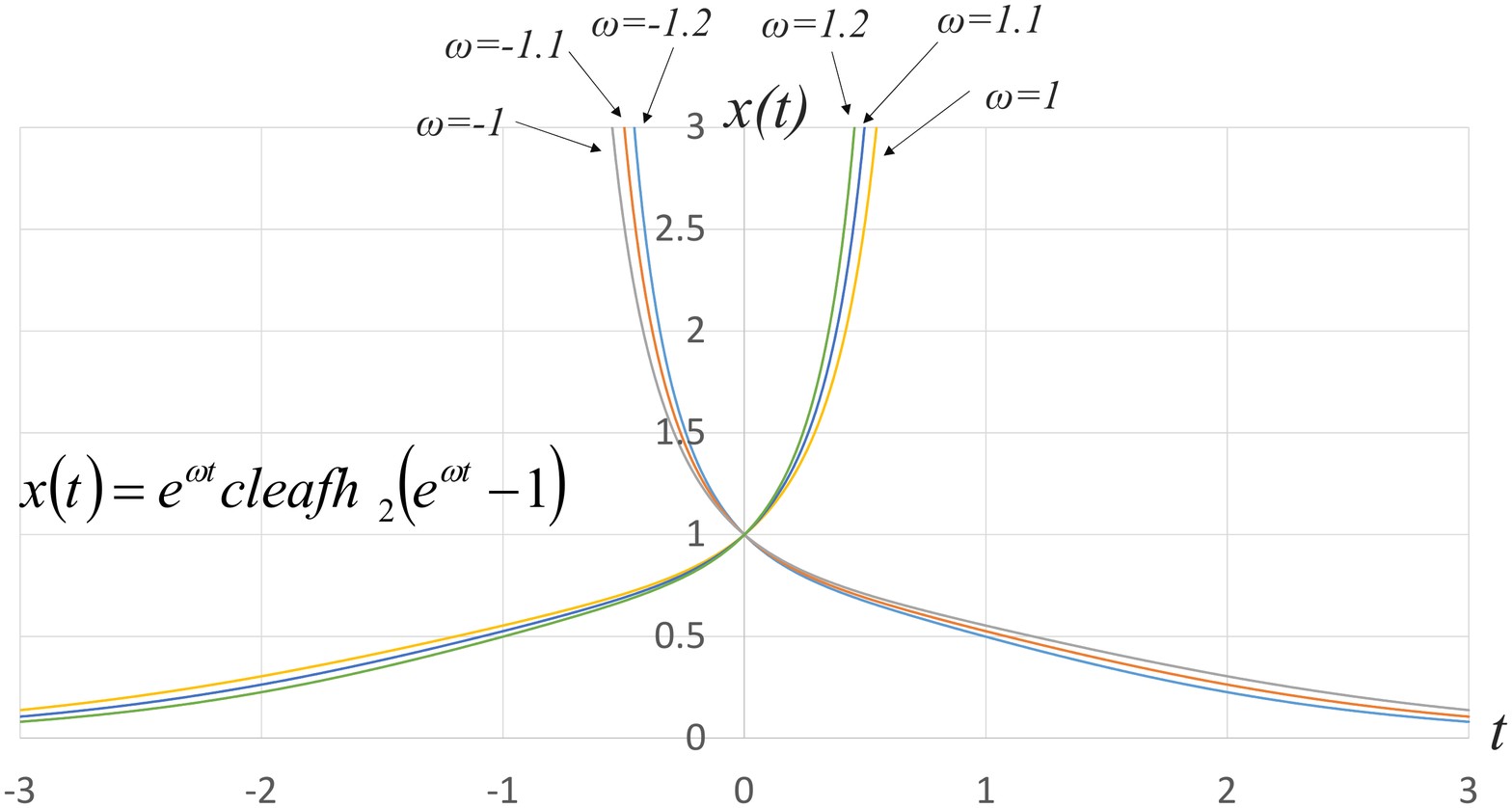}
\end{center}
\caption{ Curves obtained by the (X\hspace{-.1em}I\hspace{-.1em}I) exact solution as a function of variations in the parameter $\omega$ ($\omega=\pm1.1$, $\pm1.2$, $\pm1.3$) (Set $B=1$, $A=1$ and $\phi=-1$) }
\label{fig:18}       
\end{figure*}

\begin{figure*}[tb]
\begin{center}
\includegraphics[width=0.75 \textwidth]{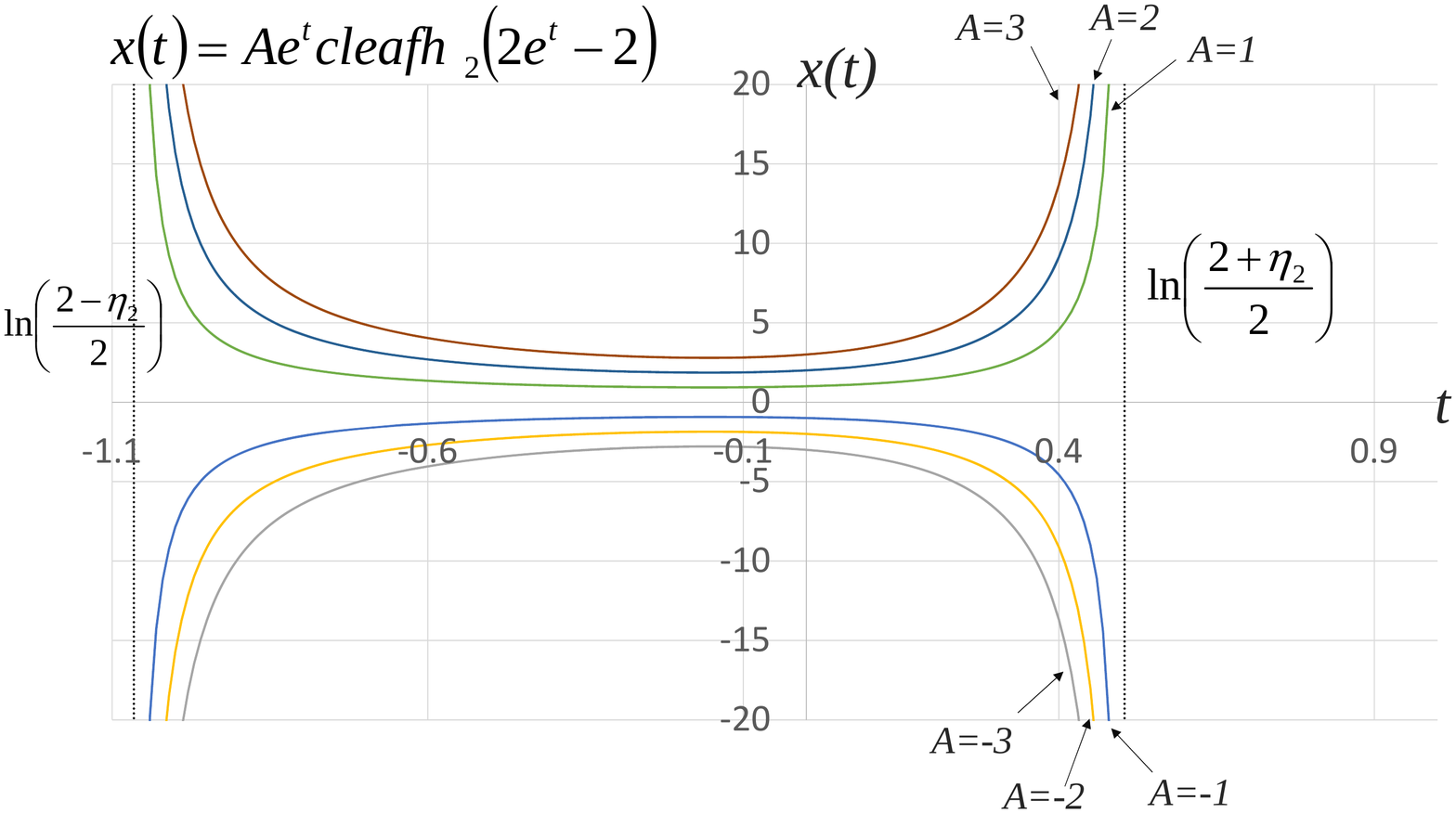}
\end{center}
\caption{ Curves obtained by the (X\hspace{-.1em}I\hspace{-.1em}I) exact solution as a function of variations in the parameter $A$  ($A=\pm1$, $\pm2$, $\pm3$) (Set $B=2$, $\omega=1$ and $\phi=-2$) }

\label{fig:19}       
\end{figure*}

Using the hyperbolic leaf function, the negative damping solution of the type (X\hspace{-.1em}I\hspace{-.1em}I) is presented in this section. In the (X\hspace{-.1em}I\hspace{-.1em}I) exact solution, the damped term for the Duffing equation $\mathrm{d}x(t)/\mathrm{d}t$ does not work in the direction of the suppressing vibration. Let us consider when $A =1$, $\omega=1$, $\phi=-1$ and $B=1$ in Eq. (\ref{3.5.1}), the exact solution of the type (X\hspace{-.1em}I\hspace{-.1em}I) is shown in Fig. \ref{fig:16}. 
For the inequality $t <0$, the variable $x(t)$ converges asymptotically to zero. The curve passes through $x(0)=1$. In the domain of the inequality $t>0$, the variable $x(t)$ diverges to infinity. The hyperbolic leaf function $\mathrm{cleafh}_2(\mathrm{e}^t-1)$ has a limit. Therefore, the range that the variable $\mathrm{e}^t-1$ can take is as follows.
\begin{equation}
-\eta_2<\mathrm{e}^t-1<\eta_2 \label{4.5.1}
\end{equation}
The following equation is obtained from the above equation.
\begin{equation}
-\eta_2+1<\mathrm{e}^t<\eta_2+1 \label{4.5.2}
\end{equation}
The constants $\eta_2$ is obtained by the Eq. (\ref{2.6}) (See Ref. \cite{Kaz_ch}).
In Eq. (\ref{4.5.2}), the inequality $ - \eta_2+1=-0.311\cdots<\mathrm{e}^t$ always holds for arbitrary time $t$. In the inequality $\mathrm{e}^t<\eta_2+1$ in Eq. (\ref{4.5.2}), we can take the logarithm of both sides. Solving for the variable $t$ yields the following equation.
\begin{equation}
t<\mathrm{ln}(\eta_2+1)=\mathrm{ln}(2.31102 \cdots)=0.837693 \cdots \label{4.5.3}
\end{equation}
The above inequality represents the domain of the variable $t$. In case that the curves vary with the amplitude $A$ under the conditions $\omega=1$, $\phi=-1$ and $B=1$, the curves obtained by the (X\hspace{-.1em}I\hspace{-.1em}I) solution are  shown in Fig. \ref{fig:17}. The limits do not vary even if the parameter $A$ varies.
For the inequality $A> 0$, the exact solution $x(t)$ for (X\hspace{-.1em}I\hspace{-.1em}I) diverges to positive infinity. In case that the curves vary with the phase $\omega$ under the conditions $A=1$, $\phi=-1$ and $B=1$, the curves obtained by the (X\hspace{-.1em}I\hspace{-.1em}I) solution are shown in Fig.\ref{fig:18}. 
Let us consider the limit of the exact solution for (X\hspace{-.1em}I\hspace{-.1em}I) when the parameter $\omega$ varies. The limit value for the variable $t$ is obtained using the following equation.
\begin{equation}
\mathrm{e}^{\omega t}-1=\pm\eta_2=\pm 1.31102 \cdots  \label{4.5.4}
\end{equation}
When the right side of Eq. (\ref{4.5.4}) is $1.31102 \cdots$, solving the above equation for the variable $t$ yields the following equation.
\begin{equation}
t=\frac{1}{\omega} \mathrm{ln}(\eta_2+1)=\frac{1}{\omega} \mathrm{ln}(1.311102 \cdots+1)=\frac{0.363 \cdots }{\omega} \label{4.5.5}
\end{equation}
When $\omega> 0$ in Eq. (\ref{4.5.4}), the limit of the variable $t$ exists in the domain $t>0$. When $\omega<0$, the limit of the variable $t$ exists in the domain $t<0$. As the absolute value of $\omega$ increases in the above equation, the limit of the variable $t$ approaches zero. The domain of possible values of $t$ becomes narrow. As the absolute value of $\omega$ decreases, the limit of the variable $t$ goes away from $0$. Therefore, the domain of the variable $t$ widens. When the right side of equation (\ref{4.5.4}) is -$1.31102 \cdots$, solving the above equation for the variable $t$ yields the following equation.
\begin{equation}
t=\frac{1}{\omega} \mathrm{ln}(-\eta_2+1)=\frac{1}{\omega} \mathrm{ln}(-1.311102 \cdots+1)=\frac{\mathrm{ln}(-0.311 \cdots) }{\omega} 
\label{4.5.6}
\end{equation}
Because the logarithm take the negative value $(-0.311 \cdots)$ in parenthesis, the variable $t$ that satisfies the above equation does not exist. In case that the curves vary with the phase $A$ under the conditions $\omega=1$, $\phi=-2$ and $B=2$, the curves obtained by the (X\hspace{-.1em}I\hspace{-.1em}I) solution are shown in Fig. \ref{fig:19}. 
The possible range of  $2\mathrm{e}^t-2$ in the above equation is as follows:
\begin{equation}
-\eta_2<2\mathrm{e}^t-2<\eta_2 \label{4.5.7}
\end{equation}
The above equation is transformed as an inequality of the variable $t$ to yield the following equation.
\begin{equation}
\mathrm{ln}(\frac{2-\eta_2}{2})<t<\mathrm{ln}(\frac{2+\eta_2}{2}) \label{4.5.8}
\end{equation}
Because the constant $\eta_2$ is defined by Eq. (\ref{2.6}), the following equation is obtained.
\begin{equation}
-1.065690 \cdots<t<0.504109 \cdots \label{4.5.9}
\end{equation}
One or two limits depend on the variables $B$ and $\phi$. The exact solution $x(t)$ for (X\hspace{-.1em}I\hspace{-.1em}I) passes through $x(0)=A$ . In a situation where $A\geqq 0$, the exact solution diverges to plus infinity near both $t=-1.065690 \cdots$ and $t=0.504109 \cdots$. In a case where $A\leqq0$, the exact solution diverges to minus infinity near both $t=-1.065690 \cdots$ and $t=0.504109 \cdots$. These limits do not vary even if the parameter $A$ varies. In case that the curves vary with the amplitude $\omega$ under the conditions $A=1$, $\phi=-2$ and $B=2$, the curves obtained by the (X\hspace{-.1em}I\hspace{-.1em}I) solution are shown in Fig. 20. 
The possible range of the variable $t$ in the above equation is as follows:
\begin{equation}
\frac{1}{\omega} \mathrm{ln}(\frac{2-\eta_2}{2})<t<\frac{1}{\omega} \mathrm{ln}(\frac{2+\eta_2}{2}) \label{4.5.10}
\end{equation}
As the absolute value of $\omega$ increases in the above equation, the domain of possible values for $t$ becomes narrower and narrower as shown in the Eq. (\ref{4.5.10}). As the absolute value of $\omega$ decreases, the limit of the variable $t$ goes away from $0$. Therefore, the domain of the variable $t$ becomes wider and wider.

%
%
%
%
%
\begin{figure*}[tb]
\begin{center}
\includegraphics[width=0.75 \textwidth]{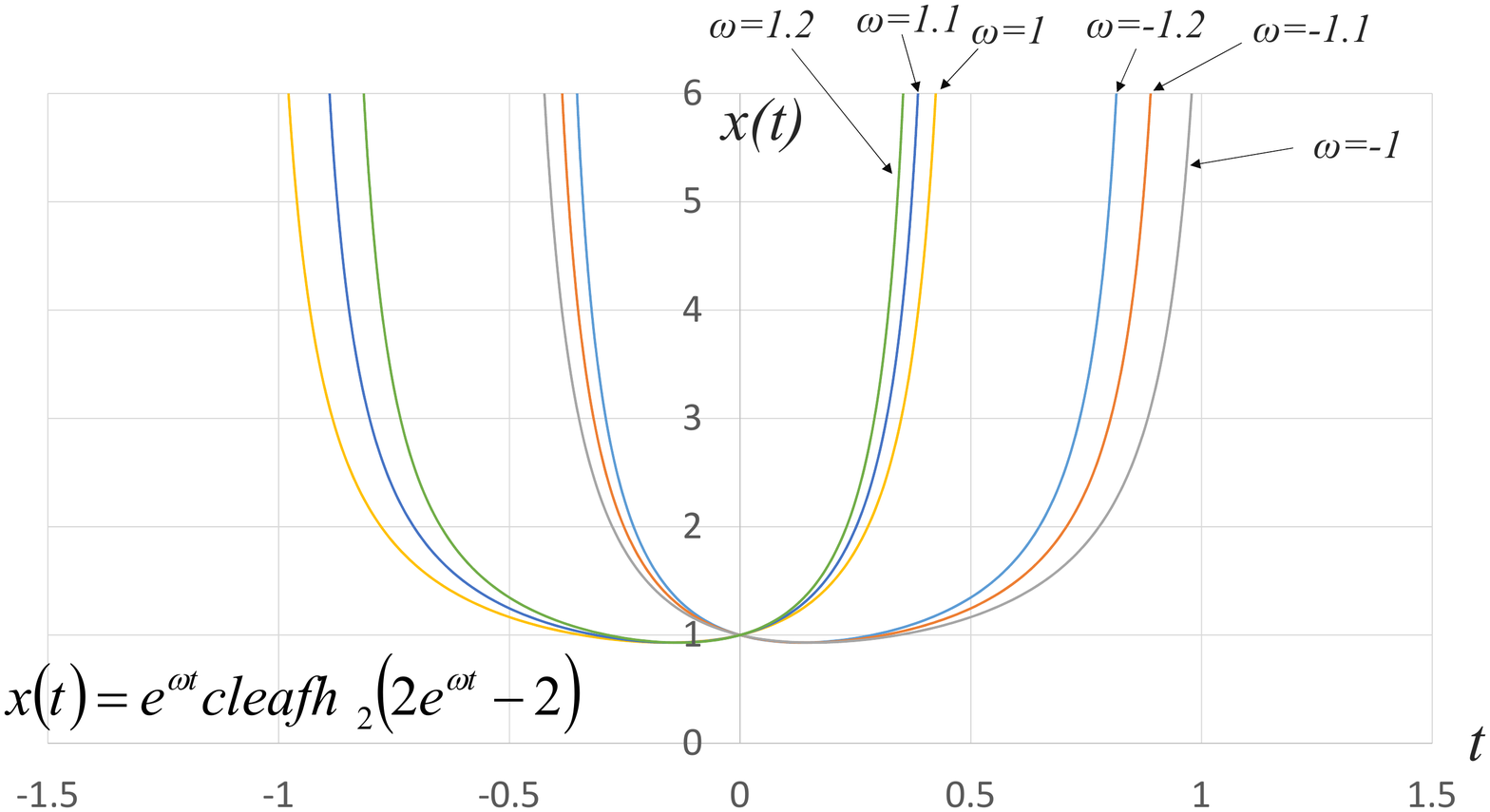}
\end{center}
\caption{ Curves obtained by the (X\hspace{-.1em}I\hspace{-.1em}I) exact solution as a function of variations in the parameter $\omega$ ($\omega=\pm1.1$, $\pm1.2$, $\pm1.3$) (Set $B=2$, $A=1$ and $\phi=-2$) }
\label{fig:20}       
\end{figure*}

\begin{table}
\caption{ Numerical data for the (X\hspace{-.1em}I\hspace{-.1em}I) solution under the conditions $A=1$, $B=1$, $\phi=-1$ and $\omega=1$}
\label{tab:6}    
\centering
\begin{tabular}{cccc}
\hline\noalign{\smallskip}
$t$ & $x(t)$ & $x(t)^3$  & $\frac{\mathrm{d}^2x(t)}{\mathrm{d}t^2}$ \\
\noalign{\smallskip}\hline\noalign{\smallskip}
-1.3	&	0.472441888 $\cdots$	&	0.105449661$\cdots$	&	0.051898217$\cdots$	\\
-1.2	&	0.498895941 $\cdots$	&	0.124173783$\cdots$	&	0.051915902$\cdots$	\\
-1.1	&	0.525873021$\cdots$	&	0.145426205$\cdots$	&	0.056603345$\cdots$	\\
-1.0	&	0.553421826$\cdots$	&	0.169499667$\cdots$	&	0.068162077$\cdots$	\\
-0.9	&	0.581660476$\cdots$	&	0.196792554$\cdots$	&	0.089646632$\cdots$	\\
-0.8	&	0.610807381$\cdots$	&	0.227883473$\cdots$	&	0.125364194$\cdots$	\\
-0.7	&	0.641224819$\cdots$	&	0.26365194	$\cdots$&	0.181487155$\cdots$	\\
-0.6	&	0.673481459$\cdots$	&	0.305475886$\cdots$	&	0.267021052$\cdots$	\\
-0.5	&	0.708443676$\cdots$	&	0.355562526$\cdots$	&	0.395339264$\cdots$	\\
-0.4	&	0.747411388$\cdots$	&	0.417521777$\cdots$	&	0.586723946$\cdots$	\\
-0.3	&	0.792324421$\cdots$	&	0.49740383	$\cdots$&	0.872693767$\cdots$	\\
-0.2	&	0.846083923$\cdots$	&	0.605675949$\cdots$	&	1.303697695$\cdots$	\\
-0.1	&	0.913068861$\cdots$	&	0.761220712$\cdots$	&	1.963756376$\cdots$	\\
0	&	1	&	1	&	3	\\
0.1	&	1.117463203$\cdots$	&	1.395403132	$\cdots$&	4.687637715$\cdots$	\\
0.2	&	1.28278683	$\cdots$&	2.110879669	$\cdots$&	7.587902875$\cdots$	\\
0.3	&	1.525994729$\cdots$	&	3.553522755	$\cdots$&	12.98313508$\cdots$	\\
0.4	&	1.903703402$\cdots$	&	6.899186071	$\cdots$&	24.28272498$\cdots$	\\
0.5	&	2.537412837$\cdots$	&	16.33704096	$\cdots$&	52.72116914$\cdots$	\\
0.6	&	3.748216016$\cdots$	&	52.65914896	$\cdots$&	150.0200518$\cdots$	\\
0.7	&	6.779310073$\cdots$	&	311.5706171	$\cdots$&	767.2598046$\cdots$	\\
0.8	&	26.03361091$\cdots$	&	17644.25107	$\cdots$&	37347.95494$\cdots$	\\

\noalign{\smallskip}\hline
\end{tabular}
\end{table}

\subsection{Damping solution of the type (X\hspace{-.1em}I\hspace{-.1em}I\hspace{-.1em}I)}
\label{4.6}
Let us consider Eq. (\ref{3.6.1}) that satisfies the Duffing equation. Physically, the exact solution of the type (X\hspace{-.1em}I\hspace{-.1em}I\hspace{-.1em}I) corresponds to the nonlinear spring model of the damping system. 

\subsubsection{Period}
 When $A=1$, $\omega=1$, $\phi=0$ and $B=\pi_2/2$, the curves obtained by the (X\hspace{-.1em}I\hspace{-.1em}I\hspace{-.1em}I) solution are shown in Fig.21. The horizontal axis represents time $t$ whereas the vertical axis represents displacement, $x(t)$. As shown in Fig. 21, the period of the wave varies with time.  In the Fig.21, the $T_1$, $T_2$, $\cdots$, $T_6$ represents the first period, the second period, $\cdots$, and the sixth period, respectively. First, let us consider the first period $T_1$.  Substituting -$\infty$ into the variable  $t$, the value of the type (X\hspace{-.1em}I\hspace{-.1em}I\hspace{-.1em}I) can be obtained as follows:
\begin{equation}
x( -\infty)=\mathrm{e}^{- \infty} \cdot \mathrm{sleaf}_2(\frac{\pi_2}{2} \cdot \mathrm{e}^{-\infty})=0 \cdot \mathrm{sleaf}_2(\frac{\pi_2}{2} \cdot 0)=0
\end{equation} 
On the other hand, substituting $\mathrm{ln}4$ into the variable $t$, the value of the type (X\hspace{-.1em}I\hspace{-.1em}I\hspace{-.1em}I) can be obtained as follows:
\begin{equation}
x(\mathrm{ln}4)=\mathrm{e}^{\mathrm{ln}4} \cdot \mathrm{sleaf}_2(\frac{\pi_2}{2} \cdot \mathrm{e}^{\mathrm{ln}4})=4 \cdot \mathrm{sleaf}_2(2 \pi_2)=0
\end{equation} 
The first period $T_1$ is as follows:
\begin{equation}
T_1=\mathrm{ln}4-\infty=-\infty
\end{equation} 
The domain of the first period $T_1$ is as follows: 
\begin{equation}
-\infty<t \leqq \mathrm{ln}4
\end{equation} 
Next, let us consider the second period $T_2$.  Substituting $\mathrm{ln}8$ into the variable $t$, the value of the type (X\hspace{-.1em}I\hspace{-.1em}I\hspace{-.1em}I) can be obtained as follows:
\begin{equation}
x(\mathrm{ln}8)=e^{\mathrm{ln}8} \cdot \mathrm{sleaf}_2(\frac{\pi_2}{2} \cdot e^{\mathrm{ln}8})=8 \cdot \mathrm{sleaf}_2(4 \pi_2)=0
\end{equation} 
The second period $T_2$ is as follows:
\begin{equation}
T_2=\mathrm{ln}8-\mathrm{ln}4
\end{equation} 
The domain of the second period $T_2$ is as follows: 
\begin{equation}
\mathrm{ln}4 \leqq t \leqq \mathrm{ln}8
\end{equation} 
In this way, the $m^{th}$ period $T_m$ can be generalized as follows:
\begin{equation}
T_m=\mathrm{ln}(4m)-\mathrm{ln}(4m-4)=\mathrm{ln} \frac{m}{m-1} (m=1,2,3, \cdots)
\end{equation} 
The domain of the $m^{th}$ period $T_m$ is as follows: 
\begin{equation}
\mathrm{ln}(4m-4) \leqq t \leqq \mathrm{ln}(4m) (m=1,2,3, \cdots)
\end{equation} 

\subsubsection{Amplitude}
A time with respect to the convex upward of wave is obtained by using the gradient of the type (XIII). Using the condition $\mathrm{d}x(t)/\mathrm{d}t=0$ in the Eq. (\ref{A13.1}), we can obtain the equation as follows:
\begin{equation}
B^2 \cdot \mathrm{e}^{2 \omega t} (\mathrm{sleaf}_2(B \cdot \mathrm{e}^{\omega t}))^4
+(\mathrm{sleaf}_2(B \cdot \mathrm{e}^{\omega t}))^2-B^2 \cdot \mathrm{e}^{2 \omega t}=0 \label{4.6.2.1}
\end{equation} 
Unfortunately, we can not obtain the exact solution of the variable $t$ by the above equation. Using the numerical analysis, the variable $t$ that satisfy with the Eq. (\ref{4.6.2.1}) is obtained. Under the conditions $B=\pi_2/2, A=1, \omega=1$ and $\phi=0$, the amplitude $x(t)$ with respect to the time $t$ are summarized in the Table 7 and Table 8.

\begin{table}
\caption{ Exact amplitude $x(t)$ at convex upward for (X\hspace{-.1em}I\hspace{-.1em}I\hspace{-.1em}I) }
\label{tab:7}       
\centering
\begin{tabular}{ccc}
\hline\noalign{\smallskip}
Number in Fig. 21 & Time $t$ & Exact amplitude $x(t)$ \\
\noalign{\smallskip}\hline\noalign{\smallskip}
(1U)	&	0.210868709 $\cdots$	&	1.123036671 $\cdots$	\\
(2U)	&	1.620875816 $\cdots$	&	5.028840775 $\cdots$	\\
(3U)	&	2.200796752 $\cdots$	&	9.016117892 $\cdots$	\\
(4U)	&2.566666239 $\cdots$	&	13.0111741 $\cdots$	\\
(5U)	&	2.834218406 $\cdots$	&	17.00854947 $\cdots$	\\
(6U)	&	3.045181426 $\cdots$	&	21.0069228 $\cdots$	\\
\noalign{\smallskip}\hline
\end{tabular}
\end{table}

\begin{table}
\caption{  Exact amplitude $x(t)$ at convex downward for (X\hspace{-.1em}I\hspace{-.1em}I\hspace{-.1em}I) }
\label{tab:8}       
\centering
\begin{tabular}{ccc}
\hline\noalign{\smallskip}
Number in Fig. 21 & Time $t$ & Exact amplitude $x(t)$ \\
\noalign{\smallskip}\hline\noalign{\smallskip}
(1D)	&	1.129468643 $\cdots$	&	-3.047364138 $\cdots$	\\
(2D)	&	1.951794708 $\cdots$	&	-7.020686991 $\cdots$	\\
(3D)	&	2.40029079 $\cdots$	&	-11.01319903 $\cdots$	\\
(4D)	&2.709340595 $\cdots$	&	-15.00968732 $\cdots$	\\
(5D)	&	2.945243827 $\cdots$	&	-19.00765067 $\cdots$	\\
(6D)	&	3.136043671 $\cdots$	&	-23.00632133 $\cdots$	\\
\noalign{\smallskip}\hline
\end{tabular}
\end{table}

The curve $x(t)=\pm \mathrm{e}^{t}$ is added in Fig. 21. As shown in Fig. 21, the curve $\mathrm{e}^{t}$ intersects the curve $\mathrm{sleaf}_2(\frac{\pi_2}{2} \mathrm{e}^{t})$ at  the condition $\mathrm{sleaf}_2(\frac{\pi_2}{2} \mathrm{e}^{t})= \pm 1$. In the case where the condition $\mathrm{sleaf}_2(\frac{\pi_2}{2} \mathrm{e}^{t})=1$, the time $t$ is obtained as follows:
\begin{equation}
\frac{\pi_2}{2} \mathrm{e}^{t}=\frac{4k-3}{2} \pi_2 (k=1,2,3 \cdots)
\label{4.6.2.2}
\end{equation}
The following equation is obtained by the above equation.
\begin{equation}
t_k= \mathrm{ln}(4k-3) (k=1,2,3, \cdots)
\label{4.6.2.3}
\end{equation}
Substituting $t_k$ into the Eq. (\ref{3.6.1}) for the solution of type (X\hspace{-.1em}I\hspace{-.1em}I\hspace{-.1em}I), the following equation is obtained.
\begin{equation}
\begin{split}
x(t_k)&=\mathrm{exp}(t_k) \cdot \mathrm{sleaf}_2(\frac{\pi_2}{2} \cdot \mathrm{exp}(t_k))   \\
&=\mathrm{exp}(\mathrm{ln}(4k-3)) \cdot \mathrm{sleaf}_2(\frac{\pi_2}{2} \cdot \mathrm{exp}(\mathrm{ln}(4k-3))) \\
&=(4k-3) \mathrm{sleaf}_2(\frac{\pi_2}{2} \cdot (4k-3)) \\
&=4k-3 (k=1,2,3, \cdots)
\label{4.6.2.4}
\end{split}
\end{equation}
In the case where the condition $\mathrm{sleaf}_2(\frac{\pi_2}{2} \mathrm{e}^{t})=-1$, the time $t$ is obtained as follows:
\begin{equation}
\frac{\pi_2}{2} \mathrm{e}^{t}=\frac{4k-1}{2} \pi_2 (k=1,2,3 \cdots)
\label{4.6.2.5}
\end{equation}
The following equation is obtained by the above equation.
\begin{equation}
t_k= \mathrm{ln}(4k-1) (k=1,2,3, \cdots) \label{4.6.2.6}
\end{equation}
Substituting $t_k$ into the Eq. (\ref{3.6.1}), the following equation is obtained.
\begin{equation}
\begin{split}
x(t_k)&=\mathrm{exp}(t_k) \cdot \mathrm{sleaf}_2(\frac{\pi_2}{2} \cdot \mathrm{exp}( t_k))   \\
&=\mathrm{exp}(\mathrm{ln}(4k-1)) \cdot \mathrm{sleaf}_2(\frac{\pi_2}{2} \cdot \mathrm{exp}(\mathrm{ln}(4k-1))) \\
&=(4k-1) \mathrm{sleaf}_2(\frac{\pi_2}{2} \cdot (4k-1)) \\
&=-(4k-1) (k=1,2,3, \cdots)
\label{4.6.2.7}
\end{split}
\end{equation}

\begin{table}
\caption{ Approximate amplitude $x(t)$ at convex upward for (X\hspace{-.1em}I\hspace{-.1em}I\hspace{-.1em}I) }
\label{tab:9}       
\centering
\begin{tabular}{ccc}
\hline\noalign{\smallskip}
Number in Fig. 21 & Time $t$ by Eq. (\ref{4.6.2.3}) & Approximate amplitude $x(t)$ by Eq. (\ref{4.6.2.4}) \\
\noalign{\smallskip}\hline\noalign{\smallskip}
(1U)	&	0.0 (=ln(0))	&	1.0 \\
(2U)	&	1.609437912 $\cdots$	(=ln(5)) &	5.0 	\\
(3U)	&	2.197224577 $\cdots$	(=ln(9)) &	9.0 	\\
(4U)	&2.564949357 $\cdots$	(=ln(13)) &	13.0 	\\
(5U)	&	2.833213344 $\cdots$	(=ln(17)) &	17.0 	\\
(6U)	&	3.044522438 $\cdots$	(=ln(21)) &	21.0 	\\
\noalign{\smallskip}\hline
\end{tabular}
\end{table}

\begin{table}
\caption{ Approximate amplitude $x(t)$ at convex downward for (X\hspace{-.1em}I\hspace{-.1em}I\hspace{-.1em}I) }
\label{tab:10}       
\centering
\begin{tabular}{ccc}
\hline\noalign{\smallskip}
Number in Fig. 21 & Time $t$ by Eq. (\ref{4.6.2.6})& Approximate amplitude $x(t)$ by Eq. (\ref{4.6.2.7}) \\
\noalign{\smallskip}\hline\noalign{\smallskip}
(1D)	&1.098612289 $\cdots$ (=ln(3))	&	-3.0 \\
(2D)	&	1.945910149 $\cdots$	(=ln(7))&	-7.0 	\\
(3D)	&	2.397895273 $\cdots$	(=ln(11))&	-11.0 	\\
(4D)	&2.708050201 $\cdots$	(=ln(15))&	-15.0 	\\
(5D)	&	2.944438979 $\cdots$	(=ln(19))&	-19.0 	\\
(6D)	&	3.135494216 $\cdots$	(=ln(23))&	-23.0 	\\
\noalign{\smallskip}\hline
\end{tabular}
\end{table}

These data both the time $t$ (by the Eq. (\ref{4.6.2.3}) and the Eq. (\ref{4.6.2.6}) ) and the amplitude $x(t)$  (by the Eq. (\ref{4.6.2.4}) and the Eq. (\ref{4.6.2.7}))  are summarized in Table 9 and Table 10. Since the condition $\mathrm{sleaf}_2(\frac{\pi_2}{2} e^{t})= \pm 1$ is not satisfied with the Eq. (\ref{4.6.2.1}), the values $t$ obtained by the Eq. (\ref{4.6.2.3}) and the Eq. (\ref{4.6.2.6})  are not the extremal value in the type (XIII). However, comparing the numerical data in the table 7(or the table 8) with the numerical data in the table 9 (or the table 10), the approximate amplitude (the table 9 and the table 10) are almost agreed with the exact amplitude by numerical analysis (the table 7 and the table 8). As an approximate value of the amplitude, amplitude can be easily obtained from the Eq. (\ref{4.6.2.4}) and the Eq. (\ref{4.6.2.7}). Fig. 22 and Fig. 23 show the variation of waves as the parameter $A$ vary. As the parameter $A$ increases, the height of the wave increases with keeping the wave period. Let us consider the parameter $\omega$ under the conditions $B = \pi_2/2$, $\phi = 0$ and $A=1$. The curves for the (X\hspace{-.1em}I\hspace{-.1em}I\hspace{-.1em}I) solution are shown in the Fig. 24. In the negative domain of the variable $t$, the waves have no period. Conversely, in the positive domain of the variable $t$, the waves show periodicity. As the magnitude of $\omega$ increases, a period becomes shorter. As the magnitude of $\omega$ decreases, a period becomes longer. As shown in Fig 25, for $\omega<0$, the wave has no periodicity in the positive domain $t$. The exact solution for (X\hspace{-.1em}I\hspace{-.1em}I\hspace{-.1em}I) decreases monotonically as time, $t$ increases. In the negative domain of the variable $t$, the wave shows periodicity.

%
\begin{figure*}[tb]
\begin{center}
\includegraphics[width=0.75 \textwidth]{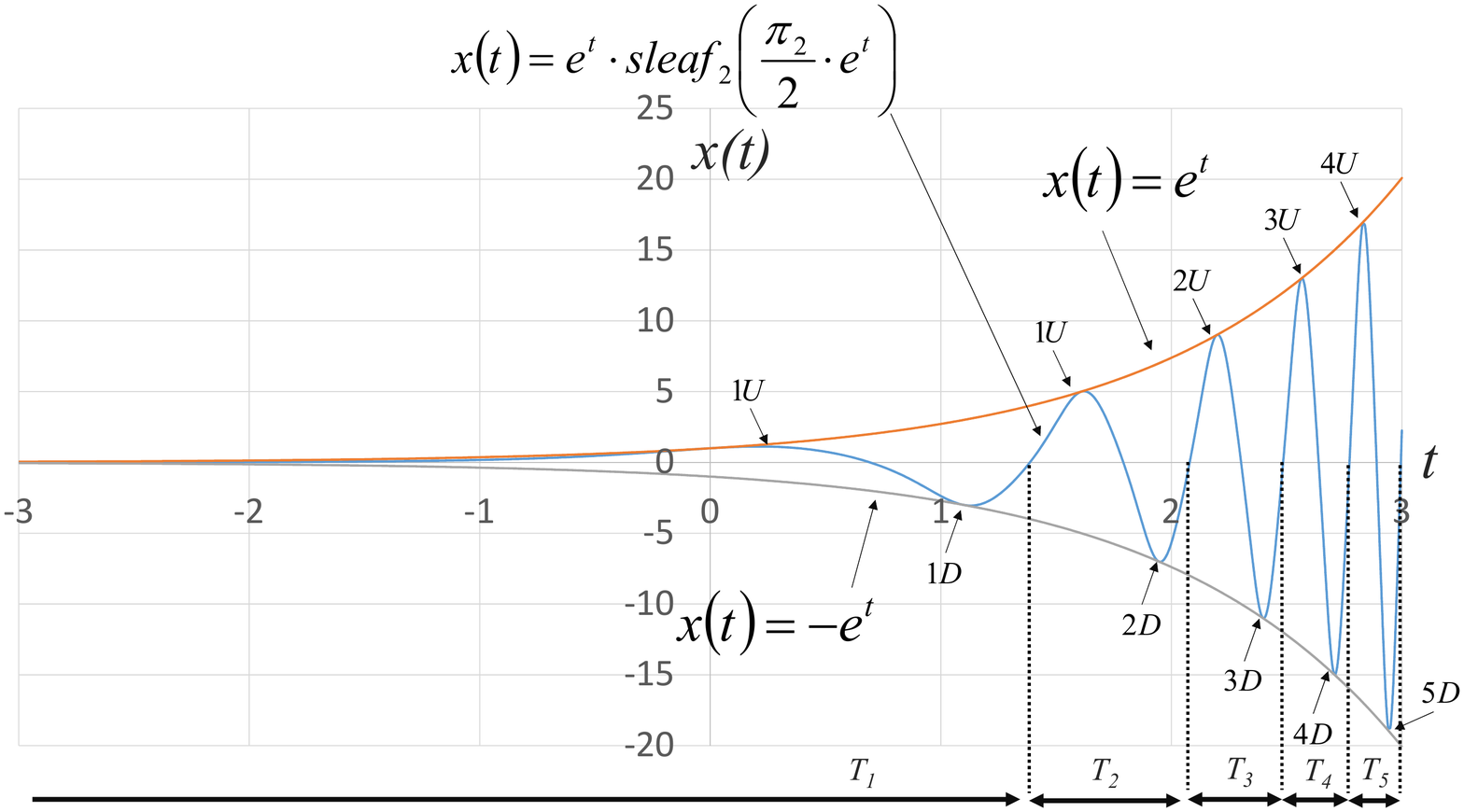}
\end{center}
\caption{ Waves obtained by the (X\hspace{-.1em}I\hspace{-.1em}I\hspace{-.1em}I) exact solution and the exponential function $\mathrm{e}^t$ (Set $B=\pi_2/2$, $A=1$, $\omega=1$ and $\phi=0$) }

\label{fig:21}       
\end{figure*}
%
\begin{figure*}[tb]
\begin{center}
\includegraphics[width=0.75 \textwidth]{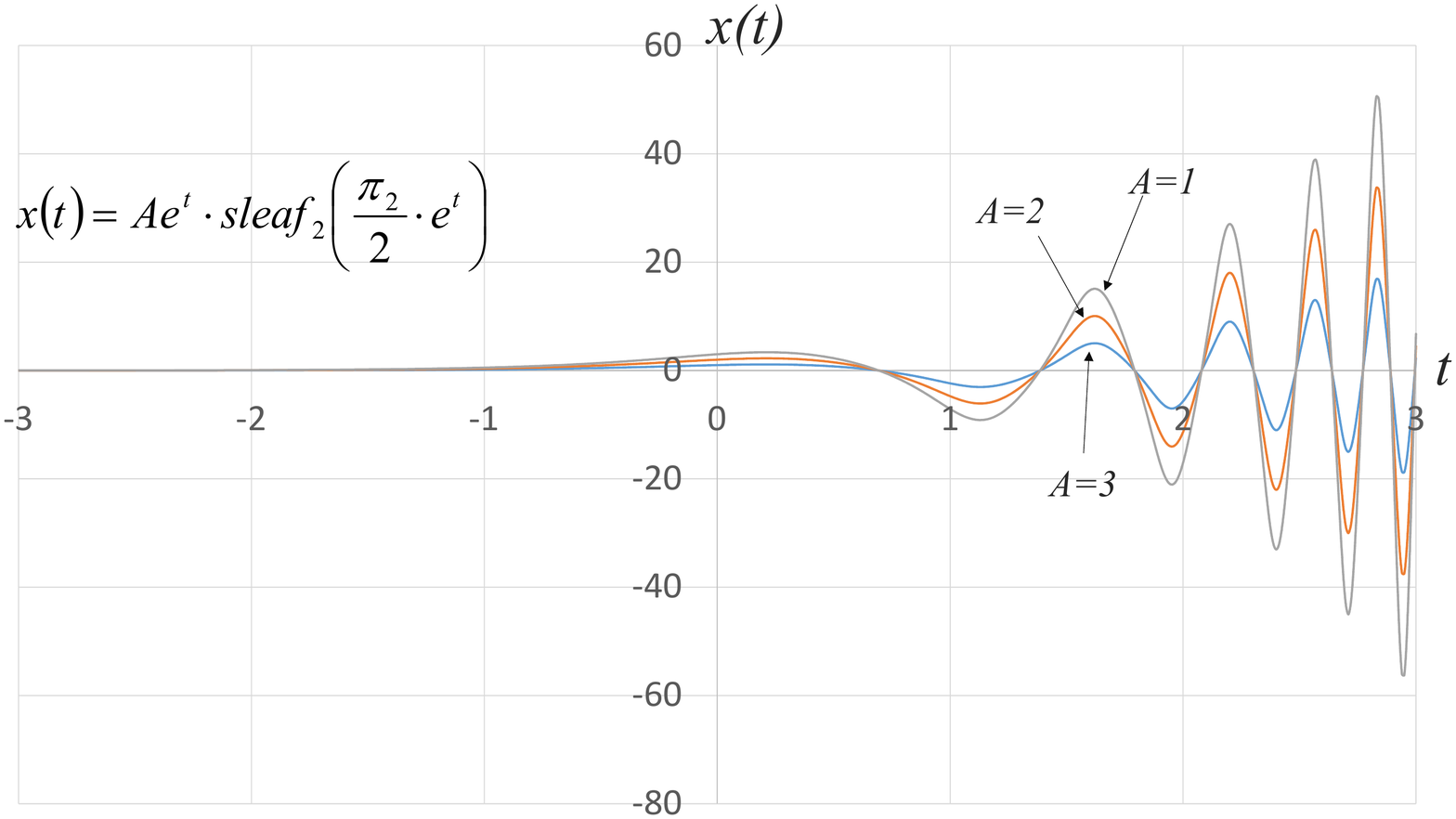}
\end{center}
\caption{ Waves obtained by the (X\hspace{-.1em}I\hspace{-.1em}I\hspace{-.1em}I) exact solution as a function of variations in the parameter $A$ ($A=1$, $2$, $3$) (Set $B=\pi_2/2$, $\omega=1$ and $\phi=0$)}

\label{fig:22}       
\end{figure*}
%
\begin{figure*}[tb]
\begin{center}
\includegraphics[width=0.75 \textwidth]{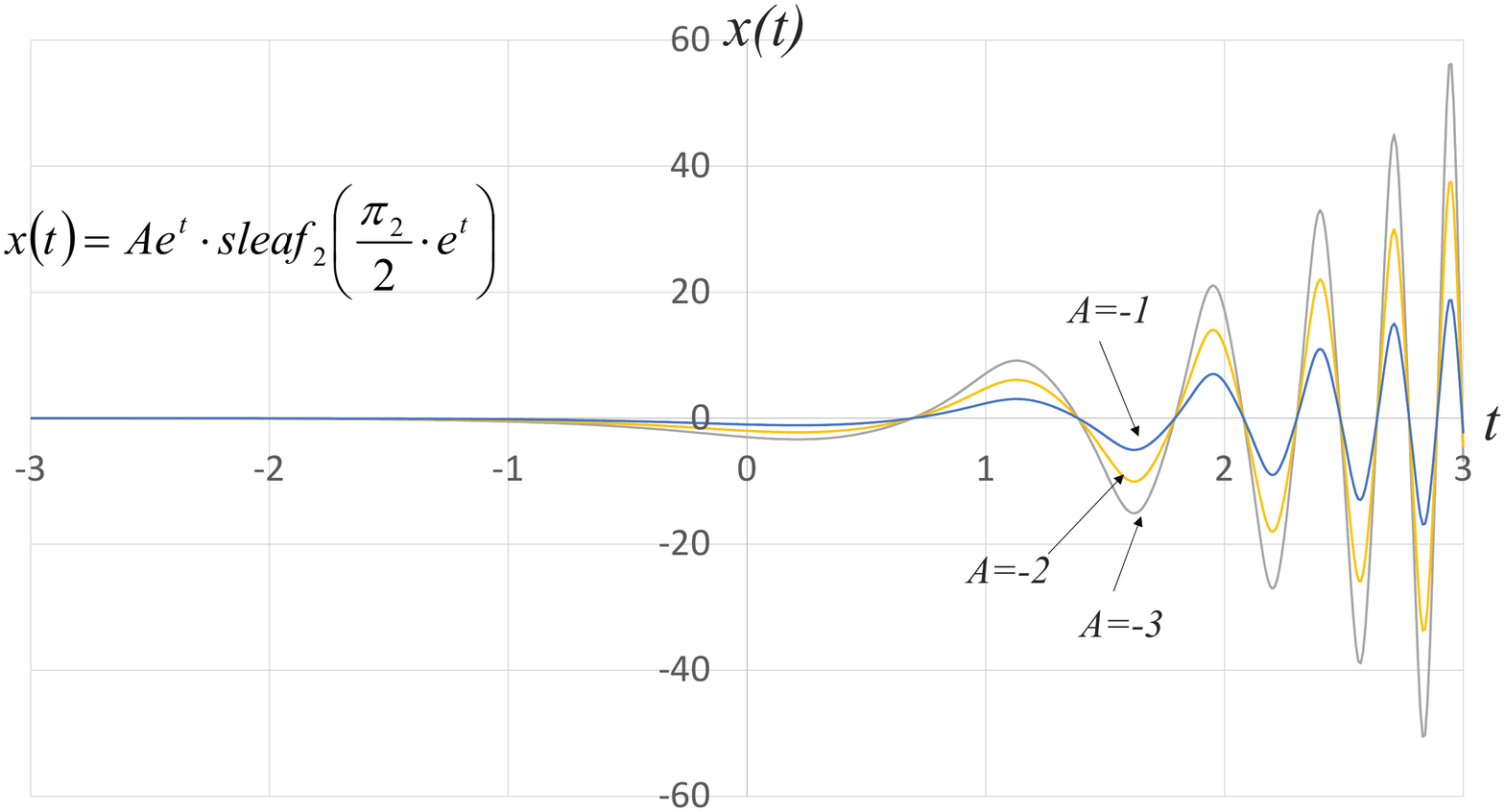}
\end{center}
\caption{ Waves obtained by the (X\hspace{-.1em}I\hspace{-.1em}I\hspace{-.1em}I) exact solution as a function of variations in the amplitude $A$ ($A=-1$, $-2$, $-3$) (Set $B=\pi_2/2$, $\omega=1$ and $\phi=0$) }
\label{fig:23}       
\end{figure*}
%
\begin{figure*}[tb]
\begin{center}
\includegraphics[width=0.75 \textwidth]{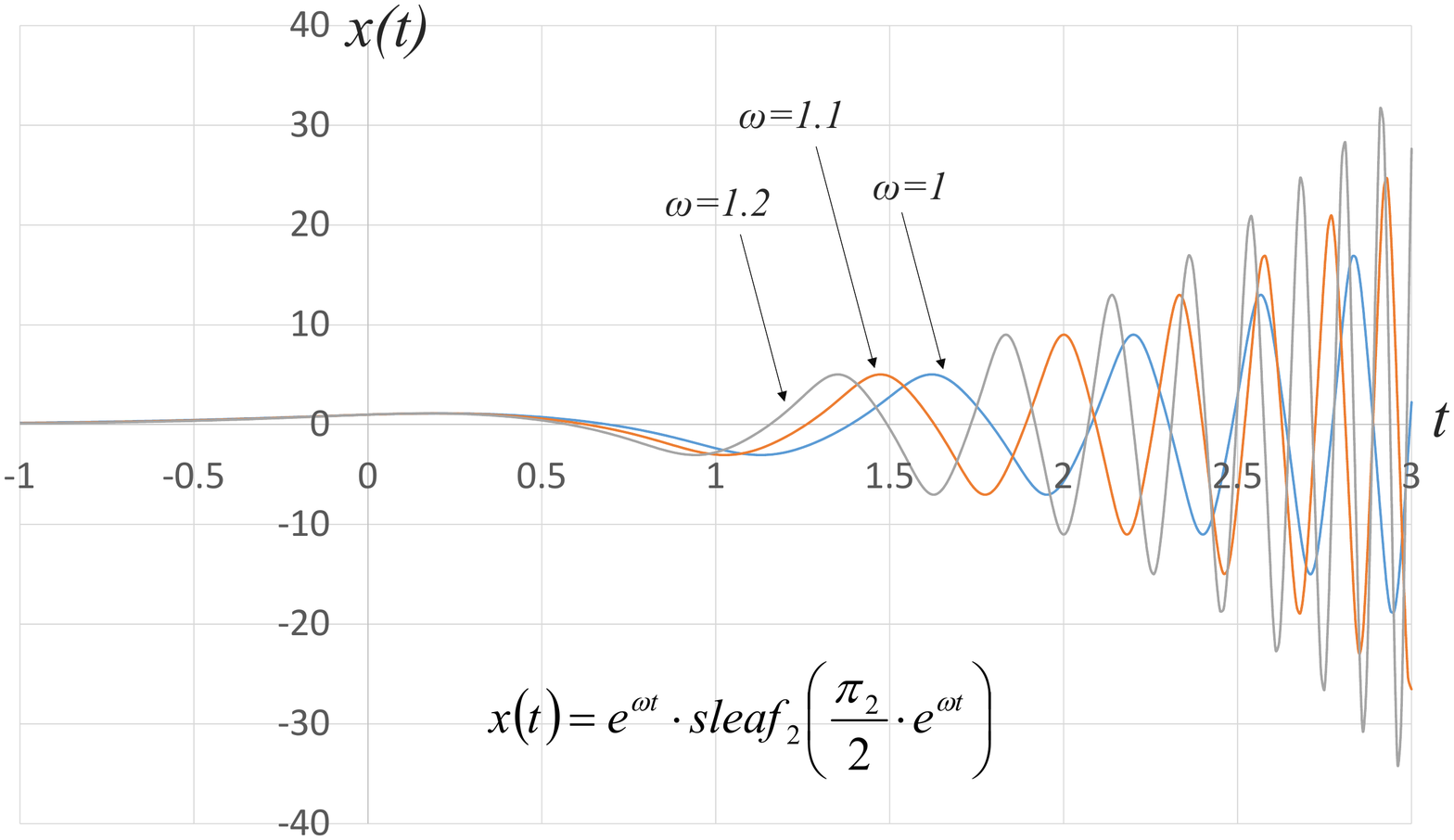}
\end{center}
\caption{ Waves obtained by the (X\hspace{-.1em}I\hspace{-.1em}I\hspace{-.1em}I) exact solution as a function of variations in the angular frequency $\omega$ ($\omega=1$, $1.1$, $1.2$) (Set $B=\pi_2/2$, $A=1$ and $\phi=0$)}

\label{fig:24}       
\end{figure*}
%
\begin{figure*}[tb]
\begin{center}
\includegraphics[width=0.75 \textwidth]{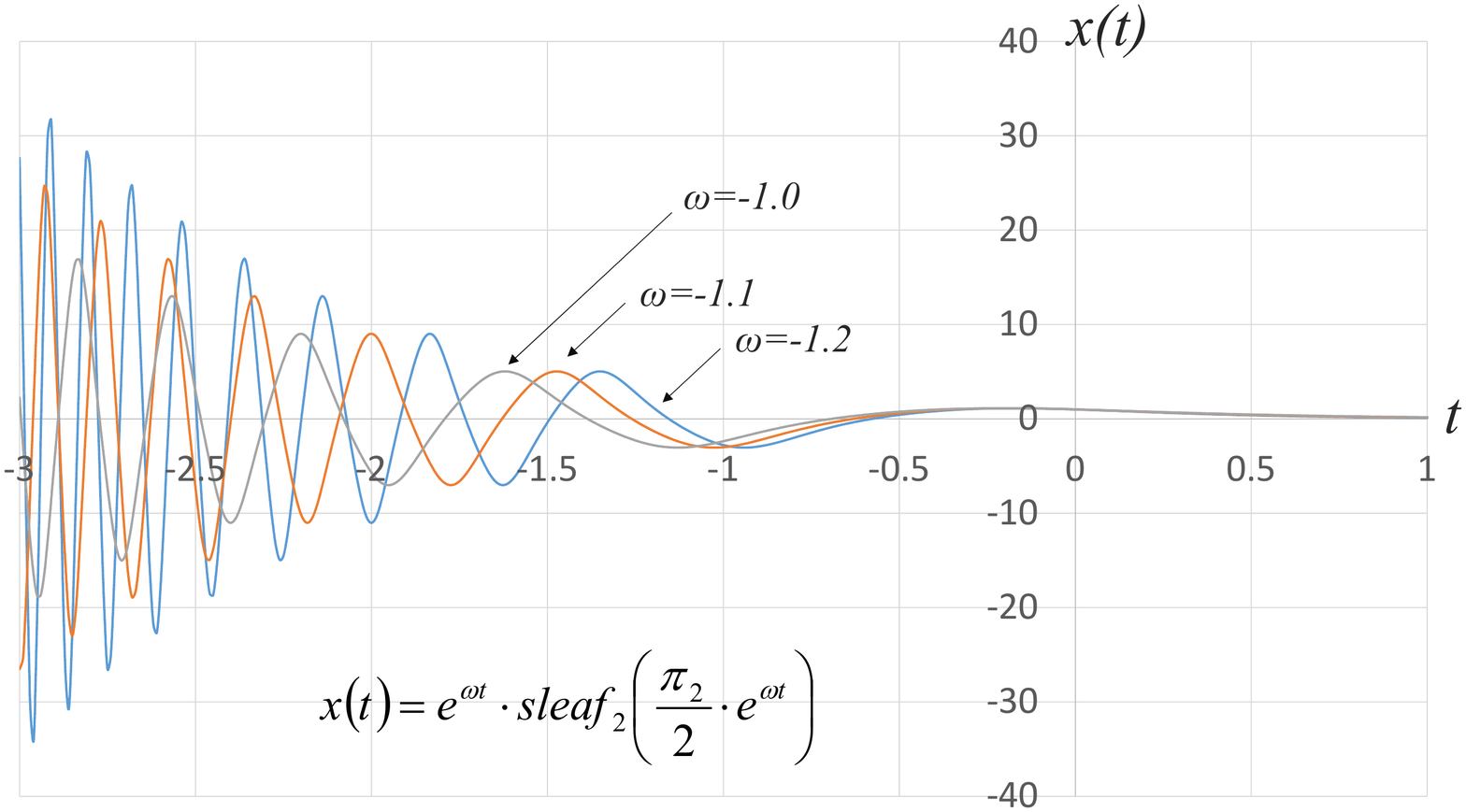}
\end{center}
\caption{ Waves obtained by the (X\hspace{-.1em}I\hspace{-.1em}I\hspace{-.1em}I) exact solution as a function of variations in the angular frequency $\omega$  ($\omega=-1$, $-1.1$, $-1.2$)  (Set $B=\pi_2/2$, $A=1$ and $\phi=0$) }
\label{fig:25}       
\end{figure*}

\begin{table}
\caption{ Numerical data for the (X\hspace{-.1em}I\hspace{-.1em}I\hspace{-.1em}I) solution under the conditions $A=1$, $B=1$, $\phi=-1$ and $\omega=1$}
\label{tab:11}
\centering
\begin{tabular}{cccc}
\hline\noalign{\smallskip}
$t$ & $x(t)$ & $x(t)^3$  & $\frac{\mathrm{d}^2x(t)}{\mathrm{d}t^2}$ \\
\noalign{\smallskip}\hline\noalign{\smallskip}
-1.0	&	-0.228880075	$\cdots$&	-0.011990132	$\cdots$&	0.169455704$\cdots$	\\
-0.8	&	-0.245174519	$\cdots$&	-0.014737574	$\cdots$&	0.362522502$\cdots$	\\
-0.6	&	-0.246594795	$\cdots$&	-0.014995181	$\cdots$&	0.668372752$\cdots$	\\
-0.4	&	-0.220730277	$\cdots$&	-0.010754389	$\cdots$&	1.140829977$\cdots$	\\
-0.2	&	-0.148394692	$\cdots$&	-0.003267797	$\cdots$&	1.868015936$\cdots$	\\
0.0	&	0	&	0	&	3	\\
0.2	&	0.270356983	$\cdots$&	0.019761175	$\cdots$&	4.700945926$\cdots$	\\
0.4	&	0.729443961	$\cdots$&	0.388128738	$\cdots$&	6.436195696$\cdots$	\\
0.6	&	1.432075704	$\cdots$&	2.936959316	$\cdots$&	3.391086916$\cdots$	\\
0.8	&	2.209335527	$\cdots$&	10.78412786	$\cdots$&	-16.83682468$\cdots$	\\
1.0	&	2.301435961	$\cdots$&	12.18980293	$\cdots$&	-37.53445744$\cdots$	\\
1.2	&	1.001645424	$\cdots$&	1.004944398	$\cdots$&	-33.94050058$\cdots$	\\
1.4	&	-1.750314905	$\cdots$&	-5.362268712	$\cdots$&	-39.49603198$\cdots$	\\
1.6	&	-4.951062236	$\cdots$&	-121.3654741	$\cdots$&	240.7147121	$\cdots$\\
1.8	&	-1.176291742	$\cdots$&	-1.627590492	$\cdots$&	111.7952215	$\cdots$\\
2.0	&	7.187993559	$\cdots$&	371.3838704	$\cdots$&	-682.6341556$\cdots$	\\
2.2	&	-1.433450358	$\cdots$&	-2.945425026	$\cdots$&	-239.8177858$\cdots$	\\
2.4	&	-5.102484331	$\cdots$&	-132.8449468	$\cdots$&	616.6525292	$\cdots$\\
2.6	&	8.555042892	$\cdots$&	626.1329714	$\cdots$&	-1741.232719$\cdots$	\\
2.8	&	-4.727862807	$\cdots$&	-105.6804362	$\cdots$&	1015.14273	$\cdots$\\
3.0	&	-14.27534157	$\cdots$&	-2909.10586	$\cdots$&	4759.409349$\cdots$	\\

\noalign{\smallskip}\hline
\end{tabular}
\end{table}

\subsection{Damping solution of the type (X\hspace{-.1em}I\hspace{-.1em}V)}
\label{4.7}
Let us consider Eq. (\ref{3.6.1}) that satisfies the Duffing equation. Physically, the exact solution of the type (X\hspace{-.1em}I\hspace{-.1em}V) corresponds to the nonlinear spring model of the damping system. 

\subsubsection{Period}
 When $A=1$, $\omega=1$, $\phi=0$ and $B=\pi_2/2$, the curves obtained by the (X\hspace{-.1em}I\hspace{-.1em}V) are shown in Fig.26. The horizontal axis represents time $t$ whereas the vertical axis represents displacement, $x(t)$. As shown in Fig. 26, the period of the wave varies with time.  In the Fig.26, the $T_1$, $T_2$, $\cdots$, $T_6$ represents the first period, the second period, $\cdots$, and the sixth period, respectively. First, let us consider the first period $T_1$.  Substituting -$\infty$ into the variable  $t$, the value of the type (X\hspace{-.1em}I\hspace{-.1em}V) can be obtained as follows:
\begin{equation}
x( -\infty)=\mathrm{e}^{- \infty} \cdot \mathrm{cleaf}_2(\frac{\pi_2}{2} \cdot e^{-\infty})=0 \cdot \mathrm{cleaf}_2(\frac{\pi_2}{2} \cdot 0)=0
\end{equation} 
On the other hand, substituting $\mathrm{ln}4$ into the variable $t$, the value of the type (X\hspace{-.1em}I\hspace{-.1em}V) can be obtained as follows:
\begin{equation}
x(\mathrm{ln}4)=e^{\mathrm{ln}4} \cdot \mathrm{cleaf}_2(\frac{\pi_2}{2} \cdot \mathrm{e}^{\mathrm{ln}4})=4 \cdot \mathrm{cleaf}_2(2 \pi_2)=4
\end{equation} 
The first period $T_1$ is as follows:
\begin{equation}
T_1=\mathrm{ln}4-\infty=-\infty
\end{equation} 
The domain of the first period $T_1$ is as follows: 
\begin{equation}
-\infty<t \leqq \mathrm{ln}4
\end{equation} 
Next, let us consider the second period $T_2$.  Substituting $\mathrm{ln}8$ into the variable $t$, the value of the type (X\hspace{-.1em}I\hspace{-.1em}V) can be obtained as follows:
\begin{equation}
x(\mathrm{ln}8)=\mathrm{e}^{\mathrm{ln}8} \cdot \mathrm{cleaf}_2(\frac{\pi_2}{2} \cdot \mathrm{e}^{\mathrm{ln}8})=8 \cdot \mathrm{cleaf}_2(4 \pi_2)=8
\end{equation} 
The second period $T_2$ is as follows:
\begin{equation}
T_2=\mathrm{ln}8-\mathrm{ln}4
\end{equation} 
The domain of the second period $T_2$ is as follows: 
\begin{equation}
\mathrm{ln}4 \leqq t \leqq \mathrm{ln}8
\end{equation} 
In this way, the $m^{th}$ period $T_m$ can be generalized as follows:
\begin{equation}
T_m=\mathrm{ln}(4m)-\mathrm{ln}(4m-4)=\mathrm{ln} \frac{m}{m-1} (m=1,2,3, \cdots)
\end{equation} 
The domain of the $m^{th}$ period $T_m$ is as follows: 
\begin{equation}
\mathrm{ln}(4m-4) \leqq t \leqq \mathrm{ln}(4m) (m=1,2,3, \cdots)
\end{equation} 

\subsubsection{Amplitude}
A time with respect to the convex upward of wave is obtained by using the gradient of type (XIV). Using the condition $\mathrm{d}x(t)/\mathrm{d}t=0$ in the Eq. (\ref{A14.1}), we can obtain the equation as follows:
\begin{equation}
B^2 \cdot \mathrm{e}^{2 \omega t} (\mathrm{cleaf}_2(B \cdot \mathrm{e}^{\omega t}))^4
+(\mathrm{cleaf}_2(B \cdot \mathrm{e}^{\omega t}))^2-B^2 \cdot \mathrm{e}^{2 \omega t}=0 \label{4.7.2.1}
\end{equation} 
We can not obtain the exact solution of the variable $t$ by the above equation. Using the numerical analysis, the variable $t$ that satisfy with the Eq. (\ref{4.7.2.1}) is obtained. The value $t$ and the amplitude $x(t)$ are summarized in the Table 12 and Table 13.

\begin{table}
\caption{ Exact amplitude $x(t)$ at convex upward for (X\hspace{-.1em}I\hspace{-.1em}V) }
\label{tab:12}       
\centering
\begin{tabular}{ccc}
\hline\noalign{\smallskip}
Number in Fig. 26 & Time $t$ & Exact amplitude $x(t)$ \\
\noalign{\smallskip}\hline\noalign{\smallskip}
(1U)	&	1.410965189 $\cdots$	&	4.064588954 $\cdots$	\\
(2U)	&2.085801435 $\cdots$	&	8.032995533 $\cdots$	\\
(3U)	&2.487749982 $\cdots$	&	12.02208815 $\cdots$	\\
(4U)	&2.774191427 $\cdots$	&	16.0165903 $\cdots$	\\
(5U)	&	2.996758996 $\cdots$	&	20.01328124 $\cdots$	\\
(6U)	&	3.178767206 $\cdots$	&	24.01107178 $\cdots$	\\
\noalign{\smallskip}\hline
\end{tabular}
\end{table}

\begin{table}
\caption{  Exact amplitude $x(t)$ at convex downward for the (X\hspace{-.1em}I\hspace{-.1em}V) }
\label{tab:13}       
\centering
\begin{tabular}{ccc}
\hline\noalign{\smallskip}
Number in Fig. 26 & Time $t$ & Exact amplitude $x(t)$ \\
\noalign{\smallskip}\hline\noalign{\smallskip}
(1D)	&	0.781810055 $\cdots$	&	-2.119912637 $\cdots$	\\
(2D)	&	1.802974309 $\cdots$	&	-6.043743891 $\cdots$	\\
(3D)	&	2.306670971 $\cdots$	&	-10.02646708 $\cdots$	\\
(4D)	&2.641148943 $\cdots$	&	-14.01894944 $\cdots$	\\
(5D)	&	2.891638805 $\cdots$	&	-18.01475276 $\cdots$	\\
(6D)	&	3.091891237 $\cdots$	&	-22.01207638 $\cdots$	\\
\noalign{\smallskip}\hline
\end{tabular}
\end{table}

The curve $x(t)=\pm \mathrm{e}^{t}$ is added in Fig. 26. As shown in Fig. 26, the curve $\mathrm{e}^{t}$ intersects the curve $\mathrm{cleaf}_2(\frac{\pi_2}{2} \mathrm{e}^{t})$ at  the condition $\mathrm{cleaf}_2(\frac{\pi_2}{2} \mathrm{e}^{t})= \pm 1$. In the case where the condition $\mathrm{cleaf}_2(\frac{\pi_2}{2} \mathrm{e}^{t})=1$, the time $t$ is obtained as follows:
\begin{equation}
\frac{\pi_2}{2} \mathrm{e}^{t}=2k \pi_2 (k=1,2,3 \cdots) \label{4.7.2.2}
\end{equation}
The following equation is obtained by the above equation.
\begin{equation}
t_k= \mathrm{ln}(4k) (k=1,2,3, \cdots) \label{4.7.2.3}
\end{equation}
Substituting $t_k$ into the Eq. (\ref{3.7.1}) for the solution of type (XIV), the following equation is obtained.
\begin{equation}
\begin{split}
x(t_k)&=\mathrm{exp}(t_k) \cdot \mathrm{cleaf}_2(\frac{\pi_2}{2} \cdot \mathrm{exp}(t_k))   \\
&= \mathrm{exp}(\mathrm{ln}(4k)) \cdot \mathrm{cleaf}_2(\frac{\pi_2}{2} \cdot \mathrm{exp}(\mathrm{ln}(4k))) \\
&=4k \cdot \mathrm{cleaf}_2(\frac{\pi_2}{2} \cdot 4k) \\
&=4k (k=1,2,3, \cdots) \label{4.7.2.4}
\end{split}
\end{equation}
In the case where the condition $\mathrm{cleaf}_2(\frac{\pi_2}{2} \mathrm{e}^{t})=-1$, the time $t$ is obtained as follows:
\begin{equation}
\frac{\pi_2}{2} \mathrm{e}^{t}=(2k-1) \pi_2 (k=1,2,3 \cdots) \label{4.7.2.5}
\end{equation}
The following equation is obtained by the above equation.
\begin{equation}
t_k= \mathrm{ln}(4k-2) (k=1,2,3, \cdots) \label{4.7.2.6}
\end{equation}
Substituting $t_k$ into the Eq. (\ref{3.7.1}), the following equation is obtained.
\begin{equation}
\begin{split}
x(t_k)&=\mathrm{exp}(t_k) \cdot \mathrm{cleaf}_2(\frac{\pi_2}{2} \cdot \mathrm{exp}( t_k))   \\
&= \mathrm{exp}(\mathrm{ln}(4k-2)) \cdot \mathrm{cleaf}_2(\frac{\pi_2}{2} \cdot \mathrm{exp}(\mathrm{ln}(4k-2))) \\
&=(4k-2) \mathrm{cleaf}_2(\frac{\pi_2}{2} \cdot (4k-2)) \\
&=-(4k-2) (k=1,2,3, \cdots) \label{4.7.2.7}
\end{split}
\end{equation}
\begin{table}
\caption{ Approximate amplitude $x(t)$ at convex upward for (XIV) }
\label{tab:14}       
\centering
\begin{tabular}{ccc}
\hline\noalign{\smallskip}
Number in Fig. 26 & Time $t$ by Eq. (\ref{4.7.2.3}) &
\begin{tabular}{c}
Approximate amplitude $x(t)$ \\
for (XIV) by Eq. (\ref{4.7.2.4})  
\end{tabular}\\ 
\noalign{\smallskip}\hline\noalign{\smallskip}
(1U)	&	1.386294361 $\cdots$ (=ln(4))	&	4.0 \\
(2U)	&	2.079441542 $\cdots$	(=ln(8)) &	8.0 	\\
(3U)	&2.48490665 $\cdots$	(=ln(12)) &	12.0 	\\
(4U)	&2.772588722 $\cdots$	(=ln(16)) &	16.0 	\\
(5U)	&	2.995732274 $\cdots$	(=ln(20)) &	20.0 	\\
(6U)	&3.17805383 $\cdots$	(=ln(26)) &	24.0 	\\
\noalign{\smallskip}\hline
\end{tabular}
\end{table}
\begin{table}
\caption{ Approximate amplitude $x(t)$ at convex downward for (XIV) }
\label{tab:15}       
\centering
\begin{tabular}{ccc}
\hline\noalign{\smallskip}
Number in Fig. 26 & Time $t$ by Eq. (\ref{4.7.2.6})&
\begin{tabular}{c}
Approximate amplitude $x(t)$ \\
for (XIV) by Eq. (\ref{4.7.2.7})  
\end{tabular}\\ 
\noalign{\smallskip}\hline\noalign{\smallskip}
(1D)	&0.693147181 $\cdots$ (=ln(2))	&	-2.0 \\
(2D)	&	1.791759469 $\cdots$	(=ln(6))&	-6.0 	\\
(3D)	&	2.302585093 $\cdots$	(=ln(10))&	-10.0 	\\
(4D)	&2.63905733 $\cdots$	(=ln(14))&	-14.0 	\\
(5D)	&	2.890371758 $\cdots$	(=ln(18))&	-18.0 	\\
(6D)	&	3.091042453 $\cdots$	(=ln(22))&	-22.0 	\\
\noalign{\smallskip}\hline
\end{tabular}
\end{table}
These data both the time $t$ (by the Eq. (\ref{4.7.2.3}) and the Eq. (\ref{4.7.2.6})  ) and the amplitude $x(t)$  (by the Eq. (\ref{4.7.2.4}) and the Eq. (\ref{4.7.2.7}))  are summarized in Table 14 and Table 15. Since the condition $\mathrm{cleaf}_2(\frac{\pi_2}{2} \mathrm{e}^{t})= \pm 1$ is not satisfied with the Eq.  (\ref{4.7.2.1}), the values $t$ obtained by the Eq. (\ref{4.7.2.3}) and the Eq. ( \ref{4.7.2.6}) are not the extremal value in the type (XIV). However, comparing the numerical data in the table 14(or the table 15) with the numerical data in the table 12 (or the table 13), the approximate amplitude (the table 14 and the table 15) are almost agreed with the exact amplitude by numerical analysis (the table 12 and the table 13). As an approximate value of the amplitude, amplitude can be easily obtained from the Eq. (\ref{4.7.2.4}) and the Eq. (\ref{4.7.2.7}). 
Fig. 27 and Fig. 28 show the variation of waves as the parameter $A$ vary. As the parameter $A$ increases, the height of the wave increases with keeping the wave period. Let us consider the parameter $\omega$ under the conditions $B = \pi_2/2$, $\phi = 0$ and $A=1$. The curves for the (X\hspace{-.1em}I\hspace{-.1em}V) solution are shown in the Fig. 29. In the negative domain of the variable $t$, the waves have no period. Conversely, in the positive domain of the variable $t$, the waves show periodicity. As the magnitude of $\omega$ increases, a period becomes shorter. As the magnitude of $\omega$ decreases, a period becomes longer. As shown in Fig 30, for $\omega<0$, the wave has no periodicity in the positive domain $t$. The exact solution for (X\hspace{-.1em}I\hspace{-.1em}V) decreases monotonically as time, $t$ increases. In the negative domain of the variable $t$, the wave shows periodicity.
\begin{table}
\caption{ Numerical data for the (X\hspace{-.1em}I\hspace{-.1em}I\hspace{-.1em}I) solution under the conditions $A=1$, $B=1$, $\phi=-1$ and $\omega=1$}
\label{tab:16}
\centering
\begin{tabular}{ccccc}
\hline\noalign{\smallskip}
$t$ & $x(t)$ & $x(t)^3$  & $\frac{\mathrm{d}^2x(t)}{\mathrm{d}t^2}$ \\
\noalign{\smallskip}\hline\noalign{\smallskip}
-1.0	&	-0.228880075	$\cdots$&	-0.011990132	$\cdots$&	0.169455704$\cdots$	\\
-0.8	&	-0.245174519	$\cdots$&	-0.014737574	$\cdots$&	0.362522502$\cdots$	\\
-0.6	&	-0.246594795	$\cdots$&	-0.014995181	$\cdots$&	0.668372752$\cdots$	\\
-0.4	&	-0.220730277	$\cdots$&	-0.010754389	$\cdots$&	1.140829977$\cdots$	\\
-0.2	&	-0.148394692	$\cdots$&	-0.003267797	$\cdots$&	1.868015936$\cdots$	\\
0.0	&	0	&	0	&	3	\\
0.2	&	0.270356983	$\cdots$&	0.019761175	$\cdots$&	4.700945926$\cdots$	\\
0.4	&	0.729443961	$\cdots$&	0.388128738	$\cdots$&	6.436195696$\cdots$	\\
0.6	&	1.432075704	$\cdots$&	2.936959316	$\cdots$&	3.391086916$\cdots$	\\
0.8	&	2.209335527	$\cdots$&	10.78412786	$\cdots$&	-16.83682468$\cdots$	\\
1.0	&	2.301435961	$\cdots$&	12.18980293	$\cdots$&	-37.53445744$\cdots$	\\
1.2	&	1.001645424	$\cdots$&	1.004944398	$\cdots$&	-33.94050058$\cdots$	\\
1.4	&	-1.750314905	$\cdots$&	-5.362268712	$\cdots$&	-39.49603198$\cdots$	\\
1.6	&	-4.951062236	$\cdots$&	-121.3654741	$\cdots$&	240.7147121	$\cdots$\\
1.8	&	-1.176291742	$\cdots$&	-1.627590492	$\cdots$&	111.7952215	$\cdots$\\
2.0	&	7.187993559	$\cdots$&	371.3838704	$\cdots$&	-682.6341556$\cdots$	\\
2.2	&	-1.433450358	$\cdots$&	-2.945425026	$\cdots$&	-239.8177858$\cdots$	\\
2.4	&	-5.102484331	$\cdots$&	-132.8449468	$\cdots$&	616.6525292	$\cdots$\\
2.6	&	8.555042892	$\cdots$&	626.1329714	$\cdots$&	-1741.232719$\cdots$	\\
2.8	&	-4.727862807	$\cdots$&	-105.6804362	$\cdots$&	1015.14273	$\cdots$\\
3.0	&	-14.27534157	$\cdots$&	-2909.10586	$\cdots$&	4759.409349$\cdots$	\\
\noalign{\smallskip}\hline
\end{tabular}
\end{table}
%
\begin{figure*}[tb]
\begin{center}
\includegraphics[width=0.75 \textwidth]{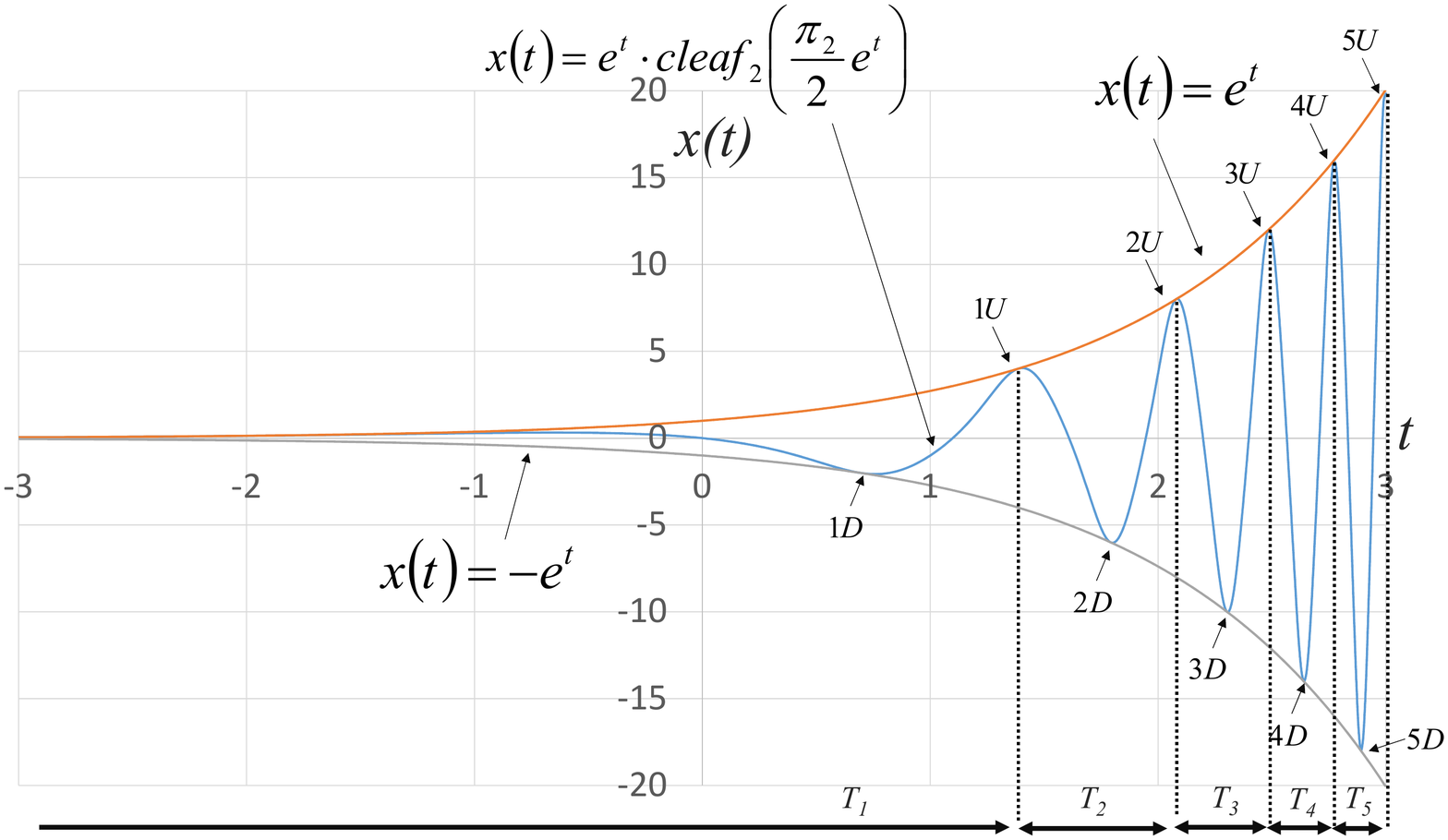}
\end{center}
\caption{ Waves obtained by the (X\hspace{-.1em}I\hspace{-.1em}V) exact solution and the exponential function $\mathrm{e}^t$ (Set $B=\pi_2/2$, $A=1$, $\omega=1$ and $\phi=0$) }
\label{fig:26}       
\end{figure*}
%
\begin{figure*}[tb]
\begin{center}
\includegraphics[width=0.75 \textwidth]{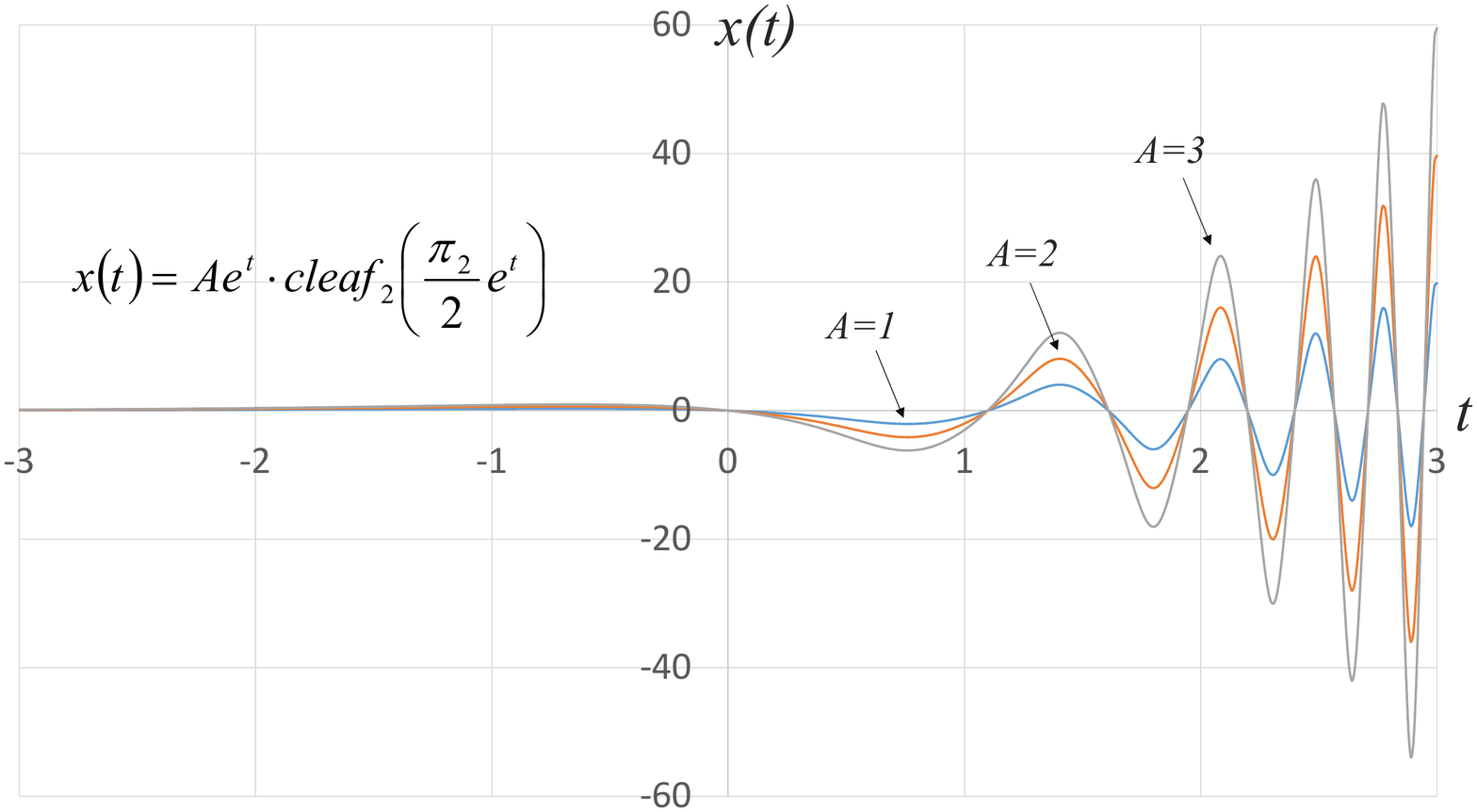}
\end{center}
\caption{ Waves obtained by the (X\hspace{-.1em}I\hspace{-.1em}V) exact solution as a function of variations in the parameter $A$  ($A=1$, $2$, $3$) (Set $B=\pi_2/2$, $\omega=1$ and $\phi=0$)}
\label{fig:27}       
\end{figure*}
%
\begin{figure*}[tb]
\begin{center}
\includegraphics[width=0.75 \textwidth]{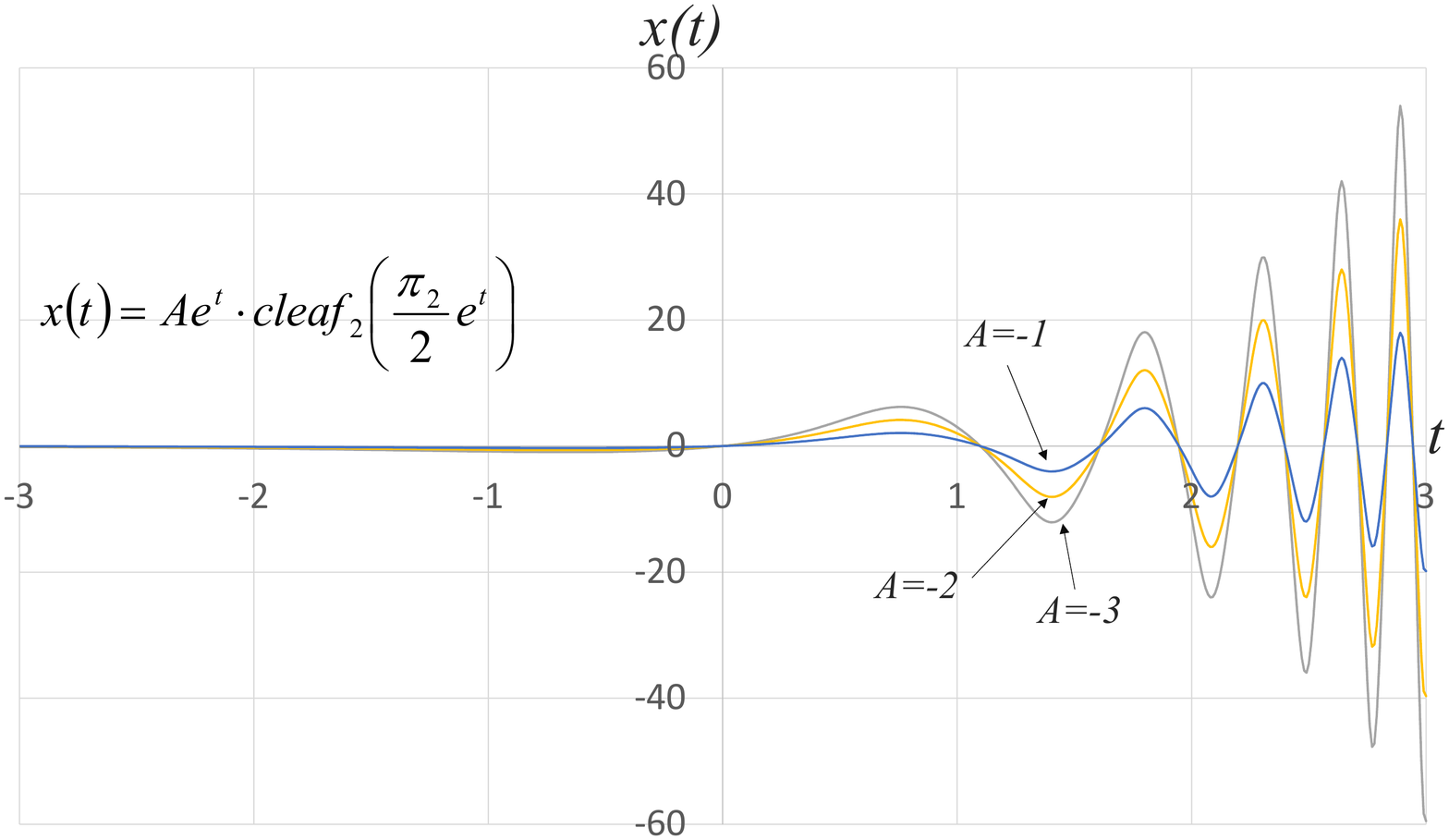}
\end{center}
\caption{ Waves obtained by the (X\hspace{-.1em}I\hspace{-.1em}V) exact solution as a function of variations in the amplitude  ($A=-1$, $-2$, $-3$) (Set $B=\pi_2/2$, $\omega=1$ and $\phi=0$) }
\label{fig:28}       
\end{figure*}
%
\begin{figure*}[tb]
\begin{center}
\includegraphics[width=0.75 \textwidth]{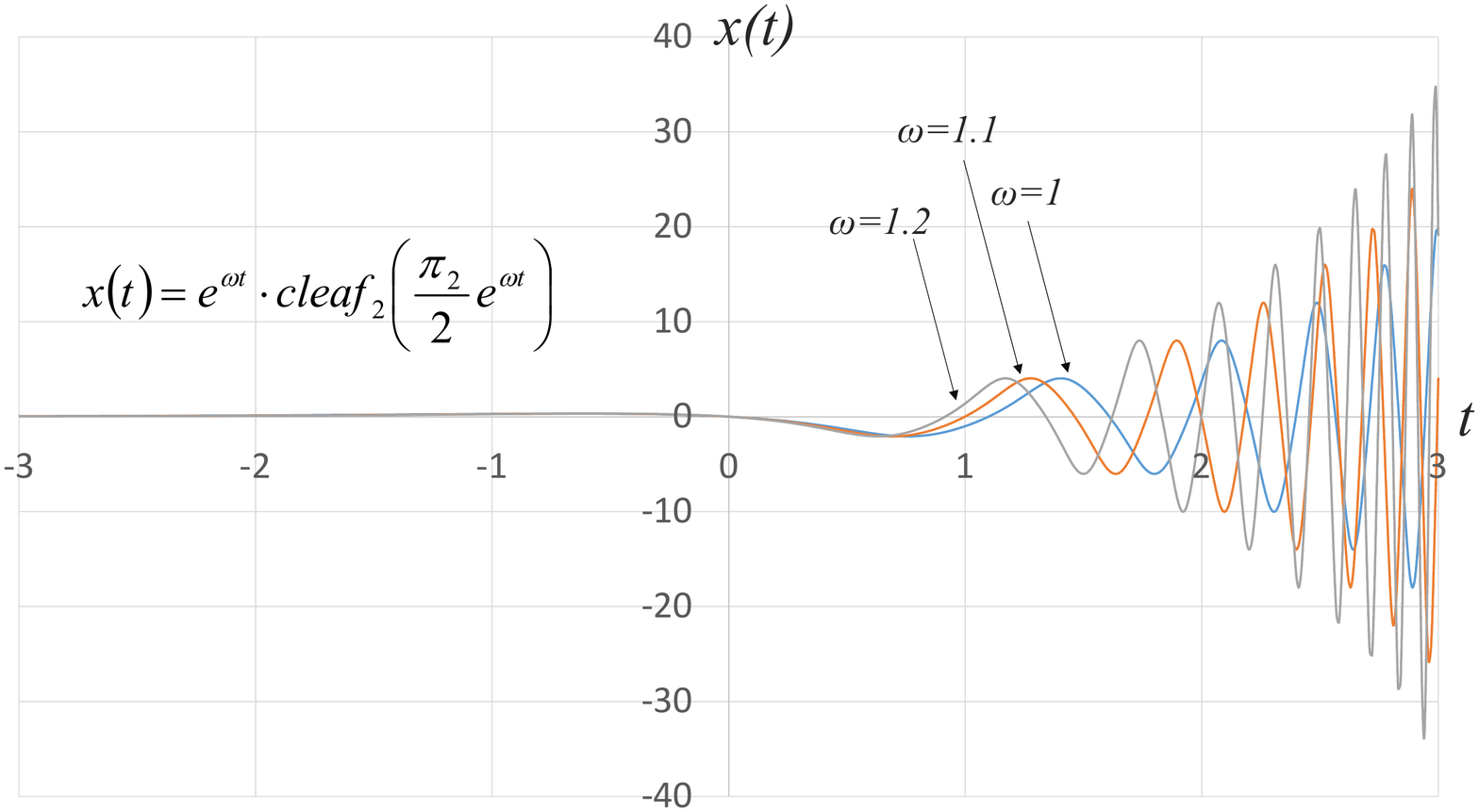}
\end{center}
\caption{ Waves obtained by the (X\hspace{-.1em}I\hspace{-.1em}V) exact solution as a function of variations in the angular frequency $\omega$ ($\omega=1.0$, $1.1$, $1.2$) (Set $B=\pi_2/2$, $A=1$ and $\phi=0$) }
\label{fig:29}       
\end{figure*}
%
\begin{figure*}[tb]
\begin{center}
\includegraphics[width=0.75 \textwidth]{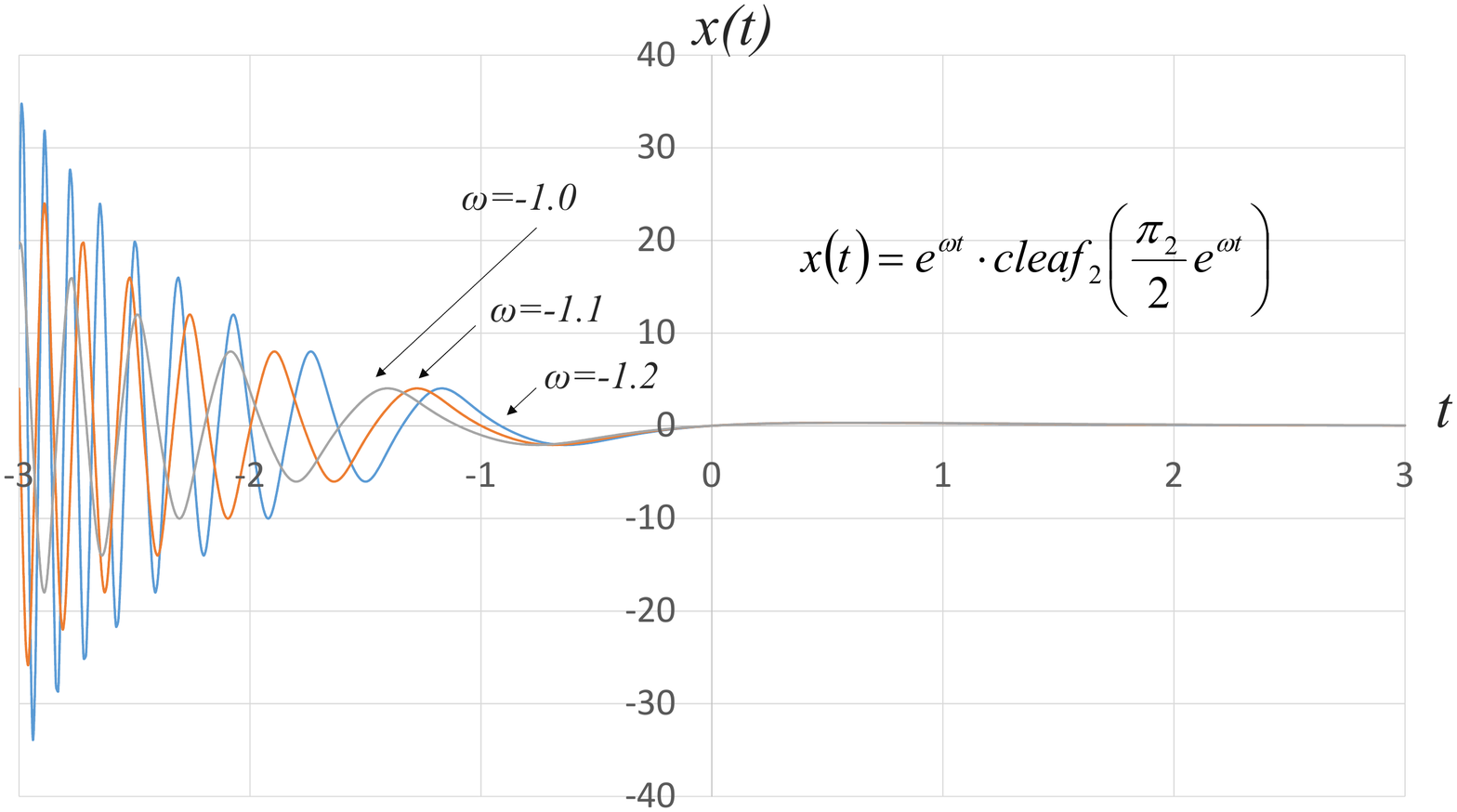}
\end{center}
\caption{ Waves obtained by the (X\hspace{-.1em}I\hspace{-.1em}V) exact solution as a function of variations in the angular frequency $\omega$   ($\omega=-1.0$, $-1.1$, $-1.2$) (Set $B=\pi_2/2$, $A=1$ and $\phi=0$) }
\label{fig:30}       
\end{figure*}
\begin{table}
\caption{ Numerical data for the (X\hspace{-.1em}I\hspace{-.1em}V) solution under the conditions $A=1$, $B=1$, $\phi=-1$ and $\omega=1$ }
\label{tab:17}  
\centering
\begin{tabular}{ccccc}
\hline\noalign{\smallskip}
$t$ & $x(t)$ & $x(t)^3$  & $\frac{\mathrm{d}^2x(t)}{\mathrm{d}t^2}$ \\
\noalign{\smallskip}\hline\noalign{\smallskip}
-1.0	&	0.244542736$\cdots$	&	0.014623937$\cdots$	&	0.579513405$\cdots$	\\
-0.8	&	0.330540411$\cdots$	&	0.036113842$\cdots$	&	0.76765176	$\cdots$\\
-0.6	&	0.447219762$\cdots$	&	0.089446419$\cdots$	&	0.943931473$\cdots$	\\
-0.4	&	0.601180255$\cdots$	&	0.217277185$\cdots$	&	0.967542439$\cdots$	\\
-0.2	&	0.792261897$\cdots$	&	0.497286086$\cdots$	&	0.503471721$\cdots$	\\
0.0	&	1	&	1	&	-1	\\
0.2	&	1.16295608	$\cdots$&	1.572859538	$\cdots$&	-3.87155113$\cdots$	\\
0.4	&	1.169058724$\cdots$	&	1.597750571	$\cdots$&	-7.295850123$\cdots$	\\
0.6	&	0.885786053$\cdots$	&	0.695002734	$\cdots$&	-10.18243951$\cdots$	\\
0.8	&	0.190255594$\cdots$	&	0.006886718	$\cdots$&	-14.68217896$\cdots$	\\
1	&	-1.103990273	$\cdots$&	-1.345537298	$\cdots$&	-20.27654537$\cdots$	\\
1.2	&	-3.030509208$\cdots$	&	-27.83215433$\cdots$	&	34.3447341	$\cdots$\\
1.4	&	-3.358521674$\cdots$	&	-37.88300891$\cdots$	&	108.3067894$\cdots$	\\
1.6	&	0.098793687	$\cdots$&	0.000964245	$\cdots$&	73.69438548$\cdots$	\\
1.8	&	5.825094124	$\cdots$&	197.6554719	$\cdots$&	-348.3448888$\cdots$	\\
2.0	&	1.227138245	$\cdots$&	1.84790855	$\cdots$&	-166.2009369	$\cdots$\\
2.2	&	-8.800138309$\cdots$	&	-681.5041324$\cdots$	&	1429.919587$\cdots$	\\
2.4	&	8.867240254	$\cdots$&	697.2129209	$\cdots$&	-1107.636432$\cdots$	\\
2.6	&	-8.774753605$\cdots$	&	-675.6235683$\cdots$	&	850.1465684	$\cdots$\\
2.8	&	15.13718215	$\cdots$&	3468.447383	$\cdots$&	-6490.862527$\cdots$	\\
3.0	&	-11.51700432$\cdots$	&	-1527.631446$\cdots$	&	4186.746094	$\cdots$\\
\noalign{\smallskip}\hline
\end{tabular}
\end{table}

\section{Conclusion}
\label{Conclusion}
In our previous paper, exact solutions that satisfy undamped and unforced Duffing equations were presented using the leaf function. The waveform obtained from the exact solution exhibits periodicity. The waveform has a different waveform from the sine and cosine waveforms.
In this paper, exact solutions for divergence and damped Duffing solutions were presented using leaf functions, hyperbolic leaf functions, and exponential functions.
When the solution diverges with time, the original function also needs to diverge with time.
However, there is no example in the literature showing the solution of the Duffing equation using the lemniscate function $\mathrm{slh}(t)$. The function $\mathrm{slh}(t)$ diverges as time, $t$ increases. To represent the exact solution for divergence in the Duffing equation, the integral function of the hyperbolic leaf function is implemented in the phase of the trigonometric function. The hyperbolic leaf function has a limit.
Depending on the limit of the hyperbolic leaf function, the exact solution for the divergent type also has a limit. One or two limits occur depending on the initial conditions.
In the exact solution for the damped Duffing equation, the exponential function was implemented in the phase of the leaf function to describe the exact solution. Over time, one period of the wave changes with time, depending on the condition of the parameter. Finally, the periodic state of the wave transitions to a state where there is no periodicity.

\appendix
\def\thesection{Appendix}
\section{ (V\hspace{-.1em}I\hspace{-.1em}I\hspace{-.1em}I)}

\label{Appendix8}
The exact solution of the type (V\hspace{-.1em}I\hspace{-.1em}I\hspace{-.1em}I) satisfies the cubic Duffing equation. The solution diverges with time, $t$.
The first derivative of the Eq. (\ref{3.1.1}) is obtained as follows:
\begin{equation}
\frac{\mathrm{d}x(t)}{\mathrm{d}t}=A \cdot \mathrm{sinh}(\mathrm{CLH}_2(\omega t + \phi)) \cdot \omega \cdot \mathrm{cleafh}_2(\omega t+ \phi) \label{A8.1}
\end{equation}
The second derivative of the above equation is obtained as follows.
\begin{equation}
\begin{split}
\frac{\mathrm{d}^2x(t)}{\mathrm{d}t^2}=A \cdot \mathrm{cosh}(\mathrm{CLH}_2(\omega t + \phi)) \cdot \omega^2 \cdot (\mathrm{cleafh}_2(\omega t+ \phi))^2 \\
+A \cdot \mathrm{sinh}(\mathrm{CLH}_2(\omega t + \phi)) \cdot \omega^2 \cdot \sqrt{(\mathrm{cleafh}_2(\omega t + \phi))^4-1 } \label{A8.2}
\end{split}
\end{equation}
The following equation is derived from Eq. (63) in Ref. \cite{Kaz_ch}.
\begin{equation}
(\mathrm{cleafh}_2(\omega t +\phi))^2=\mathrm{cosh}(2 \cdot \mathrm{CLH}_2(\omega t + \phi))  \label{A8.3}
\end{equation}
Using the above equation, the equation can be transformed as follows.
\begin{equation}
\begin{split}
\sqrt{(\mathrm{cleafh}_2(\omega t +\phi))^4-1} &
=\sqrt{(\mathrm{cosh}(2 \cdot \mathrm{CLH}_2(\omega t + \phi))^2-1} \\
&=\mathrm{sinh}(2 \cdot \mathrm{CLH}_2(\omega t + \phi)) \label{A8.4}
\end{split}
\end{equation}
Substituting the above two equations into equation (\ref{A8.2}) yields the following equation.
\begin{align}
\begin{split}
\frac{\mathrm{d}^2x(t)}{\mathrm{d}t^2}  &=  A \omega^2 \mathrm{cosh}(\mathrm{CLH}_2(\omega t + \phi)) \cdot \mathrm{cosh}(2 \cdot \mathrm{CLH}_2(\omega t + \phi)) \\
&+ A \omega^2 \mathrm{sinh}(\mathrm{CLH}_2(\omega t + \phi)) \cdot  \mathrm{sinh}(2 \cdot \mathrm{CLH}_2(\omega t + \phi)) \\
& = A\omega^2 \mathrm{cosh}(\mathrm{CLH}_2(\omega t + \phi)+2 \cdot \mathrm{CLH}_2(\omega t + \phi)) \\
& = A \omega^2 \mathrm{cosh}(3 \cdot \mathrm{CLH}_2(\omega t + \phi)) \\
& = -3 A \omega^2 cosh(\mathrm{CLH}_2(\omega t + \phi))
+4 (\frac{\omega}{A})^2 \Bigl\{ A \cdot 
\mathrm{cosh}(\mathrm{CLH}_2(\omega t + \phi)) \Bigl\}^3  \label{A8.5}
\end{split}
\end{align}
Substituting equation (\ref{3.1.1}) into the above equation leads to the following equation.
\begin{align}
\begin{split}
\frac{\mathrm{d}^2x(t)}{\mathrm{d}t^2}  = -3 \omega^2 x(t) + 4 (\frac{\omega}{A})^2 (x(t))^3 \label{A8.6}
\end{split}
\end{align}
The coefficients $\alpha$ and $\beta$ in Eq. (\ref{Duffing}) are as follows.
\begin{align}
\begin{split}
\alpha=3 \omega^2 \label{A8.7}
\end{split}
\end{align}
\begin{align}
\begin{split}
\beta=-4 (\frac{\omega}{A})^2  \label{A8.8}
\end{split}
\end{align}

\appendix
\def\thesection{Appendix}
\section{ (I\hspace{-.1em}X)}
\label{Appendix9}
The exact solution of the type (I\hspace{-.1em}X) satisfies the cubic Duffing equation. The solution diverges with time, $t$.
The first derivative of the above equation is obtained as follows.
\begin{equation}
\frac{\mathrm{d}x(t)}{\mathrm{d}t}=A \cdot \mathrm{cosh}(\mathrm{CLH}_2(\omega t + \phi)) \cdot \omega \cdot \mathrm{cleafh}_2(\omega t+ \phi) \label{A9.1}
\end{equation}
The second derivative of the above equation is obtained as follows.
\begin{equation}
\begin{split}
\frac{\mathrm{d}^2x(t)}{\mathrm{d}t^2}=A \cdot \mathrm{sinh}(\mathrm{CLH}_2(\omega t + \phi)) \cdot \omega^2 \cdot (\mathrm{cleafh}_2(\omega t+ \phi))^2 \\
+A \cdot \mathrm{cosh}(\mathrm{CLH}_2(\omega t + \phi)) \cdot \omega^2 \cdot \sqrt{(\mathrm{cleafh}_2(\omega t + \phi))^4-1 } \label{A9.2}
\end{split}
\end{equation}
Substituting (\ref{A8.3}) and (\ref{A8.4}) into the above equation gives the following equation.
\begin{align}
\begin{split}
\frac{\mathrm{d}^2x(t)}{\mathrm{d}t^2}  &=  A\omega^2 \mathrm{sinh}(\mathrm{CLH}_2(\omega t + \phi)) \cdot \mathrm{cosh}(2 \cdot \mathrm{CLH}_2(\omega t + \phi)) \\
&+ A \omega^2 \mathrm{cosh}(\mathrm{CLH}_2(\omega t + \phi)) \cdot \mathrm{sinh}(2 \cdot \mathrm{CLH}_2(\omega t + \phi)) \\
& = A\omega^2 \mathrm{sinh}(\mathrm{CLH}_2(\omega t + \phi)+2 \cdot \mathrm{CLH}_2(\omega t + \phi)) \\
& = A \omega^2 \mathrm{sinh}(3 \cdot \mathrm{CLH}_2(\omega t + \phi)) \\
& = 3 A \omega^2 \mathrm{sinh}(\mathrm{CLH}_2(\omega t + \phi))
+4 (\frac{\omega}{A})^2 \Bigl\{A \cdot \mathrm{sinh}(\mathrm{CLH}_2(\omega t + \phi)) \Bigl\}^3   \label{A9.3}
\end{split}
\end{align}
Substituting (\ref{3.2.1}) into the above equation leads to the following equation.
\begin{align}
\begin{split}
\frac{\mathrm{d}^2x(t)}{\mathrm{d}t^2}  = 3 \omega^2 x(t) + 4 (\frac{\omega}{A})^2 (x(t))^3 \label{A9.4}
\end{split}
\end{align}
The coefficients $\alpha$ and $\beta$ for the Duffing equation(Eq.(\ref{Duffing})) are as follows.
\begin{align}
\begin{split}
\alpha=-3 \omega^2 \label{A9.5}
\end{split}
\end{align}
\begin{align}
\begin{split}
\beta=-4 (\frac{\omega}{A})^2 \label{A9.6}
\end{split}
\end{align}

\appendix
\def\thesection{Appendix}
\section{ (X)}
\label{Appendix10}
The exact solution of the type (X) satisfies the cubic Duffing equation. The solution diverges with time, $t$.
The first derivative of the above equation is obtained as follows.
\begin{align}
\begin{split}
\frac{\mathrm{d}x(t)}{\mathrm{d}t} &=-A \cdot \mathrm{sin}(\mathrm{SL}_2(\omega t+\phi)) \cdot \omega \cdot \mathrm{sleaf}_2(\omega t+ \phi) \\
&-A \cdot \mathrm{cos}(\mathrm{SL}_2(\omega t+\phi)) \cdot \omega \cdot \mathrm{sleaf}_2(\omega t+ \phi) \\
&+\sqrt{2} A \cdot \mathrm{sinh}(\mathrm{CLH}_2(\omega t + \phi)) \cdot \omega \cdot \mathrm{cleafh}_2(\omega t+ \phi) \\
\label{A10.1}
\end{split}
\end{align}
The second derivative of the above equation is obtained as follows.
\begin{align}
\begin{split}
\frac{\mathrm{d}^2x(t)}{\mathrm{d}t^2} &=-A \cdot \mathrm{cos}(\mathrm{SL}_2(\omega t+\phi)) \cdot \omega^2 \cdot (\mathrm{sleaf}_2(\omega t+ \phi))^2 \\
&- A \cdot \mathrm{sin}(\mathrm{SL}_2(\omega t+\phi)) \cdot \omega^2 \cdot \sqrt{1-(\mathrm{sleaf}_2(\omega t + \phi))^4} \\
&+A \cdot \mathrm{sin}(\mathrm{SL}_2(\omega t+\phi)) \cdot \omega^2 \cdot (\mathrm{sleaf}_2(\omega t+ \phi))^2 \\
&- A \cdot \mathrm{cos}(\mathrm{SL}_2(\omega t+\phi)) \cdot \omega^2 \cdot \sqrt{1-(\mathrm{sleaf}_2(\omega t + \phi))^4} \\
&+\sqrt{2} A \cdot  \mathrm{cosh}(\mathrm{CLH}_2(\omega t + \phi)) \cdot \omega^2 \cdot (\mathrm{cleafh}_2(\omega t+ \phi))^2 \\
&+\sqrt{2} A \cdot \mathrm{sinh}(\mathrm{CLH}_2(\omega t + \phi)) \cdot \omega^2 \cdot \sqrt{(\mathrm{cleafh}_2(\omega t + \phi))^4-1} \label{A10.2}
\end{split}
\end{align}
Using Eq. (III.4) in Ref. \cite{Kaz_duf}, the above equation can be transformed using the following equation.
\begin{align}
\begin{split}
\frac{\mathrm{d}^2x(t)}{\mathrm{d}t^2} &=-A \omega^2 \mathrm{cos}(\mathrm{SL}_2(\omega t+\phi)) \cdot \mathrm{sin}(2 \mathrm{SL}_2(\omega t+\phi)) \\
&- A \omega^2  \mathrm{sin}(\mathrm{SL}_2(\omega t+\phi)) \cdot \mathrm{cos}(2 \mathrm{SL}_2(\omega t+\phi)) \\
&+A \omega^2 \mathrm{sin}(\mathrm{SL}_2(\omega t+\phi)) \cdot  \mathrm{sin}(2\mathrm{SL}_2(\omega t+\phi)) \\
&- A \omega^2 \mathrm{cos}(\mathrm{SL}_2(\omega t+\phi)) \cdot \mathrm{cos}(2\mathrm{SL}_2(\omega t+\phi)) \\
&+\sqrt{2} \omega^2 A \mathrm{cosh}(\mathrm{CLH}_2(\omega t + \phi)) \cdot \mathrm{cosh}(2 \mathrm{CLH}_2(\omega t + \phi)) \\
&+\sqrt{2} A \omega^2 \mathrm{sinh}(\mathrm{CLH}_2(\omega t + \phi)) \cdot \mathrm{sinh}(2\mathrm{CLH}_2(\omega t + \phi)) \\
\label{A10.3}
\end{split}
\end{align}
Using the addition theorem, we can summarize the equation as follows.
\begin{align}
\begin{split}
\frac{\mathrm{d}^2x(t)}{\mathrm{d}t^2} &=- A \omega^2  \mathrm{sin}(2 \cdot \mathrm{SL}_2(\omega t+\phi) + \mathrm{SL}_2(\omega t+\phi)) \\
&-A \omega^2 \mathrm{cos}(2 \cdot \mathrm{SL}_2(\omega t+\phi) + \mathrm{SL}_2(\omega t+\phi) ) \\
&+\sqrt{2} \omega^2 A \mathrm{cosh}(\mathrm{CLH}_2(\omega t + \phi) + 2 \cdot \mathrm{CLH}_2(\omega t + \phi)) \\
&=-A \omega^2 \mathrm{cos}(3 \cdot \mathrm{SL}_2(\omega t+\phi))
-A \omega^2 \mathrm{sin}(3 \cdot \mathrm{SL}_2(\omega t+\phi)) \\
&+\sqrt{2} A \omega^2 \mathrm{cosh}(3 \cdot \mathrm{CLH}_2(\omega t + \phi)) \label{A10.4}
\end{split}
\end{align}
Using the triple-angle formula, the above equation can be transformed as follows.
\begin{align}
\begin{split}
\frac{\mathrm{d}^2x(t)}{\mathrm{d}t^2} &=
-4(\frac{\omega}{A})^2 \Bigl\{ A \cdot \mathrm{cos}( \mathrm{SL}_2(\omega t+\phi) \Bigr\}^3  + 3A \omega^2 \mathrm{cos}( \mathrm{SL}_2(\omega t+\phi)) \\
&+4(\frac{\omega}{A})^2 \Bigl\{ A \cdot \mathrm{sin}( \mathrm{SL}_2(\omega t+\phi) \Bigr\}^3 - 3A \omega^2 \mathrm{sin}( \mathrm{SL}_2(\omega t+\phi)) \\
&+4\sqrt{2}(\frac{\omega}{A})^2 \Bigl\{ A \cdot \mathrm{cosh}( \mathrm{CLH}_2(\omega t + \phi) \Bigr\}^3  - 3\sqrt{2} \omega^2 \cdot A \cdot \mathrm{cosh}( \mathrm{CLH}_2(\omega t + \phi))  \label{A10.5}
\end{split}
\end{align}
Consider the following equation (See Eq. (F15) in Ref. \cite{Kaz_ch}).
\begin{align}
\begin{split}
-(\mathrm{sleaf}_2(\omega t+\phi))^2+(\mathrm{cleafh}_2(\omega t+\phi))^2-(\mathrm{sleaf}_2(\omega t+\phi))^2(\mathrm{cleafh}_2(\omega t+\phi))^2=1
\label{A10.6}
\end{split}
\end{align}
Using the Eq. (III.4) in Ref. \cite{Kaz_duf} and the Eq. (H5) in Ref. \cite{Kaz_ch}, the following equation is obtained.
\begin{align}
\begin{split}
(1-\mathrm{sin}(2 \cdot \mathrm{SL}_2(\omega t+\phi))(1+\mathrm{cosh}(2 \cdot \mathrm{CLH}_2(\omega t+\phi))=2 \label{A10.6.a1}
\end{split}
\end{align}
The following equation is obtained from the above equation.
\begin{align}
\begin{split}
&\Bigl\{ \mathrm{sin}(\mathrm{SL}_2(\omega t+\phi))^2+\mathrm{cos}(\mathrm{SL}_2(\omega t+\phi))^2
-2\mathrm{sin}(\mathrm{SL}_2(\omega t+\phi) \cdot \mathrm{cos}(\mathrm{SL}_2(\omega t+\phi))\Bigr\} \\ 
& \cdot \Bigl\{ 1+2(\mathrm{cosh}(\mathrm{CLH}_2(\omega t+\phi))^2-1 \Bigr\} =2 \label{A10.6.a2}
\end{split}
\end{align}
The following equation is obtained from the above equation.
\begin{align}
\begin{split}
\mathrm{cosh}( \mathrm{CLH}_2(\omega t + \phi))
 \left\{ \mathrm{cos}( \mathrm{SL}_2(\omega t+\phi))-\mathrm{sin}( \mathrm{SL}_2(\omega t+\phi))
\right\}=1
\label{A10.7}
\end{split}
\end{align}
The equation is expanded as follows.
\begin{align}
\begin{split}
& \left\{-A \cdot \mathrm{sin}( \mathrm{SL}_2(\omega t+\phi)+A \cdot \mathrm{cos}( \mathrm{SL}_2(\omega t+\phi))\right\}^3 \\
&= -A^3 \Bigl\{ \mathrm{sin}( \mathrm{SL}_2(\omega t+\phi)) \Bigr\}^3 
+3A^3 \Bigl\{ \mathrm{sin}( \mathrm{SL}_2(\omega t+\phi)) \Bigr\}^2 \cdot \mathrm{cos}( \mathrm{SL}_2(\omega t+\phi)) \\
&-3A^3 \mathrm{sin}( \mathrm{SL}_2(\omega t+\phi)) 
\cdot \Bigl\{ \mathrm{cos}( \mathrm{SL}_2(\omega t+\phi))\Bigr\}^2
+A^3 \Bigl\{\mathrm{cos}( \mathrm{SL}_2(\omega t+\phi))\Bigr\}^3 \\
&= -A^3 \Bigl\{\mathrm{sin}( \mathrm{SL}_2(\omega t+\phi))\Bigr\}^3 
+A^3 \Bigl\{\mathrm{cos}( \mathrm{SL}_2(\omega t+\phi))\Bigr\}^3 \\
&+3A^3 \mathrm{cos}( \mathrm{SL}_2(\omega t+\phi)) 
-3A^3 \Bigl\{\mathrm{cos}( \mathrm{SL}_2(\omega t+\phi))\Bigr\}^3 \\
&-3A^3 \mathrm{sin}( \mathrm{SL}_2(\omega t+\phi)) 
+3A^3 \Bigl\{\mathrm{sin}( \mathrm{SL}_2(\omega t+\phi))\Bigr\}^3
\label{A10.8}
\end{split}
\end{align}
The variable $x(t)^3$ is expanded as follows.
\begin{align}
\begin{split}
&(x(t))^3=A^3\Bigl\{\mathrm{cos}( \mathrm{SL}_2(\omega t+\phi))
-\mathrm{sin}( \mathrm{SL}_2(\omega t+\phi))\Bigr\}^3 
+2\sqrt{2}A^3 \Bigl\{ \mathrm{cosh}( \mathrm{CLH}_2(\omega t + \phi))\Bigr\}^3 \\
&+3\sqrt{2} A^3 \mathrm{cosh}( \mathrm{CLH}_2(\omega t + \phi))
\Bigl\{\mathrm{cos}( \mathrm{SL}_2(\omega t+\phi))
-\mathrm{sin}( \mathrm{SL}_2(\omega t+\phi))\Bigr\}^2 \\
&+6A^3 \Bigl\{ \mathrm{cosh}( \mathrm{CLH}_2(\omega t + \phi))\Bigr\}^2
\Bigl\{\mathrm{cos}( \mathrm{SL}_2(\omega t+\phi))
-\mathrm{sin}( \mathrm{SL}_2(\omega t+\phi))
\Bigr\} \label{A10.9}
\end{split}
\end{align}
The above equation is applied to Eq. (\ref{A10.7}). The following equation is obtained from the above equation.
\begin{align}
\begin{split}
&(x(t))^3=-2A^3(\mathrm{cos}( \mathrm{SL}_2(\omega t+\phi))^3
+2 A^3 \mathrm{sin}( \mathrm{SL}_2(\omega t+\phi)))^3 \\
&+2\sqrt{2}A^3 (\mathrm{cosh}( \mathrm{CLH}_2(\omega t + \phi)))^3 
-3(1+\sqrt{2})A^3\mathrm{sin}(\mathrm{SL}_2(\omega t+\phi))) \\
&+3(1+\sqrt{2})A^3\mathrm{cos}(\mathrm{SL}_2(\omega t+\phi)))
+6A^3\mathrm{cosh}(\mathrm{CLH}_2(\omega t + \phi))) \label{A10.10}
\end{split}
\end{align}
The ordinary differential equation (\ref{3.3.2}) are obtained.

\appendix
\def\thesection{Appendix}
\section{ (X\hspace{-.1em}I)}
\label{Appendix11}
The exact solution of the type (X\hspace{-.1em}I) satisfies the cubic Duffing equation. The solution diverges with time, $t$.
The first derivative of the equation is obtained as follows.
\begin{align}
\begin{split}
\frac{\mathrm{d}x(t)}{\mathrm{d}t}
&=A \omega \mathrm{e}^{\omega t} \mathrm{sleafh}_2(B e^{\omega t} +\phi) 
+A B \omega \mathrm{e}^{2 \omega t} 
\cdot \sqrt{1+(\mathrm{sleafh}_2(Be^{\omega t}+\phi))^4}
\label{A11.1}
\end{split}
\end{align}
The second derivative of the above equation is obtained as follows.
\begin{align}
\begin{split}
\frac{\mathrm{d}^2x(t)}{\mathrm{d}t^2}
&=A \omega^2 \mathrm{e}^{\omega t} \mathrm{sleafh}_2(B \cdot \mathrm{e}^{\omega t} +\phi) 
+A \omega \mathrm{e}^{\omega t} \sqrt{1+(\mathrm{sleafh}_2(B\mathrm{e}^{\omega t}+\phi))^4}
(B \omega \mathrm{e}^{\omega t} ) \\
&+ 2 A B \omega^2 \mathrm{e}^{2 \omega t} \sqrt{1+(\mathrm{sleafh}_2(B\mathrm{e}^{\omega t}+\phi))^4} 
+2 A B \omega \mathrm{e}^{2 \omega t} (\mathrm{sleafh}_2(B\mathrm{e}^{\omega t}+\phi))^3 \cdot (B\omega \cdot \mathrm{e}^{\omega t}) \\
&=A \omega^2 \mathrm{e}^{\omega t} \mathrm{sleafh}_2(B \cdot \mathrm{e}^{\omega t} +\phi)
+2A B^2 \omega^2 \mathrm{e}^{3 \omega t} (\mathrm{sleafh}_2(B \mathrm{e}^{\omega t} +\phi))^3 \\
&+3AB \omega^2 \mathrm{e}^{2 \omega t} \sqrt{1+(\mathrm{sleafh}_2(B\mathrm{e}^{\omega t}+\phi))^4} 
\label{A11.2}
\end{split}
\end{align}
Substituting these equations into the Duffing equation (\ref{Duffing}) yields the following equation.
\begin{align}
\begin{split}
&A \omega^2 \mathrm{e}^{\omega t} \mathrm{sleafh}_2(B\mathrm{e}^{\omega t}+\phi)
+2AB^2 \omega^2 \mathrm{e}^{3\omega t} (\mathrm{sleafh}_2(B\mathrm{e}^{\omega t}+\phi))^3 \\
& +3AB \omega^2 \mathrm{e}^{2\omega t}\sqrt{1+(\mathrm{sleafh}_2(B \mathrm{e}^{\omega t}+\phi))^4}  \\
&+\delta(
A\omega \mathrm{e}^{\omega t} \mathrm{sleafh}_2(B \mathrm{e}^{\omega t} +\phi) 
+AB\omega \mathrm{e}^{2\omega t} \sqrt{1+(\mathrm{sleafh}_2(B\mathrm{e}^{\omega t}+\phi))^4}) \\
&+\alpha A \mathrm{e}^{\omega t} \mathrm{sleafh}_2(B \mathrm{e}^{\omega t} +\phi)
+\beta( A \mathrm{e}^{\omega t} \mathrm{sleafh}_2(B \mathrm{e}^{\omega t} +\phi) )^3=0
\label{A11.3}
\end{split}
\end{align}
The above equation is summarized as follows.
\begin{align}
\begin{split}
&A (\omega^2+\alpha+\delta \omega) \mathrm{e}^{\omega t} \mathrm{sleafh}_2(B\mathrm{e}^{\omega t}+\phi) + A(2B^2 \omega^2+\beta A^2) \mathrm{e}^{3\omega t} 
(\mathrm{sleafh}_2(B\mathrm{e}^{\omega t}+\phi))^3 \\
&+AB\omega(3\omega+\delta) \mathrm{e}^{2\omega t}
\sqrt{1+(\mathrm{sleafh}_2(B\mathrm{e}^{\omega t}+\phi))^4}=0
\label{A11.4}
\end{split}
\end{align}
For the above equation to hold for an arbitrary variable  $t$, it is necessary to satisfy the following equations.
\begin{equation}
2 B^2 \omega^2+\beta A^2=0 \label{A11.5}
\end{equation}
\begin{equation}
3 \omega+\delta=0 \label{A11.6}
\end{equation}
\begin{equation}
\omega^2+\alpha+\delta \omega=0 \label{A11.7}
\end{equation}
Under the conditions $A\neq0$, $B\neq0$ and $\omega\neq0$, the coefficients of the Duffing equation (\ref{Duffing}) can be obtained as follows.
\begin{equation}
\beta=-2 (\frac{\beta \omega}{A})^2 \label{A11.8}
\end{equation}
\begin{equation}
\delta=-3 \omega \label{A11.9}
\end{equation}
\begin{equation}
\alpha=2\omega^2 \label{A11.10}
\end{equation}

\appendix
\def\thesection{Appendix}
\section{ (X\hspace{-.1em}I\hspace{-.1em}I)}
\label{Appendix12}
In this section, it is shown that the exact solution of the type (X\hspace{-.1em}I\hspace{-.1em}I) satisfies the Duffing equation.
The first derivative of the Eq. (\ref{3.5.1}) is obtained as follows.
\begin{align}
\begin{split}
\frac{\mathrm{d}x(t)}{\mathrm{d}t}
&=A \omega \mathrm{e}^{\omega t} \mathrm{cleafh}_2(B \cdot \mathrm{e}^{\omega t} +\phi) 
+A \mathrm{e}^{\omega t} \cdot \sqrt{(\mathrm{cleafh}_2(B\mathrm{e}^{\omega t}+\phi))^4-1}
 (B \omega \mathrm{e}^{\omega t}) \\
&=A \omega \mathrm{e}^{\omega t} \mathrm{cleafh}_2(B \cdot \mathrm{e}^{\omega t} +\phi) 
+A B \omega \mathrm{e}^{2 \omega t} 
\cdot \sqrt{(\mathrm{cleafh}_2(B\mathrm{e}^{\omega t}+\phi))^4-1}
\label{A12.1}
\end{split}
\end{align}
The second derivative of the above equation is obtained as follows.
\begin{align}
\begin{split}
\frac{\mathrm{d}^2x(t)}{\mathrm{d}t^2}
&=A \omega^2 \mathrm{e}^{\omega t} \mathrm{cleafh}_2(B \cdot \mathrm{e}^{\omega t} +\phi) 
+A \omega \mathrm{e}^{\omega t} \sqrt{(\mathrm{cleafh}_2(B\mathrm{e}^{\omega t}+\phi))^4-1}
(B \omega \mathrm{e}^{\omega t} ) \\
&+2 A B \omega^2 \mathrm{e}^{2 \omega t} \sqrt{(\mathrm{cleafh}_2(B\mathrm{e}^{\omega t}+\phi))^4-1}
+2 A B \omega^2 \mathrm{e}^{2 \omega t} (\mathrm{cleafh}_2(B\mathrm{e}^{\omega t}+\phi))^3
(B\omega \cdot \mathrm{e}^{\omega t}) \\
&=A \omega^2 \mathrm{e}^{\omega t} \mathrm{cleafh}_2(B \cdot \mathrm{e}^{\omega t} +\phi)
+2A B^2 \omega^2 \mathrm{e}^{3 \omega t} (\mathrm{cleafh}_2(B \cdot \mathrm{e}^{\omega t} +\phi))^3 \\
&+3AB \omega^2 \mathrm{e}^{2 \omega t} \sqrt{(\mathrm{cleafh}_2(B\mathrm{e}^{\omega t}+\phi))^4-1}
\label{A12.2}
\end{split}
\end{align}
Substituting these equations into the Duffing equation (\ref{Duffing}) yields the following equation.
\begin{align}
\begin{split}
&A \omega^2 \mathrm{e}^{\omega t} \mathrm{cleafh}_2(B \mathrm{e}^{\omega t}+\phi)
+2AB^2 \omega^2 \mathrm{e}^{3\omega t} (\mathrm{cleafh}_2(B \mathrm{e}^{\omega t}+\phi))^3 \\
&+3AB \omega^2 \mathrm{e}^{2\omega t}\sqrt{(\mathrm{cleafh}_2(B \mathrm{e}^{\omega t}+\phi))^4-1} +\delta(
A\omega \mathrm{e}^{\omega t} \mathrm{cleafh}_2(B \mathrm{e}^{\omega t} +\phi) \\
&+AB\omega \mathrm{e}^{2\omega t} \sqrt{(\mathrm{cleafh}_2(B \mathrm{e}^{\omega t}+\phi))^4-1}) \\
&+\alpha A \mathrm{e}^{\omega t} \mathrm{cleafh}_2(B \mathrm{e}^{\omega t} +\phi)
+\beta( A \mathrm{e}^{\omega t} \mathrm{cleafh}_2(B \mathrm{e}^{\omega t} +\phi) )^3=0 \label{A12.3}
\end{split}
\end{align}
The above equation is summarized as follows.
\begin{align}
\begin{split}
&A (\omega^2+\alpha+\delta \omega) \mathrm{e}^{\omega t} \mathrm{cleafh}_2(B \mathrm{e}^{\omega t}+\phi) + A(2B^2 \omega^2+\beta A^2) \mathrm{e}^{3\omega t} 
(\mathrm{cleafh}_2(B \mathrm{e}^{\omega t}+\phi))^3 \\
&+AB\omega(3\omega+\delta) e^{2\omega t}
\sqrt{(\mathrm{cleafh}_2(B \mathrm{e}^{\omega t}+\phi))^4-1}=0
\label{A12.4}
\end{split}
\end{align}
For the above equation to hold for an arbitrary parameter $t$, it is necessary to satisfy the following equations.
\begin{equation}
2 B^2 \omega^2+\beta A^2=0 \label{A12.5}
\end{equation}
\begin{equation}
3 \omega+\delta=0 \label{A12.6}
\end{equation}
\begin{equation}
\omega^2+\alpha+\delta \omega=0 \label{A12.7}
\end{equation}
Under the conditions $A\neq0$, $B\neq0$ and $\omega\neq0$, the coefficients of the Duffing equation can be derived as follows.
\begin{equation}
\beta=-2 (\frac{B \omega}{A})^2 \label{A12.8}
\end{equation}
\begin{equation}
\delta=-3 \omega \label{A12.9}
\end{equation}
\begin{equation}
\alpha=2\omega^2 \label{A12.10}
\end{equation}

\appendix
\def\thesection{Appendix}
\section{ (X\hspace{-.1em}I\hspace{-.1em}I\hspace{-.1em}I)}
\label{Appendix13}
In this section, it is shown that the exact solution of the type (X\hspace{-.1em}I\hspace{-.1em}I\hspace{-.1em}I) satisfies the Duffing equation. The first derivative of the Eq. (\ref{3.6.1}) is obtained as follows:
\begin{align}
\begin{split}
\frac{\mathrm{d}x(t)}{\mathrm{d}t}
&=A \omega \mathrm{e}^{\omega t} \mathrm{sleaf}_2(B \cdot \mathrm{e}^{\omega t} +\phi) 
+A \mathrm{e}^{\omega t} \cdot \sqrt{1-(\mathrm{sleaf}_2(B\mathrm{e}^{\omega t}+\phi))^4}
 (B \omega \mathrm{e}^{\omega t}) \\
&=A \omega \mathrm{e}^{\omega t} \mathrm{sleaf}_2(B \cdot \mathrm{e}^{\omega t} +\phi) 
+A B \omega \mathrm{e}^{2 \omega t} 
\cdot \sqrt{1-(\mathrm{sleaf}_2(B\mathrm{e}^{\omega t}+\phi))^4}
\label{A13.1}
\end{split}
\end{align}
The second derivative of the above equation is obtained as follows.
\begin{align}
\begin{split}
&\frac{\mathrm{d}^2x(t)}{\mathrm{d}t^2}
=A \omega^2 \mathrm{e}^{\omega t} \mathrm{sleaf}_2(B \cdot \mathrm{e}^{\omega t} +\phi) 
+A \omega \mathrm{e}^{\omega t} \sqrt{1-(\mathrm{sleaf}_2(B\mathrm{e}^{\omega t}+\phi))^4}
(B \omega \mathrm{e}^{\omega t} ) \\
&+2 A B \omega^2 \mathrm{e}^{2 \omega t} \sqrt{1-(\mathrm{sleaf}_2(B\mathrm{e}^{\omega t}+\phi))^4}
+A B \omega \mathrm{e}^{2 \omega t} (-2(\mathrm{sleaf}_2(B\mathrm{e}^{\omega t}+\phi))^3)(B\omega \cdot \mathrm{e}^{\omega t}) \\
&=A \omega^2 \mathrm{e}^{\omega t} \mathrm{sleaf}_2(B \cdot \mathrm{e}^{\omega t} +\phi)
-2A B^2 \omega^2 \mathrm{e}^{3 \omega t} (\mathrm{sleaf}_2(B \cdot \mathrm{e}^{\omega t} +\phi))^3 \\
&+3AB \omega^2 \mathrm{e}^{2 \omega t} \sqrt{1-(\mathrm{sleaf}_2(B\mathrm{e}^{\omega t}+\phi))^4} \label{A13.2}
\end{split}
\end{align}
Substituting these equations into the Duffing equation (\ref{Duffing}) yields the following equation.
\begin{align}
\begin{split}
&A \omega^2 \mathrm{e}^{\omega t} \mathrm{sleaf}_2(B\mathrm{e}^{\omega t}+\phi)
-2AB^2 \omega^2 \mathrm{e}^{3\omega t} (\mathrm{sleaf}_2(B\mathrm{e}^{\omega t}+\phi))^3 \\
&+3AB \omega^2 \mathrm{e}^{2\omega t}\sqrt{1-(\mathrm{sleaf}_2(B\mathrm{e}^{\omega t}+\phi))^4} 
+\delta(
A\omega \mathrm{e}^{\omega t} \mathrm{sleaf}_2(B \mathrm{e}^{\omega t} +\phi) \\
&+AB\omega \mathrm{e}^{2\omega t} \sqrt{1-(\mathrm{sleaf}_2(B\mathrm{e}^{\omega t}+\phi))^4}) \\
& +\alpha A \mathrm{e}^{\omega t} \mathrm{sleaf}_2(B \mathrm{e}^{\omega t} +\phi)
+\beta( A \mathrm{e}^{\omega t} \mathrm{sleaf}_2(B \mathrm{e}^{\omega t} +\phi) )^3=0
\label{A13.3}
\end{split}
\end{align}
The above equation is summarized as follows.
\begin{align}
\begin{split}
&A (\omega^2+\alpha+\delta \omega) \mathrm{e}^{\omega t} \mathrm{sleaf}_2(B \mathrm{e}^{\omega t}+\phi) + A(-2B^2 \omega^2+\beta A^2) \mathrm{e}^{3\omega t} (\mathrm{sleaf}_2(B \mathrm{e}^{\omega t}+\phi))^3 \\
&+AB\omega(3\omega+\delta) \mathrm{e}^{2\omega t}
\sqrt{1-(\mathrm{sleaf}_2(Be^{\omega t}+\phi))^4}=0
\label{A13.4}
\end{split}
\end{align}
For the above equation to hold for an arbitrary variable  $t$, it is necessary to satisfy the following equations.
\begin{equation}
-2 A B^2 \omega^2+\beta A^3=0 \label{A13.5}
\end{equation}
\begin{equation}
AB\omega(3\omega+\delta)=0 \label{A13.6}
\end{equation}
\begin{equation}
A(\omega^2+\alpha+\delta \omega)=0 \label{A13.7}
\end{equation}
Using the conditions $A\neq0$, $B\neq0$ and $\omega\neq0$, the coefficients of the Duffing equation can be derived as follows.
\begin{equation}
\beta=2 (\frac{\beta \omega}{A})^2 \label{A13.8}
\end{equation}
\begin{equation}
\delta=-3 \omega \label{A13.9}
\end{equation}
\begin{equation}
\alpha=2\omega^2 \label{A13.10}
\end{equation}

\appendix
\def\thesection{Appendix}
\section{ (X\hspace{-.1em}I\hspace{-.1em}V)}
\label{Appendix14}
In this section, it is shown that the exact solution of the type (X\hspace{-.1em}I\hspace{-.1em}V) satisfies the Duffing equation.
The first derivative of the Eq. (\ref{3.7.1}) is obtained as follows.
\begin{align}
\begin{split}
\frac{\mathrm{d}x(t)}{\mathrm{d}t}
&=A \omega \mathrm{e}^{\omega t} \mathrm{cleaf}_2(B \cdot \mathrm{e}^{\omega t} +\phi) 
-A \mathrm{e}^{\omega t} \cdot \sqrt{1-(\mathrm{cleaf}_2(B \mathrm{e}^{\omega t}+\phi))^4}
 (B \omega \mathrm{e}^{\omega t}) \\
&=A \omega \mathrm{e}^{\omega t} \mathrm{cleaf}_2(B \cdot \mathrm{e}^{\omega t} +\phi) 
-A B \omega \mathrm{e}^{2 \omega t} 
\cdot \sqrt{1-(\mathrm{cleaf}_2(Be^{\omega t}+\phi))^4}
\label{A14.1}
\end{split}
\end{align}
The second derivative of the above equation is obtained as follows.
\begin{align}
\begin{split}
\frac{\mathrm{d}^2x(t)}{\mathrm{d}t^2}
&=A \omega^2 \mathrm{e}^{\omega t} \mathrm{cleaf}_2(B \cdot \mathrm{e}^{\omega t} +\phi) 
-A \omega \mathrm{e}^{\omega t} \sqrt{1-(\mathrm{cleaf}_2(B \mathrm{e}^{\omega t}+\phi))^4}
(B \omega \mathrm{e}^{\omega t} ) \\
&-2 A B \omega^2 \mathrm{e}^{2 \omega t} \sqrt{1-(\mathrm{cleaf}_2(B \mathrm{e}^{\omega t}+\phi))^4}
-2A B \omega \mathrm{e}^{2 \omega t} ((\mathrm{cleaf}_2(B \mathrm{e}^{\omega t}+\phi))^3)
(B\omega \cdot \mathrm{e}^{\omega t}) \\
&=A \omega^2 \mathrm{e}^{\omega t} \mathrm{cleaf}_2(B \cdot \mathrm{e}^{\omega t} +\phi)
-2A B^2 \omega^2 \mathrm{e}^{3 \omega t} (\mathrm{cleaf}_2(B \cdot \mathrm{e}^{\omega t} +\phi))^3 \\
&-3AB \omega^2 \mathrm{e}^{2 \omega t} \sqrt{1-(\mathrm{cleaf}_2(B \mathrm{e}^{\omega t}+\phi))^4} \label{A14.2}
\end{split}
\end{align}
Substituting these equations into the Duffing equation (\ref{Duffing}) yields the following equation.
\begin{align}
\begin{split}
&A \omega^2 \mathrm{e}^{\omega t} \mathrm{cleaf}_2(B \mathrm{e}^{\omega t}+\phi)
-2AB^2 \omega^2 \mathrm{e}^{3\omega t} ( \mathrm{cleaf}_2(B \mathrm{e}^{\omega t}+\phi))^3
-3AB \omega^2 \mathrm{e}^{2\omega t}\sqrt{1-(\mathrm{cleaf}_2(B \mathrm{e}^{\omega t}+\phi))^4} \\
&+\delta(
A\omega \mathrm{e}^{\omega t} \mathrm{cleaf}_2(B \mathrm{e}^{\omega t} +\phi)
-AB\omega \mathrm{e}^{2\omega t} \sqrt{1-(\mathrm{cleaf}_2(B\mathrm{e}^{\omega t}+\phi))^4}) \\
&+\alpha A \mathrm{e}^{\omega t} \mathrm{cleaf}_2(B \mathrm{e}^{\omega t} +\phi) 
+\beta( A \mathrm{e}^{\omega t} \mathrm{cleaf}_2(B \mathrm{e}^{\omega t} +\phi) )^3=0 \label{A14.3}
\end{split}
\end{align}
The above equation is summarized as follows.
\begin{align}
\begin{split}
&A (\omega^2+\alpha+\delta \omega) \mathrm{e}^{\omega t} \mathrm{cleaf}_2(B \mathrm{e}^{\omega t}+\phi) + A(-2B^2 \omega^2+\beta A^2) \mathrm{e}^{3\omega t} (\mathrm{cleaf}_2(B \mathrm{e}^{\omega t}+\phi))^3 \\
&-AB\omega(3\omega+\delta) \mathrm{e}^{2\omega t}
\sqrt{1-(\mathrm{cleaf}_2(B \mathrm{e}^{\omega t}+\phi))^4}=0
\label{A14.4}
\end{split}
\end{align}
For the above equation to hold for an arbitrary variable $t$, it is necessary to satisfy the following equations.
\begin{equation}
-2 A B^2 \omega^2+\beta A^3=0 \label{A14.5}
\end{equation}
\begin{equation}
-AB\omega(3\omega+\delta)=0 \label{A14.6}
\end{equation}
\begin{equation}
A(\omega^2+\alpha+\delta \omega)=0 \label{A14.7}
\end{equation}
Using the conditions $A\neq0$, $B\neq0$ and $\omega\neq0$, the coefficients of the Duffing equation can be derived as follows.
\begin{equation}
\beta=2 (\frac{\beta \omega}{A})^2 \label{A14.8}
\end{equation}
\begin{equation}
\delta=-3 \omega \label{A14.9}
\end{equation}
\begin{equation}
\alpha=2\omega^2 \label{A14.10}
\end{equation}

\nocite{*}

\bibliographystyle{unsrt}  
\bibliography{references}  

\end{document}